\documentclass[11pt]{amsart}
\usepackage{amssymb, latexsym}
\usepackage{amscd,amsmath,amssymb,amsfonts}
\usepackage{url,enumerate,color,hyperref}
\usepackage{dsfont}
\usepackage[cmtip, all]{xy}

\usepackage[top=30mm,right=30mm,bottom=30mm,left=20mm]{geometry}

\theoremstyle{plain}
\newtheorem{thm}{Theorem}
\newtheorem{lem}[thm]{Lemma}
\newtheorem{cor}[thm]{Corollary}
\newtheorem{prop}[thm]{Proposition}
\theoremstyle{definition}

\newtheorem{rmk}[thm]{Remark}

\numberwithin{thm}{section}
\numberwithin{equation}{thm}

\newcommand{\al}{\alpha}
\newcommand{\om}{\varpi}

\newcommand{\inj}{\hookrightarrow}

\newcommand{\rank}{{\rm rank}}
\newcommand{\End}{{\rm End}}

\newcommand{\Hom}{{\rm Hom}}

\newcommand{\Lie}{{\rm Lie}}

\newcommand{\Char}{{\rm char}}

\newcommand{\Trace}{{\rm Trace}}
\newcommand{\RawTrace}{{\rm RawTrace}}
\newcommand{\Stab}{{\rm Stab}}
\newcommand{\Ker}{{\rm Ker}}

\newcommand{\Aut}{\mathrm{Aut}}
\newcommand{\Out}{\mathrm{Out}}
\newcommand{\Irr}{\mathrm{Irr}}

\newcommand{\eps}{\epsilon}

\newcommand{\Ind}{\mathrm{Ind}}

\newcommand{\Swan}{\mathsf {Swan}}
\newcommand{\Reg}{\mathsf {Reg}}

\def\irr#1{{\rm Irr}(#1)}

\def\cent#1#2{{\bf C}_{#1}(#2)}

\def\nor{\triangleleft\,}

\def\norm#1#2{{\bf N}_{#1}(#2)}

\def\sbs{\subseteq}

\newcommand{\SL}{\mathrm{SL}}
\newcommand{\PSL}{\mathrm{PSL}}
\newcommand{\GL}{\mathrm{GL}}
\newcommand{\PGL}{\mathrm{PGL}}

\newcommand{\SU}{\mathrm{SU}}
\newcommand{\PSU}{\mathrm{PSU}}

\newcommand{\Sp}{\mathrm{Sp}}
\newcommand{\PSp}{\mathrm{PSp}}
\newcommand{\SO}{\mathrm{SO}}

\newcommand{\Sz}{\mathrm{Sz}}

\newcommand{\Norm}{{\rm Norm}}

\newcommand{\sA}{{\mathcal A}}
\newcommand{\sB}{{\mathcal B}}
\newcommand{\sC}{{\mathcal C}}

\newcommand{\sF}{{\mathcal F}}
\newcommand{\sG}{{\mathcal G}}
\newcommand{\sH}{{\mathcal H}}

\newcommand{\sK}{{\mathcal K}}
\newcommand{\sL}{{\mathcal L}}

\newcommand{\sO}{{\mathcal O}}

\newcommand{\A}{{\mathbb A}}

\newcommand{\C}{{\mathbb C}}

\newcommand{\E}{{\mathbb E}}
\newcommand{\F}{{\mathbb F}}
\newcommand{\G}{{\mathbb G}}

\renewcommand{\P}{{\mathbb P}}
\newcommand{\Q}{{\mathbb Q}}
\newcommand{\R}{{\mathbb R}}

\newcommand{\Z}{{\mathbb Z}}

\newcommand{\id}{{\rm id}}

\newcommand{\ppd}{{\rm ppd}}
\newcommand{\ZB}{{\mathbf Z}}
\newcommand{\CB}{{\mathbf C}}

\newcommand{\OB}{{\mathbf O}}

\newcommand{\ABS}{\mathsf{A}}
\newcommand{\SBS}{\mathsf{S}}

\newcommand{\ari}{{\mathrm {arith}}}
\newcommand{\geo}{{\mathrm {geom}}}

\newcommand{\CSP}{\mathrm {({\bf S+})}}
\newcommand{\LP}{{\Lambda^{+}}}
\newcommand{\lam}{\lambda}
\newcommand{\triv}{{\mathds{1}}}
\newcommand{\obar}{\bar{\mathsf o}}
\newcommand{\meo}{\mathsf{meo}}
\newcommand{\dl}{\mathfrak{d}}
\newcommand{\tw}[1]{{}^#1\!}
\renewcommand{\Char}{{\rm Char}}
\newcommand{\Frob}{\mathrm{Frob}}
\newcommand{\FT}{\mathrm{FT}}

\makeatletter
\let\@@pmod\pmod
\DeclareRobustCommand{\pmod}{\@ifstar\@pmods\@@pmod}
\def\@pmods#1{\mkern4mu({\operator@font mod}\mkern 6mu#1)}
\makeatother

\begin{document}
\title{Local systems and Suzuki groups}
\author[Alp\"oge, Katz, Navarro, O'Brien, and Tiep]
{Levent Alp\"oge, Nicholas M. Katz, Gabriel Navarro, E.A. O'Brien, and Pham Huu Tiep}
\address{L. Alp\"oge, Department of Mathematics, Harvard University, Cambridge, MA 02138, U.S.A.}
\email{alpoge@math.harvard.edu}
\address{N.M. Katz, Department of Mathematics, Princeton University, Princeton, NJ 08544,U.S.A.}
\email{nmk@math.princeton.edu}
\address{G. Navarro, Departament de Matem\`atiques, Universitat de Val\`encia, 46100
Burjassot, Val\`encia, Spain}
\email{gabriel@uv.es}
\address{E.A. O'Brien, University of Auckland,  Auckland, New Zealand}
\email{e.obrien@auckland.ac.nz}
\address{P.H. Tiep, Department of Mathematics, Rutgers University, Piscataway, NJ 08854, U.S.A.}
\email{tiep@math.rutgers.edu}
\thanks{The first author gratefully acknowledges the support of the NSF (grant DMS-2002109) and the Society of Fellows.}
\thanks{The third author gratefully acknowledges the support of the Grant PID2019-103854GB-I00 funded by MCIN/AEI\-/10.13039\-/501100011033.}
\thanks{The fourth author gratefully acknowledges  the support of the Marsden Fund of New Zealand via  grant UOA 107.}
\thanks{The fifth author gratefully acknowledges the support of the NSF (grants DMS-1840702 and DMS-2200850), the Simons Foundation, and the Joshua Barlaz Chair in Mathematics.}
\thanks{The authors are grateful to the referee for careful reading and helpful comments on the paper.}

\maketitle

\begin{abstract}
We study geometric monodromy groups $G_{\geo,\sF_q}$ of the local systems $\sF_q$ on the affine line over $\F_2$ 
of rank $D=\sqrt{q}(q-1)$, $q=2^{2n+1}$, constructed in \cite{Ka-ERS}. The main result of the paper shows that $G_{\geo,\sF_q}$ is either the Suzuki simple group
$\tw2 B_2(q)$, or the special linear group $\SL_D$. We also show that $\sF_8$ has geometric monodromy group
$\tw2B_2(8)$, and arithmetic monodromy group $\Aut(\tw2 B_2(8))$ over $\F_2$, thus establishing \cite[Conjecture 2.2]{Ka-ERS}
in full in the case $q=8$.
\end{abstract}

\tableofcontents


\section*{Introduction}
In an earlier paper \cite{Ka-ERS}, one of us, inspired by a paper \cite{Gross} of Gross, defined, for each $n \ge 1$, a local system on $\A^1/\F_2$ of rank  $2^n(2^{2n+1}-1)$, whose geometric monodromy group we conjectured to be the Suzuki group $\Sz(q):=\tw2B_2(q),\ q=2^{2n+1}$, in one of its two lowest dimensional nontrivial irreducible representations. These representations, complex conjugates of each other, are of dimension $2^n(2^{2n+1}-1)$ and have traces in $\Z[i]$, the Gaussian integers.

The definition involved $p$-Witt vectors  of length $2$ for $p=2$, with values in $\F_2$-algebras. We identify $W_2(\F_2)$ with $\Z/4\Z$, by the map
$[a,b] \mapsto a^2+2b$, with the usual convention that we first lift $a,b$ to $\Z$ and then reduce $a^2+2b$ modulo $4$. 
We take the additive
character of $\Z/4\Z$ given by $n \mapsto i^n$, and view it as the additive character $\psi_2:W_2(\F_2) \rightarrow \mu_4(\Z[i])$ given by
$$\psi_2([a,b]):= i^{a^2+2b}.$$

Attached to a Witt vector of length $2$ with coefficients in $\F_2[x]$, say $[a(x),b(x)]$, is the Artin-Schreier-Witt sheaf $\sL_{\psi_2([a(x),b(x)])} $ on $\A^1/\F_2$. It is lisse of rank one, and its trace function is given as follows. For $k/\F_2$ a finite extension, and $x \in k$,
$$\Trace(\Frob_{x,k}|\sL_{\psi_2([a(x),b(x)])})=\psi_2(\Trace_{W_2(k)/W_2(\F_2)}([a(x),b(x)]).$$
If instead we take $a(x), b(x) \in k_0[x]$ for some finite extension $k_0/\F_2$, then the Artin-Schreier-Witt sheaf $\sL_{\psi_2([a(x),b(x)])} $ is lisse of rank one on $\A^1/k_0$, and for $k/k_0$ a finite extension and $x \in k$, the trace of $\Frob_{x,k}$ is given by the same formula (which only makes sense when $k$ is an extension of $k_0$).


The unique nontrivial additive character $\psi$ of $\F_2$ is related to $\psi_2$ by the formula 
$$\psi(b)=\psi_2([0,b]);$$
this is simply the identity $(-1)^b = i^{2b}$.

Under Witt vector addition, $[a,b]=[a,0]+[0,b]$. Thus we have the factorization
$$\sL_{\psi_2([a(x),b(x)])} = \sL_{\psi_2([a(x),0])} \otimes \sL_{\psi(b(x))}.$$

When both $a(x)$ and $b(x)$ are polynomials of degree prime to $p=2$, $\sL_{\psi_2([a(x),0])}$ has Swan conductor 
$\Swan_\infty=p\deg(a(x))$, while $\sL_{\psi(b(x))}$ has $\Swan_\infty=\deg(b(x))$. So in this case, $\sL_{\psi_2([a(x),b(x)])}$ has $\Swan_\infty=\max(p\deg(a(x)), \deg(b(x)))$.

Quite generally, for an Artin-Schreier-Witt sheaf $\sL:=\sL_{\psi_2([a(x),b(x)])} $ on $\A^1/k_0$ with Swan conductor 
$\Swan_\infty =n \ge 2$, its Fourier transform
$\FT_\psi(\sL)$ on $\A^1/k_0$ is an {\it Airy sheaf} in the sense of  \v{S}uch \cite{Such}. It is lisse of rank $n-1$, and all its $\infty$-slopes are
$\frac{n}{n-1}$. It is pure of weight one. Its trace function is given as follows: for $k/k_0$ a finite extension, and $t \in k$, 
$$\Trace(\Frob_{t,k}|\FT_\psi(\sL))=-\sum_{x \in k}\psi_2\bigl(\Trace_{W_2(k)/W_2(\F_2)}([a(x),b(x)+tx])\bigr).$$
Some key facts about Airy sheaves and their monodromy groups are due to \v{S}uch \cite{Such}, and are fundamental for the investigations reported on here.

We now turn to the Suzuki ``candidates" of \cite{Ka-ERS}. Here 
$$n \ge 1\ , q_0:=2^n,\ q:=2q_0^2,\  t(q):=q+1-2q_0.$$ 
We take the Witt vector
$$[x^{t(q)},\sum_{i=1}^n x^{(1+2^i)t(q)}],$$
form the  Artin-Schreier-Witt sheaf 
$$\sL:=\sL_{\psi_2([x^{t(q)},\sum_{i=1}^n x^{(1+2^i)t(q)}])}, $$
form its $\FT_\psi(\sL)$, and twist by the constant field twist 
$(\frac{1}{1-(-1)^n\,i})^{\deg}$, to arrive at the local system $\sF_q$ on $\A^1/\F_2$, whose trace function is given as follows. For $k/\F_2$ a finite extension, and $t \in k$,
$$\Trace(\Frob_{t,k}|\sF_q) =\frac{-1}{\bigl(1-(-1)^n\,i\bigr)^{\deg(k/\F_2)}}\sum_{x \in k}\psi_2\bigl(\Trace_{W_2(k)/W_2(\F_2)}([x^{t(q)},\sum_{i=1}^n x^{(1+2^i)t(q)}+tx])\bigr).$$

A key fact about $\sF_q$ is that the input Witt vector $[x^{t(q)},\sum_{i=1}^n x^{(1+2^i)t(q)}]$ is a function of $x^{t(q)}$. So in the trace formula for $\sF_q$, if we restrict to $t \neq 0$ and make the substitution $x \mapsto x/t$, the formula becomes
$$\Trace(\Frob_{t,k}|\sF_q) =\frac{-1}{\bigl(1-(-1)^n\,i\bigr)^{\deg(k/\F_2)}}\sum_{x \in k}\psi_2\bigl(\Trace_{W_2(k)/W_2(\F_2)}([\frac{x^{t(q)}}{t^{t(q)}},\sum_{i=1}^n \frac{x^{(1+2^i)t(q)}}{t^{(1+2^i)t(q)}}+x])\bigr),$$
which is a function of $t^{t(q)}$.
Thus $\sF_q|\G_m$ has a descent to a lisse sheaf $\sG_q$ on $\G_m/\F_2$ whose trace function is given by
$$\Trace(\Frob_{t,k}|\sF_q) =\frac{-1}{\bigl(1-(-1)^n\,i\bigr)^{\deg(k/\F_2)}}\sum_{x \in k}\psi_2\bigl(\Trace_{W_2(k)/W_2(\F_2)}([\frac{x^{t(q)}}{t},\sum_{i=1}^n \frac{x^{(1+2^i)t(q)}}{t^{(1+2^i)}}+x])\bigr),$$
and such that under the $t(q)$ Kummer pullback, 
$$[t(q)]^\star \sG_q=\sF_q|\G_m.$$

For the case $n=1$, i.e.\ for $\tw2 B_2(8)$, we prove in Theorem \ref{main1} that $\sF_8$ has the predicted geometric monodromy group 
$G_{\geo,\sF_8}=\tw2 B_2(8)$ and arithmetic monodromy group $G_{\ari,\sF_8,\F_2} = \Aut(\tw2 B_2(8))$ over 
$\F_2$, and thus establish \cite[Conjecture 2.2]{Ka-ERS} in full in this case. For each $n \ge 1$, we show that 
$G_{\geo,\sF_q}$ is either $\tw2 B_2(q)$ or $\SL_D$ for $D=\rank(\sF_q)=q_0(q-1)$. A (huge) calculation of fourth moment for $\sF_8$ 
shows that $G_{\geo,\sF_8}$ cannot be $\SL_{14}$. It remains an open problem to prove (or disprove) that each $\sF_q$ has  $G_{\geo,\sF_q}=\tw2 B_2(q)$ when $q \geq 32$.

In fact, we show in Theorem \ref{main2} that we have such a dichotomy of possible geometric monodromy groups $G_\geo$, either $\tw2 B_2(q)$ or $\SL_D,\ D=q_0(q-1)$, for  a general class of local systems of ``the same shape" as $\sF_q$ (see 
Remark \ref{inf-ex} and Theorem \ref{inf-g} for examples of such local systems with $G_\geo=\SL_D$).
Key to our investigation is the fact that these sheaves all satisfy the condition $\CSP$ of  \cite[Definition 1.2]{KT5}. Somewhat to our surprise, the most difficult part of establishing condition $\CSP$ was to show the sheaves in question are geometrically primitive, i.e.\ that the representation of their $G_{\geo}$ is not induced. 
Our proof utilizes the existence of primitive prime divisors of the integer $t(q)$, cf.\ Theorem \ref{prime-suzuki1}. We give a second proof which will be useful in future studies of geometric monodromy groups of quite general Airy sheaves.
Condition $\CSP$ implies that for each of these sheaves, either $G_{\geo}$ is a finite, almost quasisimple group, or
has $G_{\geo}^\circ$ acting irreducibly. This initial dichotomy, plus a substantial group-theoretic analysis, in which Theorem \ref{prim} plays a key role, is what leads to the  
$\tw2 B_2(q)/\SL_D$ dichotomy. Pursuing the study of primitive prime divisors, we prove in Theorems \ref{prime-suzuki1}, \ref{prime-ree}, and 
Corollary \ref{prime-suzuki2} the existence of primitive prime divisors in the orders of maximal tori of the Suzuki-Ree groups 
$\tw2 B_2(q)$, $\tw2 G_2(q)$, and $\tw2 F_4(q)$. We also extend in \S\ref{low-dim} the classification of low-dimensional representations 
of classical groups in characteristic $p \geq 0$, beyond the bounds in \cite{KlL} and \cite{Lu}. These results will be useful in other situations as well. Finally, the structure of arithmetic monodromy groups, assuming finiteness, is determined in Theorem \ref{main3}.

\section{Descents of Suzuki candidates and moment calculations}\label{suz-desc}
For each odd power $q=2^{2n+1}$ of $2$, starting with $q=8$, the finite simple group $\tw2B_2(q)$ has two (complex conjugate) lowest dimensional nontrivial irreducible representations, of dimension $d(q):=q_0(q-1)$, with $q_0:=2^n$. In each case we have a factorization
$$1+d(q) =(q_0+1)t(q)\ \ {\rm with \ }t(q):=q+1-2q_0.$$

In \cite{Ka-ERS}, for each such $q$ there is proposed an Airy sheaf in the sense of  \v{S}uch \cite{Such}, call it $\sF_q$ on $\A^1/\F_2$ which is lisse of rank $d(q)$ and with $\Swan_\infty(\sF_q)=1+d(q)$. Its $I(\infty)$-representation is irreducible of dimension $d(q)$, with all $\infty$-slopes $\frac{1+d(q)}{d(q)}$. Moreover, $\sF_q$ is given \cite[Section 4]{Ka-ERS} with an explicit
descent to a lisse sheaf $\sG_q$ on $\G_m/\F_2$ whose Kummer pullback by $t(q)^{\mathrm {th}}$ power  is (the restriction to $\G_m/\F_2$ of) $\sF_q$:
$$[t(q)]^\star\sG_q \cong \sF_q|\G_m.$$
Thus $\sG_q$ is lisse on $\G_m/\F_2$, tame at $0$ with $I(0)$-representation a direct sum of Kummer characters of order dividing $t(q)$,
and whose $I(\infty)$-representation is irreducible of dimension $d(q)$, with all $\infty$-slopes $\frac{q_0+1}{d(q)}=\frac{q_0+1}{q_0(q-1)}$.

Next we give a slight improvement of \cite[Theorem 6.5]{Ka-RL-T-spor}.
\begin{thm}\label{Mab}Let $\sH_0$ be a lisse sheaf on $\G_m/\F_q$ which is tame at $0$ and pure of weight zero.
Let $a,b$ be nonnegative integers, and consider the moment $M_{a,b}$.
Denote by
$$\sH_0^{a,b}:=\sH_0^{\otimes a}\otimes (\sH_0^\vee)^{\otimes b}.$$
Denote by $A,B,C,D$ the following constants.
$$\begin{aligned}C &:={\rm dimension\ of\ the\ space\ of\ }I(0)\mbox{-\rm{invariants\ in\ }}\sH_0^{a,b},\\
   D&:={\rm dimension\ of\ the\ space\ of\ }I(\infty)\mbox{-\rm{invariants\ in\ }}\sH_0^{a,b},\\
  B & :=\Swan_\infty(\sH_0^{a,b}) +M_{a,b},\\
  A & :=B+M_{a,b}-C-D.\end{aligned}$$
Then we have the following estimate.
$$\biggl{|}\frac{1}{q-1}\sum_{u \in \F_q^\times}\Trace(\Frob_{u,\F_q}|\sH_0^{a,b})\biggr{|} \le \frac{q}{q-1}M_{a,b} + \frac{A\sqrt{q}}{q-1} +\frac{B-A}{q-1}.$$
\end{thm}
\begin{proof}For any lisse sheaf $\sF$ on $\G_m$, the Lefschetz trace formula gives
$$\sum_{u \in \F_q^\times}\Trace(\Frob_{u,\F_q}|\sF) =\Trace\bigl(\Frob_{\F_q}|H^2_c(\G_m/\overline{\F_q},\sF)\bigr)-\Trace\bigl(\Frob_{\F_q}|H^1_c(\G_m/\overline{\F_q},\sF)\bigr).$$
If $\sF$ is pure of weight zero, then $H^2_c$ is pure of weight $2$, and $H^1_c$ is mixed of weight $\le 1$, indeed
$$H^1_c = H^1_c({\rm {wt}}=1) \oplus H^1_c({\rm {wt}} \le 0).$$
Thus for $\sF$ pure of weight zero, 
$$\bigl{|}\sum_{u \in \F_q^\times}\Trace(\Frob_{u,\F_q}|\sF)\bigr{|} \le qh^2_c + \sqrt{q}h^1_c({\rm {wt}}=1) +h^1_c({\rm {wt}} \le 0)  .$$
When $\sF$ is tame at $0$, the Euler-Poincar\'{e} formula gives
$$\Swan_\infty(\sF) = h^1_c - h^2_c.$$

To compute the dimension of $H^1_c({\rm {wt}}=1)$, we use the fact that for the inclusion $j:\G_m \subset \P^1$, the group
$H^1(\P^1/\overline{\F_q},j_\star \sF)$ is pure of weight one. We exploit this by looking at the short exact sequence of sheaves
on $\P^1/\overline{\F_q}$ given by
$$0 \rightarrow j_!\sF \rightarrow j_\star \sF \rightarrow (\sF^{I(0)})_0 \oplus (\sF^{I(\infty)})_\infty \rightarrow 0,$$
where the last two summands are skyscraper sheaves at $0$ and $\infty$. The group 
$$H^2(\P^1/\overline{\F_q},j_\star \sF)=H^2_c(\G_m/\overline{\F_q},\sF)$$
is the Tate-twisted group of $\pi_1^{\geo}$ co-invariants in $\sF$, but as $\sF$ is pure, the action of $\pi_1^{\geo}$ is semisimple, so this is also the (Tate-twisted) group of $\pi_1^{\geo}$ invariants.
The group
$$H^0(\P^1/\overline{\F_q},j_\star \sF)=H^0(\G_m/\overline{\F_q},\sF)$$
is the space of $\pi_1^{\geo}$ invariants in $\sF$, so has dimension $h^0=h^2_c$. Consider the long exact sequence
$$0 \rightarrow H^0(\P^1/\overline{\F_q},j_\star \sF) \rightarrow \sF^{I(0)} \oplus \sF^{I(\infty)} \rightarrow H^1_c(\G_m/\overline{\F_q},\sF) \rightarrow H^1(\P^1/\overline{\F_q},j_\star\sF_0) \rightarrow 0.$$
Apply it with $\sF$ taken to be $\sH_0^{a,b}$. Then $h^ 0=h_c^2$ is $M_{a,b}$, and the Euler-Poincar\'{e} formula gives
$$h_c^1 =\Swan_\infty(\sH_0^{a,b}) + M_{a,b}.$$
Thus we have the equalities
$$h^1_c({\rm {wt}}=1)=h^1_c +h^0 -\dim  \sF^{I(0)} - \dim  \sF^{I(\infty)}=A,$$
and
$$h^1_c({\rm {wt}} \le 0)=B-A.$$
Then the estimate
$$|\sum_{u \in \F_q^\times}\Trace(\Frob_{u,\F_q}|\sF) |  \le qh^2_c + \sqrt{q}h^1_c({\rm {wt}}=1) +h^1_c({\rm {wt}} \le 0)$$
becomes
$$\bigl{|}\sum_{u \in \F_q^\times}\Trace(\Frob_{u,\F_q}|\sH_0^{a,b})\bigr{|} \le q M_{a,b}  +A\sqrt{q} + (B-A).$$
\end{proof}

We will now apply this estimate to the descent $\sG_8$ to $\G_m/\F_2$ of the lisse sheaf $\sF_8$ on $\A^1/\F_2$. Recall from \cite[\S2]{Ka-ERS} that $\sF_8$ is the Fourier transform of the Artin-Schreier-Witt sheaf $\sL_{\psi_2([x^{5},x^{15}])}$, with a constant field twist by $\frac{1}{1+i}$.
Thus for $k/\F_2$ a finite extension, and $t \in k$,
$$\Trace(\Frob_{t,k}|\sF_8)=-(\frac{1}{1+i})^{\deg(k/\F_2)}\sum_{x \in k}\psi_2(\Trace_{W_2(k)/W_2(\F_2)}([x^5,x^{15}+tx])).$$
The descent $\sG_8$ is the lisse sheaf on $\G_m/\F_2$ whose trace function is given as follows. For  $k/\F_2$ a finite extension, and $t \in k^\times$,
$$\Trace(\Frob_{t,k}|\sG_8)=-(\frac{1}{1+i})^{\deg(k/\F_2)}\sum_{x \in k}\psi_2(\Trace_{W_2(k)/W_2(\F_2)}([x^5/t,x^{15}/t^3+x])).$$
Thus one visibly has, for $t \in k^\times$, the identity
$$\Trace(\Frob_{t,k}|\sF_8)=\Trace(\Frob_{t^5,k}|\sG_8),$$
simply by the substitution $x \mapsto x/t$ in the formula for $\sF_8$.

In any extension field $k/\F_2$ such that $\gcd(5, \#k^\times)=1$, the map $t \mapsto t^5$ is bijective on $k^\times$. Such $k/\F_2$ are precisely those whose degree over $\F_2$ is not divisible by $4$. For such a $k/\F_2$,  the traces $\Trace(\Frob_{t,k}|\sG_8)$ as $t$ runs over $k^\times$ are precisely the traces $\Trace(\Frob_{t,k}|\sF_8)$ as $t$ runs over $k^\times$. An extensive 
calculation shows that over $\F_{2^{18}}^\times$, the seven traces which occur, with their multiplicities, are
\begin{itemize}
\item $-2i$, {\rm multiplicity\ } 16256,
\item $-2$, {\rm multiplicity\ } 4095,
\item $-1$, {\rm multiplicity\ } 52429,
\item $0$, {\rm multiplicity\ } 112347,
\item $1$, {\rm multiplicity\ } 60495,
\item $2i$, {\rm multiplicity\ } 16512,
\item $14$, {\rm multiplicity\ } 9.
\end{itemize}
[We have not been able to find any conceptual explanation for the multiplicities of these traces, let alone for the fact that all traces are algebraic integers. Had we known that the traces of all Frobenii over all finite extensions of $\F_2$ are algebraic integers, we would have been able to conclude that $\sG_8$ has finite geometric monodromy group $G_\geo$, and the proof of the main result Theorem \ref{main1} would have been much simpler.]
 
This computation was carried out using {\sc Magma} 2.26-6 \cite{Magma} at
the University of Auckland. It exploited a new feature of {\sc Magma}
which supports running tasks in parallel on multiple processors:
each task performed the necessary computation
for a given $t \in k^\times$.
Using 55 3.0GHz processors, the computation
completed in 5.75 days, taking about 7500 hours of CPU time.

The empirical $M_{2,2}$ for $\sG_8$ over $\F_{2^{18}}^\times$ is thus  approximately 3.99963378766551080898593515753.
Applying Theorem \ref{Mab}, we find
\begin{cor}\label{m4}
For the lisse sheaf $\sG_8$ on $\G_m/\F_2$, $M_{2,2} >2$.
\end{cor}
\begin{proof}The sheaf $\sF_8$ is a geometrically irreducible Airy sheaf, lisse on $\A^1/\F_2$ of rank $14$, whose $I(\infty)$ representation is
irreducible, with all slopes $15/14$. Its descent $\sG_8$ is thus lisse and geometrically irreducible on $\G_m$, its $I(0)$ representation is the direct sum with multiplicities of the characters of order dividing $5$, and its $I(\infty)$ representation is irreducible, with all slopes $3/14$.

Let us denote by 
$$\sK :=(\sG_8)^{\otimes 2}\otimes (\sG_8^\vee)^{\otimes 2}.$$
In the proof of Theorem \ref{Mab}, the estimate
$$\bigl{|}\sum_{t \in \F_q^\times}\Trace(\Frob_{t,\F_q}|\sK) \bigr{|}  \le qh^2_c + \sqrt{q}h^1_c({\rm {wt}}=1) +h^1_c{(\rm {wt}} \le 0)$$
can first be weakened to 
$$\bigl{|}\sum_{x \in \F_q^\times}\Trace(\Frob_{u,\F_q}|\sK) \bigr{|}  \le qh^2_c + \sqrt{q}h^1_c({\rm {wt}}=1) +h^1_c.$$
The equality  
$$h^1_c({\rm {wt}}=1)=\Swan_{\infty}(\sK) +2M_{2,2}-\dim \sK^{I(0)} -\dim \sK^{I(\infty)} $$
can be weakened to 
$$h^1_c({\rm {wt}}=1) \le \Swan_\infty(\sK) +M_{2,2}-\dim \sK^{I(0)} ,$$
simply because in the exact sequence, the space $H^0$ of global invariants of $\sK$ injects into the space $\sK^{I(\infty)}$ of $I(\infty)$-invariants.

Thus $\sK$ is lisse on $\G_m$, its $I(0)$ representation is the direct sum with multiplicities of the characters of order dividing $5$, and all its
 $I(\infty)$ slopes are $\le 3/14.$ So we have the crude estimate
 $$\Swan_\infty(\sK) \le \rank(\sK)({\rm biggest\ slope}) \le 14^4(3/14) =8232.$$
 The $I(0)$ representation of $\sK$ is $\End(\End(\mbox{the }I(0)\mbox{-representation\ of\ }\sG_8)).$
 
 We now turn to the  $I(0)$ representation of $\sG_8$. The action of $I(0)$ is through $\mu_5$.  So in terms of a fixed character $\chi$ of order $5$, it is a direct sum
 $$a\triv +b\chi +c\chi^2 +d\chi^3 +e\chi^4,$$
 with non-negative integers $a,b,c,d,e$ which sum to $14$. We also know from Deligne's ``independence of $\ell$" result \cite[Theorem 9.8]{De-Const} or from Serre-Tate \cite[Theorem 2(ii)]{Se-Ta},  that the character of this representation of $I(0)$ has values in the field $\Q(i)$ (because $\sG_8$ is part of a compatible system over $\Q(i)$). On the other hand, this character has values in $\Q(\zeta_5)$. But the intersection $\Q(i)\cap \Q(\zeta_5)$ is just $\Q$, so the trace has values in $\Q$. In other words, the quantity 
 $$a +b\zeta_5 +c\zeta_5^2 +d\zeta_5^3 +e\zeta_5^4$$
 lies in $\Q$. This in turn forces $b=c=d=e$, and so  $a +b\zeta_5 +b\zeta_5^2 +b\zeta_5^3 +b\zeta_5^4=a-b$. But $a+4b=14$, so there are only  four possibilities for the character, namely
 $$2\triv +3\chi +3\chi^2 +3\chi^3 +3\chi^4,\ ,6\triv +2\chi +2\chi^2 +2\chi^3 +2\chi^4, 10\triv +1\chi +1\chi^2 +1\chi^3 +1\chi^4, \ 14\triv.$$
 More intrinsically, let us denote by $\Reg$ the regular representation of $\mu_5$:
 $$\Reg =\triv +\chi +\chi^2 +\chi^3 +\chi^4.$$
 Then  the  $I(0)$ representation of $\sG_8$ is one of
 $$-\triv +3\Reg, \ 4\triv+2\Reg, \ 9\triv+\Reg, \ 14\triv.$$
 Because the character takes real (in fact integer) values, 
 $\dim \sK^{I(0)}$ is the coefficient of $\triv$ in the fourth power of the character. For each of our four candidates, this is easily computed by hand, because $\Reg^2=5\Reg$, $\Reg^3=5^2\Reg$, $\Reg^4=5^3\Reg$. The least of the possible values of  $\dim \sK^{I(0)}$ is $7684$, the value attained by the candidate $-\triv +3\Reg$.
 So the weakened estimate becomes
 $$\bigl{|}\sum_{t \in \F_{2^{18}}^\times}\Trace(\Frob_{t,\F_{2^{18}}}|\sK)\bigr{|} \le 2^{18}M_{2,2} + (8232 +M_{2,2} -7684)2^9 +(8232 +M_{2,2}).$$
 Dividing by $2^{18}-1$, we get
 $$3.999 \le \frac{2^{18}}{2^{18}-1}M_{2,2} + \frac{(548+M_{2,2})2^9}{2^{18}-1} + \frac{8232 +M_{2,2}}{2^{18}-1}.$$
 Here $2^{18}-1=262143$, and $2^9/(2^{18}-1) \le (2^9+1)/(2^{18}-1)=1/(2^9-1) =1/511$. 
 
 Thus if $M_{2,2}$ were $2$, we would get
 $$3.999 \le (1 + (1/262143))2 +550/511 + 8332/262143,$$
 which is nonsense. 
\end{proof}

\section{Background results on determinants, rationality, and slopes}

\begin{thm}\label{dets}Let $p$ be a prime, $k/\F_p$ a finite extension, $\nu \ge 1$ an integer, and $\sF$ a lisse $\overline{\Q_\ell}$ sheaf on a smooth, geometrically connected scheme $U/k$ which is pure of weight zero and part of a compatible system of lisse sheaves on $U$ whose trace functions take values in
$\Q(\zeta_{p^\nu})$. Denote by $G_{\geo}\leq G_{\ari}$ the geometric and arithmetic monodromy groups of $\sF$. Then we have the following results.
\begin{enumerate}[\rm(i)]
\item $\det(G_{\ari})$ has finite order dividing $2p^\nu$ if $p$ is odd, dividing $2^\nu$ if $p=2$.
\item $\det(G_{\geo})$ has finite order dividing $2p^\nu$ if $p$ is odd, dividing $2^\nu$ if $p=2$.
\item Suppose that $U$ is a dense open set of $\P^1$, and that at each point $x \in \P^1(\overline{k})\setminus U(\overline{k})$, all $I(x)$-slopes of $\sF$ are $<1$. Then $\det(G_{\geo})$ has order dividing $2$, and is trivial if $p=2$.
\end{enumerate}
\end{thm}

\begin{proof}
(i) For any finite extension $L/k$ and any point $t \in U(L)$, $\det(\Frob_{t,L}|\sF)$ lies in $\Q(\zeta_{p^\nu})$. It has absolute value $1$ at every archimedean place of $\Q(\zeta_{p^\nu})$ (purity of weight zero), and is a $\lambda$-adic unit at every finite place $\lambda$ of residue characteristic $\neq p$ (being part of a compatible system).  Because $\Q(\zeta_{p^\nu})$ has a unique place above $p$, the product formula tells us that this determinant is a unit
at all finite places. Thus it is a root of unity in $\Q(\zeta_{p^\nu})$, so of order dividing $2p^\nu$ for $p$ odd, and of order dividing $2^\nu$ if $p=2$.

As $G_{\geo}$ is a subgroup of $G_{\ari}$, we trivially obtain (ii). 

To prove (iii), the slope hypotheses imply that $\det$ as a character of $G_{\geo}$ has slope $<1$, hence $0$, at each missing point, and thus is everywhere tame, and hence (being on a dense open set of $\P^1$) has order prime to $p$.
\end{proof}

\begin{thm}\label{ssandI(x)}
Let $p$ be a prime, $k/\F_p$ a finite extension, $\nu \ge 1$ an integer, and $\sF$ a lisse $\overline{\Q_\ell}$ sheaf on a smooth, geometrically connected scheme $U/k$ which is pure of weight zero and part of a compatible system of lisse sheaves on $U$ whose trace functions take values in a number field $E$. Suppose further that there exists a proper smooth curve with geometrically connected fibres
$$\pi: \sC \rightarrow U,$$
a finite group $\Gamma$ of automorphisms of $\sC/U$, and a linear character 
$$\chi: \Gamma \rightarrow E^\times$$
such $\sF$ is isomorphic to the $\chi$-component of $R^1\pi_\star(\overline{\Q_\ell})$. Then we have the following results.
\begin{enumerate}[\rm(i)]
\item For any finite extension $L/k$ and any point $t \in U(L)$, the action of $\Frob_{t,L}|\sF$ is semisimple.
\item If $U$ is a curve, with complete nonsingular model $X$, then for each point $x \in X(\overline{k})\setminus U(\overline{k})$, the character
of the action of the inertia group $I(x)$ acting on $\sF$ has values in $K$.
\end{enumerate}
\end{thm}

\begin{proof}The semisimplicity of Frobenii on $H^1$ of curves goes back to Weil. By Serre-Tate \cite[Theorem 2(ii)]{Se-Ta}, the character of $I(x)$ on $R^1\pi_\star(\overline{\Q_\ell})$ has values in $\Z$. When we project the $I(x)$ action onto the $\chi$ component, the character of the resulting $I(x)$ action has values in $\Q(\chi)$,
a subfield of $K$.
\end{proof}

\begin{rmk} 
Theorem \ref{ssandI(x)} applies to the Airy sheaves of \v{S}uch \cite{Such} and any of their descents, where the family of curves in question is
a family of Artin-Schreier-Witt coverings of $\A^1$ (compactified by adding its one point at $\infty$).
\end{rmk}

Next we record a general result on $G_\geo$ of lisse sheaves on open sets of $\A^1$:

\begin{thm}\label{slope}
Let $U$ be a dense open set of $\A^1/\overline{\F_p}$, and $\sF$ a lisse $\overline{\Q_\ell}$-sheaf on $U$, with $\ell \neq p$.
Suppose all the $\infty$-slopes of 
$\sF$ are at most $\sigma$, for some $0 < \sigma < 1$. Suppose the geometric monodromy group $G=G_\geo$ of $\sF$ admits a representation
$\Phi:G \to \GL_d(\F)$ over some algebraically closed field $\F$ of characteristic $\neq p$, of dimension $d < 1/\sigma$. Then 
$\Phi$ is tame at $\infty$.
\end{thm}

\begin{proof}
The hypothesis that all $\infty$-slopes of $\sF$ are $\le \sigma$ is that for all $y > \sigma$, the upper numbering subgroup $I(\infty)^y$ acts trivially on $\sF$, i.e.\ dies in $G_{\geo}$, and hence dies under $\Phi$. Thus all $\infty$-slopes of $\sF$ are $\le \sigma$, and hence
$\Swan_\infty(\Phi) \le d\sigma <1$. As
Swan conductors are non-negative integers, 
$\Swan_\infty(\Phi) = 0$. This means $P(\infty)$ acts trivially in $\Phi$, 
i.e.\ $\Phi$ is tame at $\infty$.
\end{proof}

Now we prove a generalization of  \cite[Theorem 4.16]{KT5}:

\begin{thm}\label{bound}
Let $U$ be a dense open set of $\A^1/\overline{\F_p}$, $\sF$ a lisse $\overline{\Q_\ell}$-sheaf on $U$, with $\ell \neq p$,
and let $G$ be the geometric monodromy group
of $\sH$. Suppose that the following hold:
\begin{enumerate}[\rm(a)]
\item All $\infty$-slopes of $\sH$ are at most $\sigma$ for some $0 < \sigma < 1$, and $\sH$ is not tame at $\infty$; 
\item $G$ is a finite almost quasisimple group: $S \lhd G/\ZB(G) \leq \Aut(S)$ for some finite non-abelian simple group $S$;
\item For some normal subgroup $R$ of $G/\ZB(G)$ containing $S$, $R$ admits either a faithful $d$-dimensional linear representation 
$$\Phi:R \to \GL_d(\F),$$ 
or an $e$-dimensional projective representation 
$$\Psi:R \to \PGL_e(\F)$$
which is nontrivial over $S$, 
over some algebraically closed field $\F$ of characteristic $\neq p$.
\end{enumerate}
Then 
$$1/\sigma \leq d\cdot [G/\ZB(G):R] \leq d\cdot|\Out(S)|,$$
respectively
$$1/\sigma \leq (e^2-1) \cdot [G/\ZB(G):R] \leq (e^2-1)\cdot|\Out(S)|.$$
\end{thm}

\begin{proof}
Note that the given $\Psi$ is faithful.
Indeed, $\Ker(\Psi) \lhd R$ does not contain $S$, so it
intersects $S$ trivially by simplicity of $S$. Because both $S$ and  $\Ker(\Psi)$ are normal in $R$, the commutator $[S,  \Ker(\Psi) ]  \subset S \cap  \Ker(\Psi) =1$. Thus $\Ker(\Psi) \leq \CB_R(S) \le \CB_{\Aut(S)}(S)=1$. Hence
$R$ is embedded in $\PGL(U)$, where $U = \F^e$. Composing this embedding with the faithful action of $\PGL(U)$ on 
$\End^0(U) = \End(U)/{\rm scalars}$, we obtain a faithful action of $R$ on a module of dimension $\leq e^2-1$. 
Thus it suffices to prove the bound 
$1/\sigma \leq d\cdot [G/\ZB(G):R]$ 
when $\Phi: R \to \GL(V)$ is given.

So assume the contrary: $\Phi: R \to \GL(V)$ is faithful for some $V$ with $\dim(V)=d$, but
\begin{equation}\label{bd21}
 1/\sigma > d \cdot [G/\ZB(G):R].
\end{equation} 
Let $\tilde V$ denote the $\bar{G}$-module $\Ind^{\bar{G}}_R(V)$ for $\bar{G}:=G/\ZB(G)$. 
Note that $\bar{G}$ acts faithfully on $\tilde V$. Indeed, let $K \lhd \bar{G}$ denote the kernel of the action of $\bar{G}$ on $\tilde V$. 
By the construction of $V$ as the induced representation, the $R$-module $\tilde V$ contains $V$ as a submodule. But $S$ acts faithfully on $V$, hence
$S \cap K=1$. As $S \lhd \bar{G}$, it follows that $[S,K]=1$, and so 
$$K \leq \CB_{\bar{G}}(S) \leq \CB_{\Aut(S)}(S)=1.$$  
We also note that
$$\dim(\tilde V) = [\bar{G}:R]\cdot \dim(V) =  d \cdot [\bar{G}:R] < 1/\sigma$$
by \eqref{bd21}.

Now view $\tilde V$ as a representation of $G$, of dimension $< 1/\sigma$. By Theorem \ref{slope}, this representation is tame at $\infty$. Thus the image $Q$ in $G$ of  $P(\infty)$ acts trivially on $\tilde V$. But $G/\ZB(G)$ acts faithfully on $\tilde V$. Therefore $Q$ lands in $\ZB(G)$. Recall that $I(\infty)$ has finite image $J$ in $G$, and $J/Q$ is cyclic. As $Q \leq \ZB(J)$, it follows that 
$J$ is abelian. Thus all simple $J$-summands in $\sH$ are one-dimensional, and at least one of them is wild, as $\sH$ is not tame at
$\infty$. Each one-dimensional wild component has $\Swan$ a strictly positive integer, which is also its slope, contradicting the hypothesis that
all $\infty$-slopes of $\sH$ are $<1$. 
\end{proof}

\section{Primitive prime divisors for Suzuki-Ree groups}
The order of a finite group of Lie type $G(\F_q)$ over a field $\F_q$ is usually a product of a power of $q=p^f$ 
($p$ the defining characteristic) and the values at 
$q$ of cyclotomic polynomial $\Phi_m(q)$ for various $m$. In a number of problems on $G(\F_q)$, the existence of primitive prime divisors 
$\ppd(q,m)$ or $\ppd(p,mf)$ for certain $m$ was
helpful. Recall \cite{Zs} that for $a,m \in \Z_{\geq 2}$, a {\it primitive 
prime divisor $\ell=\ppd(a,m)$} is a prime divisor of $a^m-1$ that does not divide $\prod^{m-1}_{i=1}(a^i-1)$; such a prime divisor always exists unless $(a,m) = (2,6)$ or $m=2$ and $a+1$ is a $2$-power. For the Suzuki-Ree groups $\tw2 B_2(q)$ with $q=2^n$, 
$2 \nmid n \geq 3$, $\tw2 G_2(q)$ with $q=3^n$, $2 \nmid n \geq 3$, and $\tw2 F_4(q)$ with $q=2^n$, $2 \nmid n \geq 3$, 
some factor $\Phi_m(q)$ of $|G(\F_q)|$ decomposes further into values at $\sqrt{q}$ of polynomials over $\Z[\sqrt{2}]$ or 
$\Z[\sqrt{3}]$. More precisely,
$$\Phi_4(q) = q^2+1 = (q-\sqrt{2q}+1)(q+\sqrt{2q}+1)$$
for $\tw2 B_2(q)$, 
$$\Phi_6(q) = q^2-q+1 = (q-\sqrt{3q}+1)(q+\sqrt{3q}+1)$$
for $\tw2 G_2(q)$, and 
$$\Phi_{12}(q) = q^4-q^2+1 = (q^2-q\sqrt{2q}+q-\sqrt{2q}+1)(q^2+q\sqrt{2q}+q+\sqrt{2q}+1)$$
for $\tw2 F_4(q)$. In applications, it is desirable to prove that these factors also possess primitive prime divisors
$\ppd(p,4n)$, respectively $\ppd(p,6n)$, $\ppd(p,12n)$, whose existence does not follow from \cite{Zs}. The main
results of this section establish the existence of such prime divisors for Suzuki-Ree groups.  

\subsection*{3A. Almost equidistribution of coprime integers in congruence classes} 
For $n \in \Z_{\geq 1}$, let $\phi(n)$ denote the Euler function of $n$, let $\mu(n)$ denote the M\"obius function of 
$n$, and let $\omega(n)$ denote the number of distinct prime divisors of $n$ (not counting multiplicities). First we prove the following

\begin{prop}\label{equi1}
Let $m,n \in \Z_{\geq 1}$ be coprime integers. For any integer $0 \leq a \leq m-1$, the number $N_a$ of integers 
$1 \leq k \leq n$ such that $\gcd(k,n)=1$ and $k \equiv a \pmod*{m}$ satisfies
$$\bigl{|}N_a - \frac{\phi(n)}{m}\bigr{|} < 2^{\omega(n)}.$$ 
\end{prop}

\begin{proof}
For $k \in \Z_{\geq 1}$, define $F(k) := 0$ if $\gcd(k,n) > 1$ and $F(k):=1$ if $\gcd(k,n) = 1$. By \cite[Theorem 4.7]{MNZ}, 
$F(k) = \sum_{d|\gcd(k,n)}\mu(d)$. Now 
\begin{equation}\label{equi10}
 N_a  = \sum_{1 \leq k \leq n,~k \equiv a \pmod*{m}}F(k) 
= \sum_{1 \leq k \leq n,~k \equiv a \pmod*{m}}\sum_{d|\gcd(k,n)}\mu(d)
 = \sum_{d|n}\mu(d)N(a,d), 
\end{equation}
with 
$$N(a,d):=\sum_{1 \leq k \leq n,~k \equiv a \pmod*{m},~d|k}1.$$ 
If $d|n$, then $\gcd(d,m)=1$, so we can find $1 \leq e \leq m-1$ such that $de \equiv 1 \pmod*{m}$. Now 
write every $1 \leq k \leq n$ with $d|k$ as $k = dl$ with $1 \leq l \leq n/d$. Then the condition that $k \equiv a \pmod*{m}$ is 
equivalent to $l \equiv ea \pmod*{m}$, in which case we can write $l = s+mi$ with $0 \leq r \leq m-1$, $ea \equiv r \pmod*{m}$,
and $i \in \Z$. To count the number $N(a,d)$ of  $i$ occurring, 
write $n/d=qm+r$ with $0 \leq r \leq m-1$ and $q \in \Z_{\geq 0}$.
Certainly, every $0 \leq i \leq q-1$ works, but neither $i=-1$ nor $i=q+1$ can occur. It follows that
\begin{equation}\label{equi11}
  n/md-1 < q \leq N(a,d) \leq q+1 < n/md+1.
\end{equation}  
We also note by \cite[(4.1)]{MNZ} that $\phi(n) = \sum_{d|n}\mu(d)n/d$ and that $\sum_{d|n}|\mu(d)|$ is the number of square-free 
divisors of $n$ and hence equals to $2^{\omega(n)}$. Combining with \eqref{equi10} and \eqref{equi11}, this yields
$$\bigl{|}N_a - \frac{\phi(n)}{m}\bigr{|} = \bigl{|}\sum_{d|n}\mu(d)(N(a,d)-n/md) \bigr{|} < \sum_{d|n}|\mu(d)| = 2^{\omega(n)}.$$
\end{proof}

We will also need the following analogue of Proposition \ref{equi1}:

\begin{prop}\label{equi2}
Let $n \in \Z_{\geq 1}$ be an odd integer divisible by $3$. For any integer $0 \leq a \leq 11$ coprime to $3$, the number $N_a$ of integers 
$1 \leq k \leq n$ such that $\gcd(k,n)=1$ and $k \equiv a \pmod*{12}$ satisfies
$$\bigl{|}N_a - \frac{\phi(n)}{8}\bigr{|} < 2^{\omega(n)-1}.$$ 
\end{prop}

\begin{proof}
As in the proof of Proposition \ref{equi2}, we have
$$N_a  = \sum_{1 \leq k \leq n,~k \equiv a \pmod*{12}}F(k) 
= \sum_{1 \leq k \leq n,~k \equiv a \pmod*{12}}\sum_{d|\gcd(k,n)}\mu(d)
= \sum_{d|n}\mu(d)N(a,d),$$
with 
$$N(a,d):=\sum_{1 \leq k \leq n,~k \equiv a \pmod*{12},~d|k}1.$$
Now, if $d|n$ but $3|d$, then $N(a,d) = 0$ since $3 \nmid a$. Hence,  
\begin{equation}\label{equi20}
  N_a = \sum_{d|n,~3 \nmid d}\mu(d)N(a,d).
\end{equation} 
Next, $\mu(d)=0$ if $9|d$, and $\mu(d)=-\mu(d/3)$ if $3|d$. It follows from \cite[(4.1)]{MNZ} that 
\begin{equation}\label{equi21}
  \frac{\phi(n)}{n} = \sum_{d|n}\frac{\mu(d)}{d}  = \sum_{d|n,~3 \nmid d}\frac{\mu(d)}{d} - \sum_{3d'=d|n,3 \nmid d'}\frac{\mu(d')}{3d'} 
  = \frac{2}{3}\sum_{d|n,~3 \nmid d}\frac{\mu(d)}{d}.
\end{equation}   
If $d|n$ and $3 \nmid d$, then $\gcd(d,12)=1$, so we can find $1 \leq e \leq 11$ such that $de \equiv 1 \pmod*{12}$. 
The proof of \eqref{equi11} repeated verbatim shows that
$$n/12d-1 < q \leq N(a,d) \leq q+1 < n/12d+1.$$
Combining with \eqref{equi20} and \eqref{equi21}, this yields
$$\bigl{|}N_a - \frac{\phi(n)}{8}\bigr{|} = \bigl{|}\sum_{d|n,~3 \nmid d}\mu(d)(N(a,d)-n/12d) \bigr{|} < \sum_{d|n,~3 \nmid d}|\mu(d)| = 2^{\omega(n)-1}.$$
\end{proof}

\subsection*{3B. Primitive prime divisor for Suzuki groups}
We make the choice $\sqrt{2} > 0$. For odd $n \in \Z_{\geq 1}$ and $a=1,3$, set
$$P_{2,a}(n):= \prod_{1 \leq k < 8n,~\gcd(k,n)=1,~k \equiv a \pmod*{8}}\bigl(3 - 2\sqrt{2}\cos\frac{k\pi}{4n}\bigr).$$

\begin{prop}\label{prod-suzuki}
If $2 \nmid n \geq 2603$ and $a=1,3$, then $P_{2,a}(n) > 2n$. 
\end{prop} 

\begin{proof}
(i) First we note that
\begin{equation}\label{prod10}
  \phi(n) \geq \max\bigl(2^{2.2\omega(n)},n^{6/7}\bigr)
\end{equation}
when $2 \nmid n$ and $n \neq 1,3,9,15,21,33,45,75,105,165,195$. Indeed, suppose $s:=\omega(n) \geq 1$ and write 
$n=\prod^s_{i=1}p_i^{a_i}$ for some prime divisors $2 < p_1 < p_2 < \ldots < p_s$ of $n$. If $p_1 \geq 5$, then
$$\frac{n}{\phi(n)} = \prod^s_{i=1}\frac{p_i}{p_i-1} \leq (5/4)^s < 5^{s/7} \leq n^{1/7},$$
and so $\phi(n) > n^{6/7}$. If $p_1 \geq 7$, then 
$$\phi(n) = \prod^s_{i=1}p_i^{a_i}(1-1/p_i) \geq 6^s > 2^{2.5s}.$$ 
If $p_1=5$ and $s \geq 2$, then
$$\phi(n) = \prod^s_{i=1}p_i^{a_i}(1-1/p_i) \geq 4 \cdot 6^{s-1} > 2^{2.25s}.$$
If $p_1=5$ and $s=1$, but $n \neq 5$ then $\phi(n) \geq 20 > 2^{4s}$. 

In the rest of the proof of \eqref{prod10}, we may assume that $p_1=3$. First suppose that $s \geq 4$. As
$(11^6/10^7)^{s-3} <  (2 \cdot 4 \cdot 6)^7/(3 \cdot 5 \cdot 7)^6$, 
$$\frac{n}{\phi(n)} = \prod^s_{i=1}\frac{p_i}{p_i-1} \leq \frac{3}{2} \cdot \frac{5}{4} \cdot \frac{7}{6} \cdot \bigl(\frac{11}{10}\bigr)^{s-3} < 
   \bigl( 3 \cdot 5 \cdot 7 \cdot 11^{s-3} \bigr)^{1/7} < n^{1/7}.$$
Also, $\phi(n) \geq 2 \cdot 4 \cdot 6 \cdot 10^{s-3} > 2^{2.23s}$.    

Next suppose that $s = 3$. If $p_3 \geq 11$, then 
$$\frac{n}{\phi(n)} = \prod^s_{i=1}\frac{p_i}{p_i-1} \leq \frac{3}{2} \cdot \frac{5}{4} \cdot \frac{11}{10} < 
   \bigl( 3 \cdot 5 \cdot 11 \bigr)^{1/7} < n^{1/7}.$$ 
If $p_3 < 11$ then $(p_1,p_2,p_3)=(3,5,7)$; thus if $n \neq 105$ then $n > 3 \cdot 105$ and so
$$n/\phi(n) = (3 \cdot 5 \cdot 7)/(2 \cdot 4 \cdot 6) < n^{1/7}.$$    
Also, if $n \neq 105, 165,195$, then $\phi(n) \geq 2 \cdot 6 \cdot 10 > 2^{2.2s}$.  

Suppose now that $s = 2$. If $p_2 \geq 13$, then 
$n/\phi(n) \leq (3/2) \cdot (13/12) < (3 \cdot 13)^{1/7} \leq n^{1/7}$. 
If $p_3 = 11$ and $n \neq 33$, then $n/\phi(n) = (3/2) \cdot (11/10) < (3\cdot 3 \cdot 11)^{1/7} \leq n^{1/7}$. 
If $p_3 = 7$ and $n \neq 21$, then $n/\phi(n) = (3/2) \cdot (7/6) < (3\cdot 3 \cdot 7)^{1/7} \leq n^{1/7}$. 
If $p_3 = 5$ and $n \neq 15,45,75$, then $n/\phi(n) = (3/2) \cdot (5/4) < (6\cdot 3 \cdot 5)^{1/7} < n^{1/7}$. 
Also, if $n \neq 15, 21,33,45$, then $\phi(n) \geq 24 > 2^{2.2s}$.  

Finally, assume that $s=1$, and $n=3^a$ with $a \geq 3$. Then $n/\phi(n) = 3/2 < 3^{a/7} = n^{1/7}$, and 
$\phi(n) \geq 18 > 2^{4s}$, completing the proof of \eqref{prod10}. 

\smallskip
(ii) 
Now assume that $n \geq 2602$. 
By \eqref{prod10} 
$m:=\phi(n) \geq n^{6/7} > 846$, so $m \geq 847$,
and $2^{\omega(n)} \leq m^{5/11}$. 

Fix $a \in \{1,3\}$ and let $S_j := \{ (j-1)n \leq k < jn \mid \gcd(k,8n)=1,k \equiv a \pmod*{8}\}$ for $0 \leq j \leq 7$. 
For each $j$, observe that $k \in S_j$ if and only if $0 \leq k':= k-jn < n$ is coprime to $n$ and $k' \equiv a-jn \pmod*{8}$. By Proposition
\ref{equi1}, 
$$|S_j| = |\{ 0 \leq k' < n \mid \gcd(k,n) = 1, k' \equiv a-jn \pmod*{8}\}|$$
satisfies $\phi(n)/8-2^{\omega(n)} < |S_j| < \phi(n)/8 + 2^{\omega(n)}$. 
By the above, 
\begin{equation}\label{prod11}
  m/8-m^{5/11} < |S_j| < m/8+m^{5/11}.
\end{equation}  
Now, if $k \in S_0 \cup S_7$, then $3-2\sqrt{2}\cos(k\pi/4n) \geq 3-2\sqrt{2}$. If $k \in S_1 \cup S_6$, then 
$3-2\sqrt{2}\cos(k\pi/4n) \geq 1$. If $k \in S_2 \cup S_5$, then 
$3-2\sqrt{2}\cos(k\pi/4n) \geq 3$. Finally, if $k \in S_3 \cup S_4$, then 
$3-2\sqrt{2}\cos(k\pi/4n) \geq 5$. It follows from \eqref{prod11} that
$$P_{2,a}(n) \geq (3-2\sqrt{2})^{m/4+2m^{5/11}}\cdot 15^{m/4-2m^{5/11}}=A^{m/4}B^{-2m^{5/11}}$$
with $A:= 15(3-2\sqrt{2})$ and $B:=15(3+2\sqrt{2})$.  

Setting $f(t):= A^{t/4}B^{-2t^{5/11}}t^{-7/6}$, 
$$g(t):= \log f(t)= (t/4)\log(A)-2t^{5/11}\log(B)-(7/6)\log(t).$$ 
Now $g'(t)=\log(A)/4-\log(B)/(1.1t^{6/11})-(7/6t)$ is increasing, so $g'(t) \geq g'(847) > 0.13$ when $t \geq 847$. It follows
that $g(m) \geq g(847) > 0.75$, and so $f(m) = \exp(g(m)) > 2.11$. Thus, 
for $m=\phi(n) \geq 847$ 
$$P_{2,a}(n) \geq f(m)m^{7/6} > 2m^{7/6} > 2n,$$
as desired. 
\end{proof}

As we will see in part (iii) of the proof of the following theorem, Proposition \ref{prod-suzuki} actually holds for all odd $n \geq 7$.
 
\begin{thm}\label{prime-suzuki1}
Let $q = 2^n$ with $2 \nmid n$. If $n \geq 7$ then $t(q)=t^-(q) := q-\sqrt{2q}+1$ is divisible by 
a primitive prime divisor $\ppd(2,4n)$ of $q^2+1$. In all cases, $t^+(q):=q+\sqrt{2q}+1$ is divisible by 
a primitive prime divisor $\ppd(2,4n)$ of $q^2+1$.
\end{thm}

\begin{proof}
(i) The second statement is obvious for $n \leq 5$; also note that $t(2)=1$, $t(8)=5$. Henceforth we may assume $n \geq 7$.
Consider the sets
$$\begin{aligned}A_j := \{ 1 \leq  k \leq 4n \mid \gcd(k,2n)=1,~k \equiv j \pmod*{8}\},\\
   B_j := \{ 4n+1 \leq  k \leq 8n \mid \gcd(k,2n)=1,k \equiv j \pmod*{8} \}\end{aligned}$$
for $j = 1,3,5,7$. Then the map $k \mapsto 8n-k$ yields bijections 
$$A_1 \longleftrightarrow B_7,~A_7 \longleftrightarrow B_1,~A_3 \longleftrightarrow B_5,~A_5 \longleftrightarrow B_3,$$ 
and the map $k \mapsto k+4n$ yields bijections 
$$A_1 \longleftrightarrow B_5,~A_5 \longleftrightarrow B_1,~A_3 \longleftrightarrow B_7,~A_7 \longleftrightarrow B_3.$$ 
It follows that
\begin{equation}\label{prim10a}
  |A_1|=|A_3|=|B_5|=|B_7|,~|A_5|=|A_7|=|B_1|=|B_3|.
\end{equation}
Since $\sqcup_{j=1,3,5,7}\bigl( A_j \sqcup B_j\bigr)  = \{1 \leq k < 8n \mid \gcd(k,8n)=1\}$, we now see that 
\begin{equation}\label{prim10b}
  |A_j|+|B_j| = \phi(8n)/4 = \phi(n)
\end{equation}
for each $j = 1,3,5,7$.  

\smallskip
(ii) Make the choice of $\zeta:=\zeta_{8n} = \exp(\pi i/4n)$ and consider the cyclotomic polynomial 
$$\Phi_{8n}(X) = \prod_{1 \leq k < 8n,~\gcd(k,2n)=1}\bigl(X-\zeta^k\bigr).$$
If $1 \leq k < 4n$ then $(X-\zeta^k)(X-\zeta^{k+4n}) = (X-\zeta^k)(X+\zeta^k) = X^2-\zeta_{4n}^k$. It follows that
$$\Phi_{8n}(X) =  \prod_{1 \leq k < 4n,~\gcd(k,2n)=1}\bigl(X^2-\zeta_{4n}^k\bigr) = \Phi_{4n}(X^2).$$
In particular, 
\begin{equation}\label{prim10}
  \Phi_{4n}(2) = \Phi_{8n}(\sqrt{2}).
\end{equation}
Setting 
\begin{equation}\label{prim11}
  \Phi_{8n,a}(X):=  \prod_{1 \leq k < 8n,~\gcd(k,2n)=1,~k \equiv \pm a \pmod*{8}}\bigl(X-\zeta^k\bigr),
\end{equation}  
for $a=1,3$, we have $\Phi_{8n}(X) = \Phi_{8n,1}(X)\Phi_{8n,3}(X)$. Next, since
$$(\sqrt{2}-\zeta^k)(\sqrt{2}-\zeta^{8n-k}) = 3 -2\sqrt{2}\cos \frac{k\pi}{4n}$$
for each $1 \leq k < 8n$, using the bijection 
$k \mapsto 8n-k$ in \eqref{prim10a} 
\begin{equation}\label{prim12}
  P_{2,a}(n)=\Phi_{8n,a}(\sqrt{2}),
\end{equation}
for $a=1,3$. Using \eqref{prim10} also 
\begin{equation}\label{prim13}
  P_{2,1}(n)P_{2,3}(n) = \Phi_{4n}(2).
\end{equation} 

\smallskip
(iii) To show that $P_{2,1}(n)$ and $P_{2,3}(n)$ are integers, we use $\gcd(8,n)=1$ to 
write $1 =ns+8t$ for some $s,t \in \Z$ with $\gcd(s,8)=\gcd(t,n)=1$, and set $\zeta=\zeta_{8n}=\alpha\beta$
with $\alpha := \zeta^{ns}$, a $8^{\mathrm {th}}$ root of unity, and $\beta := \zeta^{8t}$, an $n^{\mathrm {th}}$ root of unity. 
When $k$ runs over $A_1 \sqcup B_1$, 
$\zeta^k = \alpha^k\beta^k = \alpha\beta^k$, and $k \pmod*{n}$ runs over 
units in $\Z/n\Z$, each at most once. Using \eqref{prim10b}, we see that each unit is met exactly once. Repeating the same argument
for $k \in A_7 \sqcup B_7$, we get 
$$\Phi_{8n,1}(X) :=\prod_{l {\rm \ unit\ mod\ } n}\bigl(X-\alpha\beta^l\bigr)\bigl(X-\alpha^{-1}\beta^l\bigr)
    =\prod_{l {\rm \ unit\ mod\ } n}(X^2-(\alpha+\alpha^{-1})\zeta_n^lX +\zeta_{n}^{2l}).$$
Note that $\alpha+\alpha^{-1} = \eps\sqrt{2}$ for some $\eps = \pm 1$. It follows from \eqref{prim12} that
$$P_{2,1}(n) = \Phi_{8n,1}(\sqrt{2}) = \prod_{l {\rm \ unit\ mod\ } n}(2-2\eps\zeta_n^l +\zeta_{n}^{2l})$$  
is the norm $\Norm_{\Q(\zeta_n)/\Q}$
of the algebraic integer $2-2\eps\zeta_n+\zeta_n^2$, hence it is an integer. The same arguments show
that $P_{2,3}(n)$ is the norm $\Norm_{\Q(\zeta_n)/\Q}$ of the algebraic integer $2+2\eps\zeta_n+\zeta_n^2$, hence it is an integer.   

Using this norm interpretation for $P_{2,1}(n)$ and $P_{2,3}(n)$, a calculation with
{\sc Magma} shows that $P_{2,a}(n) > 2n$ for odd integers $7 \leq n \leq 2601$. Together with Proposition \ref{prod-suzuki}, this shows
that 
\begin{equation}\label{prod-norm}
  P_{2,a}(n) > 2n
\end{equation}
for $a=1,3$ and odd $n \geq 7$.  
    
\smallskip
(iv) We also set
$$f^-(X) := X^{2n}-\sqrt{2}X^n +1,~f^+(X) := X^{2n}+\sqrt{2}X^n +1,$$
so that 
$$t^-(q)=f^-(\sqrt(2)),~t^+(q)=f^+(\sqrt(2)),~f^-(X)f^+(X)= X^{4n}+1.$$
Certainly, any root of $X^{4n}+1$ is $\zeta^k$ for some odd integer $1 \leq k < 8n$. If $k \equiv \pm 1 \pmod*{8}$,
then 
$$f^-(\zeta^k) = \exp(k\pi i/2)+1-\sqrt{2}\exp(k\pi i/4) = 0.$$
Similarly, if  $k \equiv \pm 3 \pmod*{8}$,
then 
$$f^+(\zeta^k) = \exp(k\pi i/2)+1+\sqrt{2}\exp(k\pi i/4) = 0.$$
It follows that 
$$f^-(X) = \prod_{1 \leq k < 8n,~k \equiv \pm 1 \pmod*{8}}\bigl(X-\zeta^k\bigr),~
    f^+(X) = \prod_{1 \leq k < 8n,~k \equiv \pm 3 \pmod*{8}}\bigl(X-\zeta^k\bigr).$$
Comparing to \eqref{prim11} and \eqref{prim12}, we see that 
$$f^-(\sqrt{2})/P_{2,1}(n) = \prod_{1 \leq k < 8n,~\gcd(k,2n) > 1,~k \equiv \pm 1 \pmod*{8}}\bigl(\sqrt{2}-\zeta^k\bigr)$$  
is an algebraic integer. But $f^-(\sqrt{2})=t^-(q)$ is an integer, and $P_{2,1}(n)$ is an integer by (iii). Hence
$P_{2,1}(n)$ divides $t^-(q)$. Similarly, $P_{2,3}(n)$ divides $t^+(q)$. 

\smallskip
(v) By \eqref{prod-norm}, $P_{2,1}(n) > n$. 
Consider any prime divisor $\ell$ of $P_{2,1}(n)$, which then divides 
$t^-(q)$ by the result of (iv), and divides $\Phi_{4n}(2)$ by \eqref{prim13}. Suppose that $\ell$ is not a primitive prime divisor of
$2^{4n}-1$. By \cite[Satz 1]{Lun} (cf.\ \cite[Proposition 2]{Ro})
$\ell|n$, and moreover $\ell^2 \nmid \Phi_{4n}(2)$. It follows that
the $\ell$-part of $P_{2,1}(n)$ is $\ell$. Hence, if $t^-(q)$ is not divisible by any primitive prime divisor of $2^{4n}-1$, then 
$P_{2,1}(n)$ divides $n$, a contradiction. 

The proof for $t^+(q)$ is entirely similar.
\end{proof}

\begin{cor}\label{prime-suzuki2}
Let $q = 2^n$ with $2 \nmid n$. If $n \geq 3$ then $\Phi''_{24}:=q^2-q\sqrt{2q}+q-\sqrt{2q}+1$ is divisible by 
a primitive prime divisor $\ppd(2,12n)$ of $q^4-q^2+1$. In all cases, $\Phi'_{24}:=q^2+q\sqrt{2q}+q+\sqrt{2q}+1$ is divisible by 
a primitive prime divisor $\ppd(2,12n)$ of $q^4-q^2+1$.
\end{cor}

\begin{proof}
Note that if $n=1$ then $\Phi''_{24}=1$ and $\Phi'_{24}=13$. Assume that $n \geq 3$. By Theorem \ref{prime-suzuki1},
$$t^-(q^3) = q^3-q\sqrt{2q}+1 = (q+\sqrt{2q}+1)(q^2-q\sqrt{2q}+q-\sqrt{2q}+1) = t^+(q)\Phi''_{24}$$
is divisible by a primitive prime divisor $\ell_1 = \ppd(2,12n)$. 
Since $t^+(q)|(q^2+1)$ and $\ell_1|(q^4-q^2+1)$, we see that $\ell_1 \nmid t^+(q)$, and so
$\ell_1|\Phi''_{24}$. The argument for $\Phi'_{24}$ is similar, using 
$$t^+(q^3) = q^3+q\sqrt{2q}+1 = (q-\sqrt{2q}+1)(q^2+q\sqrt{2q}+q+\sqrt{2q}+1) = t^-(q)\Phi'_{24}.$$
\end{proof}

\subsection*{3C. Primitive prime divisor for Ree groups}
The results in this subsection are not needed for the rest of the paper, however they will be used elsewhere \cite{KT21}.
We make the choice $\sqrt{3} > 0$. For odd $n \in \Z_{\geq 1}$ and $a=1,5$, set
$$P_{3,a}(n):= \prod_{1 \leq k < 12n,~\gcd(k,n)=1,~k \equiv a \pmod*{12}}\bigl(4 - 2\sqrt{3}\cos\frac{k\pi}{6n}\bigr).$$

\begin{prop}\label{prod-ree}
If $2 \nmid n \geq 3$ and $a=1,5$, then $P_{3,a}(n) > 2n$. 
\end{prop} 

\begin{proof}
A computation with Mathematica shows that $P_{3,a}(n) > 2n$ when $3 \leq n \leq 353$. 
Now assume that $n \geq 354$. 
By \eqref{prod10} $m:=\phi(n) \geq n^{6/7} > 153$, so $m \geq 154$,
and $2^{\omega(n)} \leq m^{5/11}$. 

Fix $a \in \{1,5\}$ and let $R_j := \{ (j-1)n \leq k < jn \mid \gcd(k,12n)=1,k \equiv a \pmod*{12}\}$ for $0 \leq j \leq 11$. 
For each $j$, observe that $k \in R_j$ if and only if $0 \leq k':= k-jn < n$ is coprime to $n$ and $k' \equiv a-jn \pmod*{12}$. 
If $3 \nmid n$, then according to Proposition \ref{equi1}, 
$$|R_j| = |\{ 0 \leq k' < n \mid \gcd(k,n) = 1, k' \equiv a-jn \pmod*{12}\}|$$
satisfies $\phi(n)/12-2^{\omega(n)} < |R_j| < \phi(n)/12 + 2^{\omega(n)}$. On the other hand, if $3|n$ then 
$3 \nmid (a-jn)$, so by Proposition \ref{equi2}, 
$$|R_j| = |\{ 0 \leq k' < n \mid \gcd(k,n) = 1, k' \equiv a-jn \pmod*{12}\}|$$
satisfies $\phi(n)/8-2^{\omega(n)-1} < |R_j| < \phi(n)/8 + 2^{\omega(n)-1}$. 

Setting $b:=12$ when $3 \nmid n$ and $b:=8$ when
$3|n$, by the above consideration now 
\begin{equation}\label{prod21}
  m/b-m^{5/11} < |R_j| < m/b+m^{5/11}.
\end{equation}  
Now, if $k \in S_0 \cup S_{11}$, then $4-2\sqrt{3}\cos(k\pi/6n) \geq 4-2\sqrt{3}$. If $k \in S_1 \cup S_{10}$, then 
$4-2\sqrt{3}\cos(k\pi/6n) \geq 1$. If $k \in S_2 \cup S_9$, then 
$4-2\sqrt{3}\cos(k\pi/6n) \geq 4-\sqrt{3}$. If $k \in S_3 \cup S_8$, then 
$4-2\sqrt{3}\cos(k\pi/6n) \geq 4$. If $k \in S_4 \cup S_7$, then 
$4-2\sqrt{3}\cos(k\pi/6n) \geq 4+\sqrt{3}$. Finally, if $k \in S_5 \cup S_6$, then 
$4-2\sqrt{3}\cos(k\pi/6n) \geq 7$.
It follows from \eqref{prod21} that
$$\begin{aligned}
   P_{3,a}(n) & \geq \bigl(4-2\sqrt{3}\bigr)^{2m/b+2m^{5/11}}\cdot \bigl((4-\sqrt{3})\cdot 4 \cdot (4+\sqrt{3}) \cdot 7 \bigr)^{2m/b-2m^{5/11}}\\
    & =A^{2m/b}B^{-2m^{5/11}} \geq A^{m/6}B^{-2m^{5/11}}\end{aligned}$$
with $A:= 364(4-2\sqrt{3})$ and $B:=364/(4-2\sqrt{3})$.  

Setting $f(t):= A^{t/6}B^{-2t^{5/11}}t^{-7/6}$, 
$$g(t):= \log f(t)= (t/6)\log(A)-2t^{5/11}\log(B)-(7/6)\log(t).$$ 
Now $g'(t)=\log(A)/6-\log(B)/(1.1t^{6/11})-7/6t$ is increasing, so $g'(t) \geq g'(154) > 0.48$ when $t \geq 154$. It follows
that $g(m) \geq g(154) > 0.74$, and so $f(m) = \exp(g(m)) > 2.09$. Thus, for $m=\phi(n) \geq 154$ 
$$P_{3,a}(n) \geq f(m)m^{7/6} > 2m^{7/6} > 2n,$$
as desired. 
\end{proof}

\begin{thm}\label{prime-ree}
Let $q = 3^n$ with $2 \nmid n$. If $n \geq 3$ then $t(q)=t^-(q) := q-\sqrt{3q}+1$ is divisible by 
a primitive prime divisor $\ppd(3,6n)$ of $q^2-q+1$. In all cases, $t^+(q):=q+\sqrt{3q}+1$ is divisible by 
a primitive prime divisor $\ppd(3,6n)$ of $q^2-q+1$.
\end{thm}

\begin{proof}
(i) Note that $t(3)=1$. Henceforth we will assume $n \geq 3$.
Consider the sets
$$\begin{aligned}A_j := \{ 1 \leq  k \leq 6n \mid \gcd(k,12n)=1,~k \equiv j 
\pmod*{12}\},\\
   B_j := \{ 6n+1 \leq  k \leq 12n \mid \gcd(k,12n)=1,k \equiv j \pmod*{12} \}\end{aligned}$$
for $j = 1,5,7,11$. Then the map $k \mapsto 12n-k$ yields bijections 
$$A_1 \longleftrightarrow B_{11},~A_{11} \longleftrightarrow B_1,~A_5 \longleftrightarrow B_7,~A_7 \longleftrightarrow B_5,$$ 
and the map $k \mapsto k+6n$ yields bijections 
$$A_1 \longleftrightarrow B_7,~A_7 \longleftrightarrow B_1,~A_5 \longleftrightarrow B_{11},~A_{11} \longleftrightarrow B_5.$$ 
It follows that
\begin{equation}\label{prim20a}
  |A_1|=|A_5|=|B_7|=|B_{11}|,~|A_7|=|A_{11}|=|B_1|=|B_5|.
\end{equation}
Since $\sqcup_{j=1,5,7,11}\bigl( A_j \sqcup B_j\bigr)  = \{1 \leq k < 12n \mid \gcd(k,12n)=1\}$, we now see that 
\begin{equation}\label{prim20b}
  |A_j|+|B_j| = \phi(12n)/4 = \phi(3n)/2
\end{equation}
for each $j = 1,5,7,11$.  

\smallskip
(ii) Make the choice of $\zeta:=\zeta_{12n} = \exp(\pi i/6n)$ and consider the cyclotomic polynomial 
$$\Phi_{12n}(X) = \prod_{1 \leq k < 12n,~\gcd(k,12n)=1}\bigl(X-\zeta^k\bigr).$$
If $1 \leq k < 6n$ then $(X-\zeta^k)(X-\zeta^{k+6n}) = (X-\zeta^k)(X+\zeta^k) = X^2-\zeta_{6n}^k$. It follows that
$$\Phi_{12n}(X) =  \prod_{1 \leq k < 6n,~\gcd(k,6n)=1}\bigl(X^2-\zeta_{6n}^k\bigr) = \Phi_{6n}(X^2).$$
In particular, 
\begin{equation}\label{prim20}
  \Phi_{6n}(3) = \Phi_{12n}(\sqrt{3}).
\end{equation}
Setting 
\begin{equation}\label{prim21}
  \Phi_{12n,a}(X):=  \prod_{1 \leq k < 12n,~\gcd(k,12n)=1,~k \equiv \pm a 
\pmod*{12}}\bigl(X-\zeta^k\bigr),
\end{equation}  
for $a=1,5$, we have $\Phi_{12n}(X) = \Phi_{12n,1}(X)\Phi_{12n,5}(X)$. Next, since
$$(\sqrt{3}-\zeta^k)(\sqrt{3}-\zeta^{12n-k}) = 4 -2\sqrt{3}\cos \frac{k\pi}{6n}$$
for each $1 \leq k < 12n$, using the bijection $k \mapsto 12n-k$ in \eqref{prim20a} 
\begin{equation}\label{prim22}
  P_{3,a}(n)=\Phi_{12n,a}(\sqrt{3}),
\end{equation}
for $a=1,5$. Using \eqref{prim20} also 
\begin{equation}\label{prim23}
  P_{3,1}(n)P_{3,5}(n) = \Phi_{6n}(3).
\end{equation} 

\smallskip
(iii) Here we show that $P_{3,1}(n)$ and $P_{3,5}(n)$ are integers. Clearly, they are algebraic integers in 
$\Q(\zeta)$. Hence it suffices to show that each of them is fixed by any Galois automorphism 
$\sigma: \zeta \mapsto \zeta^l$, $\gcd(l,12n) = 1$. First consider the case $j \equiv \pm 1 \pmod*{12}$. Then $\sigma$ fixes
$\sqrt{3}=\zeta_{12}+\zeta_{12}^{-1}$, and fixes each of the sets $C:=\sqcup_{j = 1,11}(A_j \cup B_j)$ and 
$D:=\sqcup_{j = 1,11}(A_j \cup B_j)$ modulo $12n$. Since 
$$P_{3,1}(n) = \prod_{k \in C}(\sqrt{3}-\zeta^k),~P_{3,5}(n) = \prod_{k \in D}(\sqrt{3}-\zeta^k),$$
it follows that $\sigma$ fixes each of $P_{3,1}(n)$ and $P_{3,5}(n)$. Now assume that $j \equiv \pm 5 \pmod*{12}$.
Then $\sigma$ sends $\sqrt{3}$ to $-\sqrt{3}$ and $\zeta^k$ to $\zeta^{kl}$, and thus 
$$\sigma\bigl(\sqrt{3}-\zeta^k\bigr) = -\sqrt{3}-\zeta^{kl} = -\bigl(\sqrt{3}-\zeta^{6n+kl}\bigr).$$ 
Note that modulo $12n$, when $k$ runs over $C$, $6n+kl$ runs over $C$, covering each element of $C$ exactly once. Also,
$|C|=\phi(3n)$ by \eqref{prim20b}, and so $|C|$ is even. It follows that $\sigma$ sends $P_{3,1}(n)$ to 
$(-1)^{|C|}P_{3,1}(n)=P_{3,1}(n)$, and similarly $\sigma$ fixes $P_{3,5}(n)$.
      
\smallskip
(iv) We also set
$$f^-(X) := X^{2n}-\sqrt{3}X^n +1,~f^+(X) := X^{2n}+\sqrt{3}X^n +1,$$
so that 
$$t^-(q)=f^-(\sqrt(3)),~t^+(q)=f^+(\sqrt(3)),~f^-(X)f^+(X)= X^{4n}-X^{2n}+1.$$
Certainly, any root of $X^{4n}-X^{2n}+1$ is $\zeta^k$ for some integer $1 \leq k < 12n$ coprime to $6$. If $k \equiv \pm 1 \pmod*{12}$,
then 
$$f^-(\zeta^k) = \exp(k\pi i/3)+1-\sqrt{3}\exp(k\pi i/6) = 0.$$
Similarly, if  $k \equiv \pm 5 \pmod*{12}$,
then 
$$f^+(\zeta^k) = \exp(k\pi i/3)+1+\sqrt{3}\exp(k\pi i/6) = 0.$$
It follows that 
$$f^-(X) = \prod_{1 \leq k < 12n,~k \equiv \pm 1 \pmod*{12}}\bigl(X-\zeta^k\bigr),~
    f^+(X) = \prod_{1 \leq k < 12n,~k \equiv \pm 5 \pmod*{12}}\bigl(X-\zeta^k\bigr).$$
Comparing to \eqref{prim21} and \eqref{prim22}, we see that 
$$f^-(\sqrt{3})/P_{3,1}(n) = \prod_{1 \leq k < 8n,~\gcd(k,12n) > 1,~k \equiv \pm 1 \pmod*{12}}\bigl(\sqrt{3}-\zeta^k\bigr)$$  
is an algebraic integer. But $f^-(\sqrt{3})=t^-(q)$ is an integer, and $P_{3,1}(n)$ is an integer by (iii). Hence
$P_{3,1}(n)$ divides $t^-(q)$. Similarly, $P_{3,5}(n)$ divides $t^+(q)$. 

\smallskip
(v) By Proposition \ref{prod-ree}, $P_{3,1}(n) > n$. 
Consider any prime divisor $\ell$ of $P_{3,1}(n)$, which then divides 
$t^-(q)$ by the result of (iv), and divides $\Phi_{6n}(3)$ by \eqref{prim23}. Since $\ell|(q^2-q+1)$, $\ell \neq 2,3$.
Suppose that $\ell$ is not a primitive prime divisor of
$3^{6n}-1$. Again by \cite[Satz 1]{Lun} $\ell|3n$, whence $\ell|n$ as $\ell \geq 5$, and moreover $\ell^2 \nmid \Phi_{6n}(3)$. It follows that the $\ell$-part of $P_{3,1}(n)$ is $\ell$. Hence, if $t^-(q)$ is not divisible by any primitive prime divisor of $3^{6n}-1$, then 
$P_{3,1}(n)$ divides $n$, a contradiction. 

The proof for $t^+(q)$ is entirely similar.
\end{proof}


\subsection*{3D. Primitive prime divisors for Suzuki-Ree groups: another approach}

For $\delta\in \Z^+$, $\alpha\in (\Z/\delta)^\times$, and $x\in \R^+$ with $x\geq 1$, let $$f_n^{(\alpha\bmod{\delta})}(x) := \prod_{a\in (\Z/\delta n)^\times : a\equiv \pm \alpha\pmod*{\delta}} (x - \zeta_{\delta n}^a).$$ Note that for such $\delta, \alpha$, and $x$, $f_n^{(\alpha\bmod{\delta})}(x)\in \R^+$ by pairing $a$ with $-a$. Note also that
\begin{align*}
f_n^{(1\bmod{8})}(\sqrt{2}) &= P_{2,1}(n) = \Phi_{8n,1}(\sqrt{2}),
\\f_n^{(3\bmod{8})}(\sqrt{2}) &= P_{2,3}(n) = \Phi_{8n,3}(\sqrt{2}),
\\f_n^{(1\bmod{12})}(\sqrt{3}) &= P_{3,1}(n) = \Phi_{12n,1}(\sqrt{3}),
\\f_n^{(5\bmod{12})}(\sqrt{3}) &= P_{3,5}(n) = \Phi_{12n,5}(\sqrt{3})
\end{align*}\noindent
in the notation above.

\begin{lem}\label{increasing lemma}
Let $n\in \Z^+$. Let $x\geq 1$. Then
$$f_n^{(\alpha\bmod{\delta})}(x)\geq f_n^{(\alpha\bmod{\delta})}(1)\cdot \left(\frac{x+1}{2}\right)^{\frac{2\cdot \phi(\delta n)}{\phi(\delta)}}.$$
\end{lem}

\begin{proof}
The claim is evident for $x = 1$, and we will show that 
$$\frac{f_n^{(\alpha\bmod{\delta})}(x)}{(x+1)^{\frac{2\cdot \phi(\delta n)}{\phi(\delta)}}}$$ is increasing in $x$. Indeed 
$$\frac{f_n^{(\alpha\bmod{\delta})}(x)}{(x+1)^{\frac{2\cdot \phi(\delta n)}{\phi(\delta)}}} = \prod_{a\in (\Z/\delta n)^\times :\, a\equiv \pm \alpha\pmod*{\delta}} \frac{|x - \zeta_n^a|}{x+1},$$ 
and it suffices to show each factor is increasing in $x$. But for all $z\in S^1$ $$\frac{|x-z|^2}{(x+1)^2} = \frac{x^2 - 2x\cdot \Re{z} + 1}{x^2 + 2x + 1} = 1 - \frac{2x}{(x+1)^2}\cdot (1 + \Re{z}),$$ which is increasing in $x$ because $\frac{x}{(x+1)^2} = \frac{1}{x+1} - \frac{1}{(x+1)^2}$ decreases in $x$ when $x\geq 1$.
\end{proof}

For $a,b\in \Z^+$, write $\gcd(a,b^\infty) := \prod_{p\vert b} p^{v_p(a)}$, and $\mathrm{rad}(a) := \prod_{p\vert a} p$.

\begin{lem}\label{calculation of the constant}
Let $n\in \Z^+$. Let $m := \frac{n}{\gcd(n, \delta^\infty)}$. Then
$$f_n^{(\alpha\bmod{\delta})}(1) = \prod_{S\subseteq \{p\vert m\}} \left|1 - \zeta_\delta^{\alpha\cdot \prod_{p\in S} p^{-1}\pmod*{\delta}}\right|^{2\cdot (-1)^{\#|S|}}.$$
\end{lem}

\begin{proof}
We first claim that if $\mathrm{rad}{\left(\frac{n}{\gcd(n,\delta^\infty)}\right)} = \mathrm{rad}{\left(\frac{n'}{\gcd(n',\delta^\infty)}\right)}$ then $f_n^{(\alpha\bmod{\delta})}(1) = f_{n'}^{(\alpha\bmod{\delta})}(1)$. This follows by repeatedly applying the following. If $p\vert n$ is such that $k := v_p(\delta n)\geq 2$, then, writing $\delta n =: p^k\cdot s$, because $X^p - Y^p = \prod_{b\in \F_p} (X - \zeta_p^b\cdot Y)$ as elements of $\Z[\zeta_p][X,Y]$,
\begin{align*}
f_n^{(\alpha\bmod{\delta})}(1) &= \prod_{a\in (\Z/p^{k-1} s)^\times : a\equiv \pm \alpha\pmod*{\delta}} \prod_{b\in \F_p} (1 - \zeta_{p^k\cdot s}^{a + p^{k-1}\cdot s\cdot b})
\\&= \prod_{a\in (\Z/p^{k-1} s)^\times : a\equiv \pm \alpha\pmod*{\delta}} (1 - \zeta_{p^{k-1}\cdot s}^a).
\end{align*}\noindent
Therefore without changing $f_n^{(\alpha\bmod{\delta})}(1)$ we may assume without loss of generality that $n$ is squarefree and such that $\gcd(n,\delta) = 1$.

Now if $p\vert n$, writing $\delta n =: p\cdot n_0$ and letting $p',n_0'\in \Z$ be such that $p p' + n_0 n_0' = 1$, \begin{align*}
f_n^{(\alpha\bmod{\delta})}(1) &= \prod_{a\in (\Z/n_0)^\times : a\equiv \pm \alpha\pmod*{\delta}} \prod_{b\in \F_p^\times} (1 - \zeta_{p n_0}^{a p p' + b n_0 n_0'})
\\&= \prod_{a\in (\Z/n_0)^\times : a\equiv \pm \alpha\pmod*{\delta}} \prod_{b\in \F_p^\times} (1 - \zeta_{n_0}^{a p'}\cdot \zeta_p^{b n_0'})
\\&= \prod_{a\in (\Z/n_0)^\times : a\equiv \pm \alpha\pmod*{\delta}} \prod_{b\in \F_p^\times} (1 - \zeta_{n_0}^{a p'}\cdot \zeta_p^b)
\end{align*}\noindent
via $b\mapsto p\cdot b$. Now we apply the identity $\frac{X^p - Y^p}{X - Y} = \prod_{b\in \F_p^\times} (X - \zeta_p^b\cdot Y)$ to find:
\begin{align*}
f_n^{(\alpha\bmod{\delta})}(1) &= \prod_{a\in (\Z/n_0)^\times : a\equiv \pm \alpha\pmod*{\delta}} \frac{(1 - \zeta_{n_0}^a)}{(1 - \zeta_{n_0}^{a p'})}
\\&= \frac{f_{\frac{n}{p}}^{(\alpha\bmod{\delta})}(1)}{f_{\frac{n}{p}}^{(\alpha p'\bmod{\delta})}(1)}.
\end{align*}

Since the lemma is evident for $n = 1$, by induction on the number of prime factors of $n$ we find that
\begin{align*}
f_n^{(\alpha\bmod{\delta})}(1) &= \prod_{S\subseteq \{\ell\vert \frac{n}{p}\}} \left(\frac{\left|1 - \zeta_\delta^{\alpha\cdot \prod_{\ell\in S} \ell^{-1}\pmod*{\delta}}\right|}{\left|1 - \zeta_\delta^{\alpha\cdot p^{-1}\cdot \prod_{\ell\in S} \ell^{-1}\pmod*{\delta}}\right|}\right)^{2\cdot (-1)^{\#|S|}}
\\&= \prod_{S\subseteq \{\ell\vert n\}} \left|1 - \zeta_\delta^{\alpha\cdot \prod_{\ell\in S} \ell^{-1}\pmod*{\delta}}\right|^{2\cdot (-1)^{\#|S|}},
\end{align*}\noindent
and we are done.
\end{proof}

\begin{cor}\label{one mod eight corollary}
Let $n\in \Z^+$ with $n\geq 3$ be odd. Then
$$f_n^{(1\bmod{8})}(1) = \begin{cases} 1 & \exists p\vert n : p\equiv \pm 1\pmod*{8},\\ (1+\sqrt{2})^{-2^{\omega(n)}} & \text{else}\end{cases}$$
and
$$f_n^{(3\bmod{8})}(1) = \begin{cases} 1 & \exists p\vert n : p\equiv \pm 1\pmod*{8},\\ (1+\sqrt{2})^{2^{\omega(n)}} & \text{else,}\end{cases}$$
whence
$$f_n^{(1\bmod{8})}(\sqrt{2})\geq \left(\frac{\sqrt{2}+1}{2}\right)^{2\cdot \phi(n)}\cdot \begin{cases} 1 & \exists p\vert n : p\equiv \pm 1\pmod*{8},\\ (1+\sqrt{2})^{-2^{\omega(n)}} & \text{else}\end{cases}$$
and
$$f_n^{(3\bmod{8})}(\sqrt{2})\geq \left(\frac{\sqrt{2}+1}{2}\right)^{2\cdot \phi(n)}\cdot \begin{cases} 1 & \exists p\vert n : p\equiv \pm 1\pmod*{8},\\ (1+\sqrt{2})^{2^{\omega(n)}} & \text{else.}\end{cases}$$
\end{cor}

\begin{proof}
If there is a $p\vert n$ with $p\equiv \pm 1\pmod*{8}$, then pair $S - \{p\}$ with $S\cup \{p\}$ in the conclusion of Lemma \ref{calculation of the constant} to see that $f_n^{(1\bmod{8})}(1) = f_n^{(3\bmod{8})}(1) = 1$. Otherwise every $p\vert n$ is such that $p\equiv \pm 3\pmod*{8}$, so that by Lemma \ref{calculation of the constant} $$f_n^{(1\bmod{8})}(1) = \left(\frac{|1 - \zeta_8|}{|1 - \zeta_8^3|}\right)^{2^{\omega(n)}} = (1 + \sqrt{2})^{-2^{\omega(n)}}.$$ By the same reasoning $$f_n^{(3\bmod{8})}(1) = \left(\frac{|1 - \zeta_8^3|}{|1 - \zeta_8|}\right)^{2^{\omega(n)}} = (1 + \sqrt{2})^{2^{\omega(n)}}.$$ We are done by Lemma \ref{increasing lemma}.
\end{proof}

\begin{cor}\label{one mod twelve corollary}
Let $n\in \Z^+$ with $n\geq 3$ be odd. Let $m := \frac{n}{\gcd(n,6^\infty)}$. Then 
$$f_n^{(1\bmod{12})}(1) = \begin{cases} 1 & \exists p\vert m : p\equiv \pm 1\pmod*{12},\\ (2+\sqrt{3})^{-2^{\omega(m)}} & \text{else}\end{cases}$$
and
$$f_n^{(5\bmod{12})}(1) = \begin{cases} 1 & \exists p\vert m : p\equiv \pm 1\pmod*{12},\\ (2+\sqrt{3})^{2^{\omega(m)}} & \text{else,}\end{cases}$$
whence
$$f_n^{(1\bmod{12})}(\sqrt{3})\geq \left(\frac{\sqrt{3}+1}{2}\right)^{\frac{\phi(12 n)}{2}}\cdot \begin{cases} 1 & \exists p\vert m : p\equiv \pm 1\pmod*{12},\\ (2+\sqrt{3})^{-2^{\omega(m)}} & \text{else}\end{cases}$$
and
$$f_n^{(5\bmod{12})}(\sqrt{3})\geq \left(\frac{\sqrt{3}+1}{2}\right)^{\frac{\phi(12 n)}{2}}\cdot \begin{cases} 1 & \exists p\vert m : p\equiv \pm 1\pmod*{12},\\ (2+\sqrt{3})^{2^{\omega(m)}} & \text{else.}\end{cases}$$
\end{cor}

\begin{proof}
If there is a $p\vert m := \frac{n}{\gcd(n,6^\infty)}$ with $p\equiv \pm 1\pmod*{12}$, then pair $S - \{p\}$ with $S\cup \{p\}$ in the conclusion of Lemma \ref{calculation of the constant} to see that $f_n^{(1\bmod{12})}(1) = f_n^{(5\bmod{12})}(1) = 1$. Otherwise every $p\vert m$ is such that $p\equiv \pm 5\pmod*{12}$, so that by Lemma \ref{calculation of the constant} $$f_n^{(1\bmod{12})}(1) = \left(\frac{|1 - \zeta_{12}|}{|1 - \zeta_{12}^5|}\right)^{2^{\omega(m)}} = (2 + \sqrt{3})^{-2^{\omega(m)}}.$$ By the same reasoning $$f_n^{(5\bmod{12})}(1) = \left(\frac{|1 - \zeta_{12}^5|}{|1 - \zeta_{12}|}\right)^{2^{\omega(m)}} = (2 + \sqrt{3})^{2^{\omega(m)}}.$$ We are done by Lemma \ref{increasing lemma}.
\end{proof}

\section{Action on $2$-groups and primitivity of local systems}
We begin this section with a group theoretic lemma which will be used in the proof of primitivity given in Theorem \ref{prim}.

\begin{lem}\label{action}
Suppose $A$ is a cyclic group of  prime order $p$ that
acts faithfully on a finite $q$-group $G$, where $p\ne q$ are primes.
Let $n$ be the order of $q$ modulo $p$.
Let $\chi \in {\rm Irr}(G)$ be $A$-invariant and faithful, and write $\chi(1)=q^a$. Then $n \le 2a$.
\end{lem}

\begin{proof}
We argue by induction on $|G|$. 
Let $C=\cent GA<G$. Let $C\le N<G$ be maximal $A$-invariant in $G$.
Since $G$ is nilpotent, $N \nor G$. Also, $G/N$ is an irreducible
$A$-module. Notice that $A$ cannot act trivially on $G/N$, because $G=[G,A]C$, by coprime action. Hence, $G/N$ is a faithful irreducible $\F_q[A]$-module.
By  \cite[Example 2.7]{MW}, say, $|G/N|=q^n$. 
Let $\theta \in {\rm Irr}(N)$ be $A$-invariant under $\chi$; such exists by \cite[Theorem 13.27]{Is}.
Now, since $G/N$ is an abelian
 chief factor of the semidirect product 
 $GA$ and $\chi$ is $GA$-invariant, by the Isaacs ``going down'' theorem \cite[Theorem 6.18]{Is},  
either
$\chi_N=\theta$, or $\chi_N=e\theta$ with $e^2=|G/N|$, or $\theta^G=\chi$.

In the third case, $\chi(1)=q^n \theta(1)=q^{n+b}=q^a$, for some $b\ge 0$.
Thus $n\le n+b=a \le 2a$. In the second case,
$q^n$ is a square, and that $\chi(1)/\theta(1)=q^{n/2}$. Then $\chi(1)=q^{{n\over 2} +c}=q^a$, for some $c\ge 0$.
Then $a={n \over 2} + c$, and again $n \le 2a$.

In the first case, $\chi_N=\theta \in \irr N$.
If $A$ does not act trivially on $N$, then $A$ acts faithfully on $N$.
Since $\theta$ is $A$-invariant and faithful, then we are done by the inductive hypothesis.
Otherwise, $A$ acts trivially on $N$, and therefore
 $N=\cent GA$ (because $N$ is maximal $A$-invariant).  
 By \cite[Lemma 2.1]{N},
$[G,A] \sbs \Ker(\chi)$. Since $\chi$ is faithful, $A$ acts trivially on $G$.
But this cannot happen. 
\end{proof}

We will now introduce a general class of Airy sheaves $\sF(q,f)$ that includes the sheaves $\sF_q$. Recall that 
$q=2^{2n+1}=2q_0^2$ and $t(q):=q-2q_0+1$.
Let $k_0/\F_2$ be a finite extension, and $f(x) \in k_0[x]$ a polynomial of degree $(1+q_0)t(q)$. Form the Artin-Schreier-Witt lisse sheaf
$$\sL(q,f):=\sL_{\psi_2([x^{t(q)},f(x)]}$$
on $\A^1/k_0$. Then $\sF(q,f)$ is the constant field twisted Fourier transform 
\begin{equation}\label{airy-g}
  \sF(q,f):=\FT_\psi(\sL(q,f))\otimes (1-(-1)^ni)^{-\deg}.
\end{equation}  
The trace function of $\sF(q,f)$ at $t \in k$, $k/k_0$ a finite extension, is
\begin{equation}\label{airy-g1}
  t \mapsto \frac{-1}{\bigl((1-(-1)^n\,i\bigr)^{\deg(k/\F_2)}}\sum_{x \in k}\psi_2(\Trace_{W_2(k)/W_2(\F_2)}([x^{t(q)},f(x)+tx])).
\end{equation}

Assume in addition that 
\begin{equation}\label{airy-s}
  f(x) =f_1(x^{t(q)})
\end{equation}
for some polynomial $f_1(x)\in k_0[x]$ of degree $1+q_0$. If we make the substitution $x \mapsto x/t$, then 
\eqref{airy-g1} for $t \in k^\times$ becomes
$$t \mapsto \frac{-1}{\bigl(1-(-1)^n\,i\bigr)^{\deg(k/\F_2)}}\sum_{x \in k}\psi_2\bigl(\Trace_{W_2(k)/W_2(\F_2)}([x^{t(q)}/t^{t(q)},f_1(x^{t(q)}/t^{t(q)})+x])\bigr),$$
For $$t(q)=rs,$$
we get a descent $\sG(q,f,r)$ of $\sF(q,f)$ to $\G_m/k_0$ whose trace function is now
$$t \mapsto \frac{-1}{\bigl(1-(-1)^n\,i\bigr)^{\deg(k/\F_2)}}\sum_{x \in k}\psi_2\bigl(\Trace_{W_2(k)/W_2(\F_2)}([x^{t(q)}/t^s,f_1(x^{t(q)}/t^s)+x])\bigr),$$
and whose Kummer pullback by $r^{\mathrm {th}}$ power  is (the restriction to $\G_m/k_0$ of) $\sF(q,f)$:
$$[r]^\star\sG(q,f,r) \cong \sF(q,f)|\G_m.$$

We next give a key lemma of {\rm \v{S}uch}, which is proved in \cite[Proposition 11.1]{Such} but not stated there explicitly.
\begin{lem} {\rm (\v{S}uch)} \label{ASinduced}
Let $\sF$ be an Airy sheaf of rank $n \ge 2$ on $\A^1/\overline{\F_p}$, i.e.\ $\sF$ is the Fourier transform $\FT(\sL)$ of a lisse, rank
one sheaf $\sL$ on $\A^1/\overline{\F_p}$ with $\Swan_\infty(\sL)=n+1$. If $\sF$ is induced, then it is induced from an Artin-Schreier covering of $\A^1/\overline{\F_p}$. In particular, if $\sF$ is induced,  then 
it is induced from a normal subgroup of index $p$ of its $G_{\geo}$.
\end{lem}
\begin{proof}
Let us recall the argument, which is contained in the proof of \cite[Proposition 11.1]{Such}. If $\sF$ is induced, then it is $g_\star \sH$, in which case 
$$\End(\sF)=(g_\star \sH)\otimes (g_\star \sH^\vee) \supseteq g_\star(\sH\otimes \sH^\vee) \supseteq g_\star\triv,$$
and hence $\End^0(\sF) \supseteq g_\star\triv/\triv$. But $g_\star\triv/\triv$ has rank $<n$, and all its slopes are $\le (n+1)/n = 1 + 1/n$. But each slope of 
 $g_\star\triv/\triv$ has denominator in lowest terms at most the rank of $g_\star\triv/\triv$, which is $< n$. Therefore each slope, being at most
 $1 + 1/n$, is in fact $\le 1$. Then by \cite[Corollary 3.3]{Such}, $g_\star\triv/\triv$ is the direct sum of various $\sL_{\psi(ax)}$.

The rest of the argument
is given in the first eight lines of the second paragraph of the proof of \cite[Proposition 11.1]{Such}.
\end{proof}

We will need the following general result, a slight generalization of \cite[Theorem 4.6]{KT5}.

\begin{thm}\label{Iinftygen}Let $\sF$ be a semisimple lisse sheaf on $\G_m/\overline{\F_p}$ which is tame at $0$. Denote by $J$  the image of $I(\infty)$ in $G:=G_{\geo,\sF}$. Then $G$ is the Zariski closure of the normal subgroup of $G$ generated by all $G$-conjugates of $J$.
\end{thm}

\begin{proof}Denote by $G_\infty$ this Zariski closure. Then $G$ is reductive, and hence its quotient $G/G_\infty$ is reductive. It suffices to show that every irreducible representation of $G/G_\infty$ is trivial. But a $d$-dimensional irreducible representation is a lisse sheaf of rank $d$ on 
$\G_m/\overline{\F_p}$ which is tame at $0$ (because $P(0)$ dies in $G$) and lisse at $\infty$. By multiplicative inversion, this is an irreducible local system
on $\A^1/\overline{\F_p}$ which is tame at $\infty$, so a representation of $\pi_1^{{\mathrm{tame}}}(\A^1/\overline{\F_p)}$, which is the trivial group.
\end{proof}

We will also need a very special case (about subgroups of index $2$) of (the second part of) the following proposition. 
Since we cannot find a reference for it, we give a proof.

\begin{prop}\label{finite-index}
\begin{enumerate}[\rm(i)]
\item Let $G$ be a Lie group and $H$ an abstract subgroup of finite index in $G$. Then $H$ is closed. 
\item Let $G$ be a reductive linear algebraic group over an algebraically closed field $k$, 
and let $H$ be an abstract subgroup of finite index in 
$G$. Then $H$ is Zariski closed.
\end{enumerate}
\end{prop}

\begin{proof}
(a) We give a proof of statement (i), which is due to Jason DeVito, posted on {\tt MathStack\-Exchange}.

\smallskip
(a1) First we show that if $G$ is connected, then $G$ is generated by its divisible subgroups.  Indeed, there is an open set 
$U \subseteq G$ containing the identity such that $U \subseteq \exp(\sL)$ for the Lie algebra $\sL$ of $G$. For any $u$ in $U$, if 
$u = \exp(X)$, then $u$ lies in the divisible subgroup $\{\exp(tX) \mid t \in \R\}$ of $G$. Since $G$ is generated by $U$, the claim follows.

\smallskip 
(a2) Next we show that if $G$ is connected and $H \leq G$ is of finite index, then $H=G$.  Indeed, by considering the action 
via left translation of $G$ on the finite set $G/H$ of left $H$-cosets, we see that $H$ contains a normal subgroup $K \lhd G$
of finite index. By (i), $G$, and so $G/K$, is generated by its divisible subgroups. But $G/K$ is finite, so $K=G$.

\smallskip
(a3) In the general case, consider the map $\iota: G^\circ/(H \cap G^\circ) \to G/H$ defined by $x(H \cap G^\circ) \mapsto xH$. Then 
$\iota$ is injective, and so $H \cap G^\circ$ has finite index in $G^\circ$. By (a2), $H \geq G^\circ$, and so $H$ is a union of 
a finite number of $G^\circ$, each of which is closed. Hence $H$ is closed. 

\smallskip
(b) For statement (ii), 
the argument in (a3) shows that it suffices to prove that if $G$ is connected reductive and $H$ has 
finite index in $G$, then $H=G$. The argument in (a2) allows us to further assume that $H \lhd G$. 
Recall \cite[p. 16]{C} that $G=T[G,G]$, where $T = \ZB(G)^\circ$ is a central torus of $G$, and the derived subgroup
$[G,G]$ is a central product $G_1 \circ \ldots \circ G_n$ of simple groups. In particular, for each $i$, $H \cap G_i$ is a normal
(abstract) subgroup of finite index of $G_i$. As $G_i$ is generated by (unipotent) root subgroups, \cite[Main Theorem]{Ti} implies
that $H \cap G_i$ is either equal to $G_i$ or contained in $\ZB(G_i)$. The finite index assumption now ensures that $H \geq G_i$ for
all $i$, and thus $H \geq [G,G]$. Next, $H \cap T$ is a normal (abstract) subgroup of finite index say $m$ in $T$. In particular,
$H$ contains $t^m$ for all $t \in T$. But the torus $T$ is $(k^\times)^r$ for some $r$, so $k = \overline{k}$ implies that any element
in $T$ is an $m^{\mathrm{th}}$ power, and hence $H \geq T$, completing the proof.
[The referee kindly pointed out a much simpler argument as follows. If $H \lhd G$ has index $N$ in $G$ then $G$ contains all powers $g^N$, $g \in G$. But $G$, being connected reductive, is generated by its maximal tori, and the map $g \mapsto g^N$ is 
surjective on each maximal torus. Hence $H=G$.] 
\end{proof}

As pointed out by the referee, in both of the cases of Proposition \ref{finite-index}, once $H$ is closed, it is also open 
(indeed, $G \smallsetminus H$ is a disjoint union of a finite number of cosets $gH$, each being closed, and so it is closed).

%
%
%

\begin{thm}\label{prim}
Let $q = 2^{2n+1}$ with $n \in \Z_{\geq 1}$ and $n \neq 2$. Under the assumption \eqref{airy-s}, the group $G_{\geo,\sF(q,f)}$ of the Airy sheaf $\sF(q,f)$ has no subgroups of index $2$. As a consequence, the Airy sheaf $\sF(q,f)$ is not geometrically induced, and hence 
none of the sheaves $\sG(q,f,r)$ is geometrically induced.
\end{thm}

\begin{proof}
Assume to the contrary that the Airy sheaf $\sF(q,f)$ is induced. Then, by Lemma \ref{ASinduced}, the underlying representation
$V$ of the geometric monodromy group $G:=G_{\geo}$ of $\sF(q,f)$ is induced from a subgroup $G_1$ of $G$ of
index $2$.

\smallskip
(i) For each divisor $r > 1$ of $t(q)$, consider the descent $\sG(q,f,r)$. Because $\sF(q,f)$ is lisse at $0$, 
the image of $I(0)$ in the geometric monodromy group $H$ of $\sG(q,f,r)$ is the cyclic group $\mu_r(\overline{\F_2})$ of order $r$. 
Because the $\infty$-slopes of $\sF(q,f)$ are $\dfrac{(2^n+1)t(q)}{(2^n+1)t(q) -1}$, if $r >1$, then the $\infty$-slopes of $\sG(q,f,r)$ are
 $(1/r) \dfrac{(2^n+1)t(q)}{(2^n+1)t(q) -1} <1.$ It then follows from \cite[Proposition 4.2]{KT5} that
$H$ is equal to the Zariski closure $H_0^{{\mathrm{Zar}}}$ in $H$ of the normal closure $H_0$ in $H$ of the image of $I(0)$. 

\smallskip
We next show that for $r >1$, $H$ has no subgroup of index $2$. We argue by contradiction. Suppose
that $H$ has a subgroup $H_1$ of index $2$. Since every $H$-conjugate of the image of $I(0)$ has odd order $r$,
all such conjugates are contained in $H_1$, and thus $H_0 \leq H_1$. But $H_1$ is closed in $H$ by Proposition \ref{finite-index}(ii),
so $H=H_0^{{\mathrm{Zar}}} \leq H_1$, and hence $H=H_1$, a contradiction. We have shown that $H$ has no subgroup of index $2$.

\smallskip
(ii) We also know that $G$ is a normal subgroup of $H$ of index dividing $r$, simply because by the $r^{\mathrm {th}}$ power map, $\G_m/\overline{\F_2}$ becomes a finite \'etale Galois covering of itself with cyclic group of order $r$, and $G$ and $H$ are respectively the Zariski closure of the images of $\pi_1(\G_m/\overline{\F_2})$ and of a normal subgroup of cyclic index $r$ of that group. Now if $H=G$, then the existence of
$G_1$ leads to a contradiction by (i). If we now take $r$ to be a prime dividing $t(q)$, 
then $H/G \cong C_r$. 

Let $J$ and $Q$ denote the images of $I(\infty)$ and of $P(\infty)$ in $H$.  By Theorem \ref{Iinftygen}, $H$ is equal to the Zariski closure $H_{\infty}^{{\mathrm{Zar}}}$ in $H$ of the normal closure $H_\infty$ in $H$ of $J$. Suppose for the moment that $r \nmid |J|$. Then every $H$-conjugate of $J$ has order 
coprime to $r$, and so they are all contained in $G$, and thus $H_{\infty} \leq G$. But $G$ is closed in $H$, so
$H_\infty^{{\mathrm{Zar}}} \leq G$, and hence $H=G$, a contradiction.   

We have shown that $r$ divides $|J|$. Recall that $J = Q \rtimes C$, with $C$ a cyclic $2'$-group that permutes cyclically and 
transitively the $q-1$ simple $Q$-submodules $V_i$ in $V$, each of dimension $q_0=2^n$. As $Q$ is a $2$-group, $C$ is of order 
divisible by $r$. 

\smallskip
(iii) Since $n \neq 2$, by Theorem \ref{prime-suzuki1} the integer $t(q)=q-\sqrt{2q}+1$ admits a prime divisor 
\begin{equation}\label{desc10}
  r = \ppd(2,4(2n+1)).
\end{equation}  
At this point, we take for $r$ a $\ppd(2,4(2n+1))$ which divides $t(q)$.
Fix 
$c \in C$ of order $r$. Since $r|t(q)$, $r \nmid (q-1)$, and so $c$ stabilizes each of the subspace 
$V_i$, and certainly normalizes $Q$. For each $i$, let $\Phi_i$ denote the representation of $\langle Q,c \rangle$ on $V_i$. 
We claim that $\Phi_i(c)$ centralizes $\Phi_i(Q)$. 
Otherwise, $\Phi_i(c)$ acts faithfully on $\Phi_i(Q)$, of prime order $r$, and $\Phi_i(Q)$ is faithful and irreducible of degree 
$2^n$. Hence by Lemma \ref{action}, the order of $2$ modulo $r$ is at most $2n$, which contradicts \eqref{desc10}.

Therefore, for each $i$, $\Phi_i(c) = \alpha_i \cdot {\mathrm {Id}}_{V_i}$ for some $\alpha_i \in \C^\times$. But the cyclic group $C$ permutes the $V_i$'s transitively and $c \in C$, so $\alpha_1 = \ldots = \alpha_{q-1} =: \alpha$, i.e.\ $c$ acts on $V$ 
as $\alpha \cdot \mathrm{Id}_V$. As $|c|=r$, $\alpha \neq 1$ is a primitive $r^{\mathrm {th}}$ root of unity. Thus we may assume 
that 
$c \in J$ has trace $\dim(V) \cdot \zeta_r$. On the other hand, by Theorem \ref{ssandI(x)} and \eqref{airy-g1}, the trace of every element in $J$ belongs to
$\Q(i)$. Thus $\zeta_r \in \Q(i)$, a contradiction since $r$ is an odd prime. 

\smallskip
(iii) Since the irreducible representation $V$ of $G=G_\geo$ of $\sF(q,f)$ is not induced, for any $r|t(q)$, 
the representation of $G_\geo$ of $\sG(q,f,r)$ on $V$, which contains $G$, is not induced. 
\end{proof}

\section{Condition $\CSP$ and autoduality for Airy sheaves}
In this section, we continue to consider Airy sheaves $\sF(q,f)$ of the same general shape \eqref{airy-g}, but we consider them only geometrically, i.e.\ as lisse sheaves on $\A^1/\overline{\F_2}$. The key insight on which the results of this section are based is due to \v{S}uch, cf.\  \cite[Proposition 11.1]{Such}.
 
 \begin{lem} {\rm (\v{S}uch)}\label{suchlemma}
 Let $\sF$ be an Airy sheaf of rank $D \ge 2$ on $\A^1/\overline{\F_p}$. Let $\sH$ be a direct factor of $\End(\sF)$ of rank $r<D$. Then $\sH$ is a direct sum of $\sL_{\psi(ax)}$ for various $a \in \overline{\F_p}$.
 \end{lem}
 \begin{proof}The $\infty$-slopes of $\sF$ are all $1 + 1/D$. All slopes of $\End(\sF)$ are therefore $\le 1 + 1/D$. Thus all slopes of $\sH$ are $\le 1+1/D$. But each slope of $\sH$, written in lowest terms, has denominator $\le r$. Therefore $\sH$ has all slopes $\le 1$. Now take an irreducible constituent $\sK$ of $\sH$. Its Fourier transform, i.e.\ $t \mapsto H_c(A^1/\overline{\F_p},\sK\otimes \sL_{\psi(tx)})$, is perverse irreducible on $A^1/\overline{\F_p}$, so is either a single delta function $\delta_{a}$ or is the extension by direct image of a lisse sheaf on a dense open set. But on a dense open set of the $t$-line, $\sK\otimes \sL_{\psi(tx)}$ has all $\infty$-slopes $1$, so $H^2_c(A^1/\overline{\F_p},\sK\otimes \sL_{\psi(tx)})=0$ and by the Euler-Poincar\'{e} formula 
$\chi_c(A^1/\overline{\F_p},\sK\otimes \sL_{\psi(tx)}) =0$, so $\FT(\sK)$ is punctual, hence a single $\delta_{a}$. This means in turn that each irreducible constituent of $\sH$ is an  $\sL_{\psi(ax)}$. Because $\sF$ is irreducible, $\End(\sF)$ is completely reducible, hence $\sH$ is completely reducible, and so it 
is the sum of its irreducible constituents.
 \end{proof}

 \begin{thm}\label{indiffLpsi}Let $\sF$ be an Airy sheaf of rank at least $2$ on $\A^1/\overline{\F_p}$. Then the following conditions are equivalent.
 \begin{itemize}
 \item[\rm (i)] $\sF$ is geometrically induced.
 \item[\rm (ii)]$\End(\sF)$ contains a summand $\sL_{\psi(ax)}$ for some $a \neq 0$ in $\overline{\F_p}$.
 \item[\rm (iii)]There exists a geometric isomorphism $\sF \cong \sF \otimes \sL_{\psi(-ax)}$ for some $a \neq 0$ in $\overline{\F_p}$.
 \end{itemize}
 \end{thm}
 \begin{proof}
Let us denote by $n$ the rank of $\sF$. If $\sF$ is induced, then it is $g_\star \sH$, in which case 
$$\End(\sF)=(g_\star \sH)\otimes (g_\star \sH^\vee) \supseteq g_\star(\sH\otimes \sH^\vee) \supseteq g_\star\triv,$$
and hence $\End^0(\sF) \supseteq g_\star\triv/\triv$. But $g_\star\triv/\triv$ has rank less than $n$. So by Lemma \ref{suchlemma},
$g_\star\triv/\triv$ is the direct sum of various $\sL_{\psi(ax)}$. As $\End^0(\sF)$ does not contain the trivial sheaf (by irreducibility of $\sF$), we find that $\End^0(\sF)$, and hence  $\End(\sF)$, contains some $\sL_{\psi(ax)}$ for some nonzero $a \in \overline{\F_2}$.

Conversely, suppose $\End(\sF)$ contains a summand $\sL_{\psi(ax)}$ for some $a \neq 0$ in $\overline{\F_p}$. Thus there is a nonzero geometric homomorphism from $\sF$ to $\sF \otimes \sL_{\psi(-ax)}$. But source and target are geometrically isomorphic,
 so every nonzero homomorphism is an isomorphism. Thus we have a geometric isomorphism
 $$\sF \cong \sF\otimes \sL_{\psi(-ax)}.$$
 That $\sF$ is geometrically induced is the special case, where $G$ is $G_{\geo}$ for $\sF \oplus  \sL_{\psi(-ax)}$, $V$ is $\sF$, $L$ is $ \sL_{\psi(-ax)}$, $N$ is $p$, and $\E$ is $\overline{\Q_\ell}$, of 
the following general statement in characteristic zero representation theory.
 \end{proof}

\begin{thm}\label{inductioncriterion}
Over an algebraically closed field  $\E$ of characteristic zero, let
$V$ be a finite dimensional irreducible representation of dimension $d \ge 2$ of a group $G$, and $L$ a one-dimensional representation of $G$ which, viewed as a linear character $\chi$ of $G$, has finite order $N>1$. Suppose that $V \cong V\otimes L$ as representations of $G$. Denote by $G_0:=\Ker(\chi)$ the kernel of $\chi$, so that $G_0 \lhd G$ with $G/G_0$ cyclic of order $N$. Then $V$ is induced from a representation of a subgroup $H$ with $G_0 \leq H < G$.
\end{thm}
 
\begin{proof}
We first reduce to the case when $N$ is prime. If $V \cong V\otimes L$, then by induction $V \cong V\otimes L^{\otimes n}$ for every integer $n$. Let $r$ be a prime dividing $N$. Replacing $L$ by $L^{\otimes (N/r)}$ and  $\Ker(\chi)$ by the overgroup 
$K:=\Ker(\chi^{(N/r)})$, we are reduced to the case when $N=r$ is prime.

Let $U$ be a simple summand of $V|_{K}$, which is semisimple since $K \lhd G$.  
By Clifford theory, cf.\ \cite[(50.5)]{C-R}, $V$ is induced so long as $V|_{K}$ is not isotypic. Hence, if $V$ is not induced, then
$V|_{K} \cong eU:= \underbrace{U \oplus U \oplus \ldots \oplus U}_{e \mathrm{~times}}$ for some $e \in \Z_{\geq 1}$.
Note that $\Ind^G_K(\E) \cong \oplus^{r}_{i=1}L^{\otimes i}$, so by Frobenius reciprocity, cf.\ \cite[(10.20)]{CR2},
$$\Ind^G_K(V|_K) \cong \Ind^G_K\bigl((V|_K) \otimes \E\bigr) \cong V \otimes \Ind^G_K(\E) \cong \oplus^r_{i=1}(V \otimes L^{\otimes i}) \cong rV.$$ 
Again by Frobenius reciprocity \cite[(10.8)]{CR2},
$$e^2=\dim \Hom_K(V|_K,V|_K) = \dim \Hom_G\bigl(V,\Ind^G_{K}(V|_K)\bigr) = \dim \Hom_G(V,rV) =r,$$
a contradiction as $r$ is a prime.
\end{proof}

The following general result is stated for ease of later reference.
\begin{thm}\label{Limpliesinduced}Let $\sF$ be a lisse irreducible $\overline{\Q_\ell}$-sheaf of rank $d \ge 2$ on a smooth, geometrically connected $X/\overline{\F_p}$. Suppose that $\End^0(\sF)$ contains a  rank one summand $\sL$.
Then $\sF$ is induced.
\end{thm}
\begin{proof}Because $\sL$ occurs in $\End^0(\sF)$, it cannot be trivial (as by irreducibility of $\sF$, $\End(\sF)$ contains $\triv$ exactly once). Exactly as in the proof of the implication (ii) $\implies$ (iii) in Theorem \ref{indiffLpsi}, we infer that $\sF \cong \sF \otimes \sL^{-1}$. Taking determinants of this isomorphism, we find that $\sL^{\otimes d}$ is trivial. Thus $\sL$ has finite order $N >1$ as a character of $\pi_1(X)$. Now apply Theorem \ref{inductioncriterion}.
\end{proof} 

 \begin{thm}\label{generalprim}Let $k_0/\F_2$ be a finite extension, and $f(x) \in k_0[x]$ a polynomial of degree $(1+q_0)t(q)$. Then 
 the Airy sheaf $\sF(q,f)$ introduced in \eqref{airy-g}  
  is geometrically primitive, i.e.\ is not geometrically induced.
 \end{thm}

\begin{proof}
(i) In view of Theorem \ref{indiffLpsi}, it suffices to show that there is no geometric isomorphism from $\sF(q,f)$ to $\sF(q,f)\otimes \sL_{\psi(ax)}$ for any $a \neq 0$ in $\overline{\F_p}$. We argue by contradiction.
 Because $\sF(q,f)$ was geometrically a Fourier transform, so is $\sF(q,f)\otimes \sL_{\psi(ax)}$; the effect of tensoring 
 with  $\sL_{\psi(ax)}$ after $\FT$ is the same as translating additively by $-a$ before $\FT$:
$$\sF(q,f)\otimes \sL_{\psi(ax)} =\FT([x \mapsto x-a]^\star \sL(q,f)).$$
 Since $\FT$ is invertible, an isomorphism of  $\sF(q,f)$ with $\sF(q,f)\otimes \sL_{\psi(ax)}$ gives an isomorphism of $\sL(q,f)$ with its additive translate by $-a$:
 $$\sL_{\psi_2([x^{t(q)},f(x)])} \cong \sL_{\psi_2([(x-a)^{t(q)},f(x-a)])}.$$
 We will show no such isomorphism exists. 
 
 \smallskip
(ii)  Suppose that $n \ge 2$. In this case, we argue as follows. If two
 lisse rank one sheaves are isomorphic, then so are their tensor squares.
 
 Quite generally, addition of Witt vector of length $2$ over an $\F_2$-algebra is given by
 $$[a,b]+[A,B]=[a+A, B+b +aA], \ [a,b]+[a,b]=[0,a^2].$$
 Thus
$$\sL_{\psi_2([g(x),f(x)])}^{\otimes 2} =\sL_{\psi_2([g(x),f(x)]+[g(x),f(x)])}=\sL_{\psi_2([0,g(x)^2])}=\sL_{\psi(g(x)^2)} \cong \sL_{\psi(g(x))}.$$
 So we would have a geometric isomorphism
 $$ \sL_{\psi(x^{t(q)}))} \cong  \sL_{\psi((x-a)^{t(q)}))}, $$
 or equivalently a geometric isomorphism
 $$ \sL_{\psi((x-a)^{t(q)}-x^{t(q)})} \cong \overline{\Q_\ell}.$$
 Now we use the explicit shape of $t(q)=q+1-2q_0 =2q_0^2 -2q_0+1=(q_0-1)(2q_0)+1$.

The key point is that $t(q)=1+dQ$ with $d \ge  2$ prime to $p=2$ (here $d=q_0-1)$ and $Q$ a strictly positive power of  $p=2$ (here $Q=2q_0$). Then
 $$\begin{aligned}(x-a)^{1+qD} & =(x-a)(x^Q-a^Q)^d\\
    & =(x-a)\bigl(x^{dQ} -da^Qx^{(d-1)Q} +({\rm terms \ of \ degree \ }\le (d-2)Q)\bigr)\\
    & =x^{1+dQ} -da^Qx^{1+(d-1)Q} +({\rm terms \ of \ degree \ }\le 1+(d-2)Q)\\
    &\ \  + ({\rm polynomial \  in \ }x^Q{\rm \ of\ degree\ }d).
    \end{aligned}$$
 Thus up to Artin-Schreier equivalence, 
 $$(x-a)^{1+dQ} -x^{1+dQ} = -da^Qx^{1+(d-1)Q} +({\rm terms\ of \ degree \ } \le 1+(d-2)Q) +({\rm a\ term\ of \ degree\  }d).$$
 Thus $\sL_{\psi((x-a)^{1+dQ} -x^{1+dQ} )}$ has $\Swan=1+(d-1)Q >0$, so is not geometrically trivial.
 
  
\smallskip  
(iii)  We now turn to the case $n=1$, which requires a more delicate analysis. What makes the $n \ge 2$ case argument work is that $q_0-1 =2^n-1$ has $n \ge 2$ binary digits. In the $n=1$ case, $t(8)=5$, so the above argument would involve examining
  $$(x-a)^5 -x^5 = (x-a)(x^4-a^4) -x^5 =ax^4-a^4x +a^5,$$
but this is  Artin-Schreier equivalent to $(a^{1/4}-a^4)x$, which for $a \in \mu_{15}$ is Artin-Schreier trivial.
  
  Thus we must look instead at 
  $$ \sL_{\psi_2([(x-a)^{5},f(x-a)])}\otimes \sL_{\psi_2([x^{5},f(x)])}^{-1}= \sL_{\psi_2([(x-a)^{5},f(x-a)]-[x^{5},f(x)])}.$$
  Here $-[a,b]=[a,b-a^2]$, implying $-[x^{5},f(x)]=[x^5,f(x)-x^{10}]$, so this is
  $$ \sL_{\psi_2([(x-a)^{5}+x^{5},f(x-a)+f(x) -x^{10} +x^5(x-a)^5])}.$$
  But $f(x)$ has degree $15$, say $f(x)=bx^{15} + cx^{14} +dx^{13} + \mbox{lower terms}$
  with $b \neq 0$.
So under Artin-Schreier equivalence 
  $$f(x) \equiv bx^{15} + dx^{13 }+ \mbox{lower terms}.$$
  Similarly, under Artin-Schreier equivalence 
  $$f(x-a)\equiv b(x-a)x^{15} + d(x-a)^{13} + \mbox{lower terms}.$$
  Thus under Artin-Schreier equivalence 
  $$f(x-a)+f(x) -x^{10} +x^5(x-a)^5 \equiv b(x-a)^{15}-bx^{15} + \mbox{lower terms}.$$
  The difference $$b(x-a)^{15}-bx^{15} =-abx^{14} +ba^2x^{13}+ \mbox{lower terms}$$
  is thus Artin-Schreier equivalent to
  $$ba^2x^{13}+ \mbox{lower terms}.$$
  Now view $\sL_{\psi_2([(x-a)^{5}+x^{5},f(x-a)+f(x) -x^{10} +x^5(x-a)^5])}$ as the tensor product (using $ [a,b]=[a,0]+[0,b]$)
  $$\sL_{\psi_2([(x-a)^{5}+x^{5},0]}\otimes \sL_{\psi(f(x-a)+f(x) -x^{10} +x^5(x-a)^5)}.$$
  The first factor has $\Swan_\infty \le 2 \times 4$, while the second factor, which is $\sL_{\psi({\rm a \ polynomial\  of \  degree\ }13)}$, has
  $\Swan_\infty=13$, and hence their tensor product has $\Swan_\infty=13$, so is not geometrically trivial.
 \end{proof}
  
Recall that a lisse sheaf $\sH$ is {\it Lie-irreducible} if, in the underlying representation of its $G_{\geo}$, the identity component  $G_{\geo}^\circ$
  acts irreducibly.  It is {\it Lie-self-dual} if the  given representation of $G_{\geo}$, when restricted to $G_{\geo}^\circ$, is self-dual.
  
\begin{thm}\label{newauto}
Suppose that a lisse sheaf $\sF(q,f)$ as in \eqref{airy-g} is Lie-irreducible and Lie-self-dual. Then $\sF(q,f)$ is self-dual.
 \end{thm}
  
  \begin{proof}The two sheaves $\sF(q,f)$ and its dual $\sF(q,f)^\vee$ are irreducible representations of  $G_{\geo}$ whose restrictions to the identity component  $G_{\geo}^\circ$ are isomorphic. So the two, as representations of  $G_{\geo}$, differ by a linear character of   $G_{\geo}/G_{\geo}^\circ$.
  Thus 
  $$\sF(q,f)^\vee \cong \sF(q,f)\otimes \sL$$ 
  for some lisse, rank one $\sL$ on $\A^1$. Both $\sF(q,f)^\vee $ and $\sF(q,f)$
  have all $\infty$-slopes $1 +1/\rank(\sF(q,f)) <2$. Therefore $\Swan_\infty(\sL) \le 1$ (otherwise $\sF(q,f)^\vee$ would have all  $\infty$-slopes $\ge 2$). 
  
  If $\Swan_\infty(\sL) =0$, then $\sL$ is lisse on $\A^1$ and tame at $\infty$, so geometrically trivial, and $\sF(q,f)$ is geometrically self-dual.
  If $\Swan_\infty(\sL) =1$, then $\sL$ is $\sL_{\psi(ax)}$ for some nonzero $a \in \overline{\F_p}$. We will show that this case cannot arise.
  
  Recall that $\sF(q,f) :=\FT_{\psi}(\sL_{\psi_2([x^{t(q)},f(x)])})$. In general, the interaction of $\FT$ with duality is a geometric isomorphism
$$D(\FT_\psi(\sH)) \cong \FT_{\overline{\psi}}(D\sH)).$$
In characteristic $2$, where $\psi$ takes values $\pm 1$, we have $\overline{\psi}=\psi$. Thus 
$$\sF(q,f)^\vee \cong \FT_\psi(\sL_{\psi_2(-[x^{t(q)},f(x)])})=\FT_\psi(\sL_{\psi_2([x^{t(q)},f(x)+(x^{t(q)})^2])} )$$
while 
$$\sF(q,f)\otimes \sL_{\psi(ax)} \cong \FT_\psi([x \mapsto x-a]^\star \sL_{\psi_2([x^{t(q)},f(x)])}=\FT_\psi( \sL_{\psi_2([(x-a)^{t(q)},f(x-a)])}).$$
By Fourier inversion, this is equivalent to a geometric isomorphism
$$\sL_{\psi_2([x^{t(q)},f(x)+(x^{t(q)})^2])} \cong \sL_{\psi_2([(x-a)^{t(q)},f(x-a)])}.$$

We first treat the case $n \ge 2$. Already the tensor squares of these two lisse rank one sheaves are not geometrically isomorphic, by the identical argument used to treat the case $n \ge 2$ in the proof of Theorem \ref{generalprim}.

In the case $n=1$, we must show that the lisse rank one sheaf
$$\sL_{\psi_2([(x-a)^{s},f(x-a)])} \otimes (\sL_{\psi_2([x^{5},f(x)+(x^{5})^2])} )^{-1}$$
is not geometrically trivial. The second tensor factor is
$$(\sL_{\psi_2([x^{5)},f(x)+(x^{5})^2])} )^{-1}=\sL_{\psi_2(-[x^{5},f(x)+(x^{5})^2])} =\sL_{\psi_2([x^{5},f(x)])}.$$
So we must show that 
$$\begin{aligned}\sL_{\psi_2([(x-a)^{5},f(x-a)])} \otimes \sL_{\psi_2([x^{5)},f(x)])}
   &  = \sL_{\psi_2([(x-a)^{5},f(x-a)] +[x^{5},f(x)])}\\
   & =\sL_{\psi_2([(x-a)^{5} +x^{5}, f(x-a)+f(x) +x^{5}(x-a)^5])}\\
   & =\sL_{\psi_2([(x-a)^{5} +x^{5},0])}\otimes \sL_{\psi(f(x-a)+f(x) +x^{5}(x-a)^{5})}\end{aligned}$$
is not geometrically trivial. Exactly as in the proof of the $n=1$ case of Theorem \ref{generalprim}, the first factor has $\Swan_\infty \le 2\times 4=8$,
while the second factor has $\Swan_\infty=13$.
    \end{proof}

\begin{thm}\label{notselfdual}No sheaf $\sF(q,f)$ introduced in \eqref{airy-g} is geometrically self-dual. 
\end{thm}
\begin{proof}Recall that $\sF(q,f) :=\FT_{\psi}\bigl(\sL_{\psi_2([x^{t(q)},f(x)])}\bigr)$. In general, the interaction of the Fourier transform $\FT$ with the duality functor $D(\cdot)$ is a geometric isomorphism
$$D(\FT_\psi(\sH)) \cong \FT_{\overline{\psi}}(D\sH).$$
In characteristic $2$, where $\psi$ takes values $\pm 1$, we have $\overline{\psi}=\psi$. Therefore $\sF(q,f) $ is self-dual if and only if 
$\sL_{\psi_2([x^{t(q)},f(x)])}$ is self-dual. Its dual is $\sL_{\psi_2(-[x^{t(q)},f(x)])}$.
We have the Witt vector addition law $[x,y]=[x,0]+[0,y]$, and, for Witt vectors over $\F_2$-algebras, $-[x,y]=[x,y+x^2]$. 
So the dual of $\sL_{\psi_2([x^{t(q)},f(x)])}$ is $\sL_{\psi_2([x^{t(q)},f(x)+(x^{t(q)})^2])}$, and 
the asserted isomorphism is
$$\sL_{\psi_2([x^{t(q)},0)])}\otimes \sL_{\psi(f(x))} \cong \sL_{\psi_2([x^{t(q)},0)])}\otimes \sL_{\psi(f(x)+(x^{t(q)})^2)}.$$
This holds if and only if  there is an isomorphism
$$\sL_{\psi(f(x))} \cong  \sL_{\psi(f(x)+(x^{t(q)})^2)},$$
or equivalently if $\sL_{\psi((x^{t(q)})^2)}$ is geometrically trivial. By the Artin-Schreier reduction we have 
$\sL_{\psi((x^{t(q)})^2)} \cong  \sL_{\psi(x^{t(q)})}$. The latter sheaf
is not geometrically trivial, because $x^{t(q)}$ is a polynomial of odd degree $t(q)$ and so the sheaf has $\Swan_\infty(\sL_{\psi(x^{t(q)})})=t(q) \neq 0$.
\end{proof}

\begin{cor}\label{cor:sd}
Suppose that the sheaf $\sF(q,f)$ introduced in \eqref{airy-g} is Lie-irreducible. Then it is not Lie-self-dual.
\end{cor}
\begin{proof}Combine Theorem \ref{newauto} and Theorem \ref{notselfdual}.
\end{proof}

\begin{thm}\label{Airytensorindec}
The sheaf $\sF(q,f)$ introduced in \eqref{airy-g} is tensor indecomposable.
\end{thm}

\begin{proof}We will show that $\sF(q,f)$ is tensor indecomposable as a representation of $\pi_1:=\pi_1(\A^1/\overline{\F_p})$. Because $\pi_1$ has cohomological dimension $\le 1$ (this is true for the $\pi_1$ of any smooth, connected affine curve over an algebraically closed field), the argument of \cite[Corollary 10.4]{Ka-RL-T-2Co1} shows that if $\sF(q,f)$ is tensor decomposable, then it is linearly tensor decomposable, i.e.\ we have an isomorphism of local systems on $\A^1/\overline{\F_p}$,
$$\sF(q,f) \cong \sA \otimes \sB,$$
with both $\rank(\sA), \rank(\sB) \ge 2$. To fix ideas, suppose  $\rank(\sA) \le \rank(\sB)$. Then
$$\End(\sF(q,f) ) =\End(\sA)\otimes \End(\sB)= (\triv + \End^0(\sA))\otimes  (\triv + \End^0(\sB)) {\rm \  contains \ }  \triv + \End^0(\sA).$$
Thus $\End^0(\sF(q,f))$ contains $ \End^0(\sA)$ as a direct factor. Now $ \End^0(\sA)$ has rank less than 
$$(\rank(\sA))^2 \le \rank(\sA)\rank(\sB)=\rank(\sF(q,f)).$$ 
So by Lemma \ref{suchlemma}, $\End^0(\sF(q,f))$ contains some $\sL_{\psi(ax)}$ with $a \neq 0$ ($a \neq 0$ 
because  $\End^0(\sF(q,f))$ only contains $\triv$ once, by irreducibility of $\sF$). The proof of Theorem \ref{generalprim} shows this is impossible.
\end{proof}
  

\begin{lem}For every integer $n \ge 2$, the integer $2^n-1$ is never a perfect power $x^m$ with $x \in \Z$ and $m \ge 2$. 
\end{lem}
\begin{proof}We argue by contradiction. If $2^n-1=x^m$, then $x$ is odd and $x^m =2^n-1 \equiv 3 \pmod*{4}$. Thus $m$ is odd, and
hence $2^n =x^m+1$ is divisible by $x+1$. The quotient $\frac{x^m+1}{x+1} >1$ is the alternating sum of $m$ powers of  the odd integer $x$, so is itself odd. Thus $\frac{x^m+1}{x+1}$ is an odd divisor of $2^n$, the desired contradiction.
\end{proof}
\begin{cor}\label{nottensorind}For $n \ge 1$, no lisse sheaf of rank $2^n(2^{2n+1}-1)$ can be tensor induced. In particular, $\sF(q,f)$ is not tensor induced.
\end{cor}

\begin{thm}\label{cond-S}
The local systems $\sF(q,f)$ on $\A^1/\overline{\F_2}$ introduced in \eqref{airy-g} all satisfy the condition
$\CSP$ of \cite[Definition 1.2]{KT5}.
\end{thm}
\begin{proof}
First, by \cite[Proposition 7.4]{Such}  the underlying representation $V$ of $G_\geo$ of $\sF(q,f)$ is irreducible.
That $\sF(q,f)$ is primitive, tensor indecomposable, and not tensor induced is the content of Theorem \ref{generalprim}, Theorem \ref{Airytensorindec} and Corollary \ref{nottensorind}. That $\det(\sF(q,f))$ has finite order results from the fact that $\sF(q,f)$ began life over a finite subfield of $\overline{\F_p}$ (in fact any subfield containing the coefficients of $f$).
\end{proof}

\begin{thm}\label{comptsize}
For the local system $\sF(q,f)$ on $\A^1/\overline{\F_2}$ introduced in \eqref{airy-g}, every irreducible constituent of $\End^0(\sF(q,f))$ has dimension $\ge \rank(\sF(q,f))$. In particular, if $G_{\geo, \sF(q,f)}$ is not finite, then $G_{\geo, \sF(q,f)}^\circ$ is a simple algebraic group of dimension $\ge \rank(\sF(q,f))$.
\end{thm}
\begin{proof}
By Lemma \ref{suchlemma}, any irreducible constituent of dimension $< D:=\rank(\sF(q,f))$ is a single $\sL_{\psi(ax)}$, while Theorem \ref{generalprim} shows that $\End^0(\sF(q,f))$ contains no $\sL_{\psi(ax)}$. Because $\sF(q,f)$ satisfies condition
$\CSP$ by Theorem \ref{cond-S}, if $G_{\geo, \sF(q,f))}$ is not finite, then its identity component  $G_{\geo, \sF(q,f))}^\circ$ is a simple algebraic group. In that case,
$\Lie(G_{\geo, \sF(q,f))}^\circ)$ is an irreducible constituent of $\End^0(\sF(q,f))$, so has dimension $\ge D$.
\end{proof}


\begin{thm}\label{tensorindecandprim}
Consider the sheaf $\sF(q,f)$ in \eqref{airy-g} subject to the condition 
\eqref{airy-s}. Then, for each $r|t(q)$, the descent $\sG(q,f,r)$ satisfies the condition
$\CSP$ of \cite[Definition 1.2]{KT5}.
\end{thm}

\begin{proof}
Because $G_{\geo,\sF(q,f)}$ is a subgroup of $G_{\geo,\sG(q,f,r)}$ of finite index, the fact that  $\sF(q,f)$ satisfies condition
$\CSP$, see Theorem \ref{cond-S}, implies that $\sG(q,f,r)$ does as well. [Condition $\CSP$ holds for $(G,V)$ if it holds for $(H,V)$ with $H$ a subgroup of $G$ of finite index.]
\end{proof}

We also record the following

\begin{thm}\label{Lieselfdual}
Under the condition \eqref{airy-s}, if the sheaf $\sF(q,f)$ in \eqref{airy-g} is Lie-irreducible, then it is not Lie-self-dual.
\end{thm}

\begin{proof}
This is just a restatement of Corollary \ref{cor:sd},
under the more restrictive hypotheses of \eqref{airy-s}.
Given Theorem \ref{notselfdual},
it suffices to show that if $\sF(q,f)$ is Lie-irreducible and Lie-self-dual, then it is self-dual.

However, there is a simpler proof of this last fact, using the descent $\sG(q,f):=\sG(q,f,t(q))$ in Theorem \ref{tensorindecandprim}. 
Exactly as in the proof of Theorem \ref{newauto}, we find an isomorphism
$$\sG(q,f)^\vee \cong \sG(q,f)\otimes \sL$$
for some lisse $\sL$ on $\G_m$ of rank one. Because $\sG(q,f)^\vee $and  $\sG(q,f)$ are both tame at $0$ and with all $\infty$-slopes $<1$,
$\sL$ is tame on $\G_m$, hence a Kummer sheaf $\sL_\chi$. Because  $\sG(q,f)^\vee $and  $\sG(q,f)$ both have $I(0)$ representations which are sums of characters of order dividing $t(q)$, $\chi$ is a ratio of characters of order  dividing $t(q)$, so $\chi$ has order dividing $t(q)$.
Pulling back this isomorphism by $t(q)^{\mathrm {th}}$ power, we get an isomorphism
$\sF(q,f)^\vee \cong \sF(q,f)$.
\end{proof}

\section{A local system for the Suzuki group $\tw2B_2(8)$}
Now we can prove the first main result of the paper, which establishes \cite[Conjecture 2.2]{Ka-ERS} in the case $q=8$:

\begin{thm}\label{main1}
Let $q=8$. Both the local systems $\sF_q$ and $\sG_q$ have geometric monodromy group $G_\geo \cong \tw2 B_2(8)$ in
one of its irreducible representation of degree $14$. Over $\F_2$, the local systems $\sF_q$ and $\sG_q$ have arithmetic monodromy group $G_\ari \cong \Aut(\tw2 B_2(8))$.
\end{thm}

\begin{proof}
(a) Let $G$, respectively $H$, denote the geometric monodromy group of $\sF_q$, respectively of $\sG_q$. 
Similarly, let $G_{\ari}$, respectively $H_{\ari}$, denote the arithmetic monodromy group of $\sF_q$, respectively of $\sG_q$, over $\F_2$.
We will use the fact that $f:=\Frob_{1,\F_2}$ has order $15$ and trivial determinant. Indeed, a {\sc Magma} calculation shows that
$$\Trace(\Frob_{1,\F_{2^{15}}}|\sG_q)=\Trace(\Frob_{1,\F_{2^{15}}}|\sF_q)=14.$$
By Theorem \ref{ssandI(x)}, $\Frob_{1,\F_{2^{15}}}|\sG_q$ is semisimple, hence (being of weight zero in a $14$-dimensional representation) is the identity. Therefore $\det(\Frob_{1,\F_{2^{15}}}|\sG_q))$ is a root of unity of order dividing both $4$ (by Theorem \ref{dets}) and $15$, so this determinant is  trivial. But $H_{\ari} = \langle f,H \rangle$, and $H$ has trivial determinant (by Theorem \ref{dets}), and hence 
 $H_{\ari}$ also has trivial determinant. As $G \lhd H$ and $G_{\ari} \leq H_{\ari} $ are subgroups, both $G$ and $G_{\ari}$ have trivial determinants.
 
Recall that $G$ is normal in $H$ of index dividing $t(q)=5$.
By Theorem \ref{cond-S}, 
$H$ satisfies $\CSP$. Since the rank of the sheaves is $14$, not a prime power, by \cite[Lemma 1.1]{KT5} this implies that 
either the identity component $H^\circ$ of $H$ is a simple algebraic group that acts irreducibly on the underlying representation $V$ of
$H$, or $H$ is a finite almost quasisimple group with $L:=H^{(\infty)}$ also acting irreducibly on $V$.

\smallskip
(b) Consider the former case. Note that if $H^\circ$ is classical of rank $r$, then either it is of type $A$ and $r \leq 13$, or 
$r \leq 7$ by \cite[Proposition 5.4.11]{KlL}. Using the tables of \cite{Lu}, we see that $H^\circ$ is of type $\SL_2$, $\SL_{14}$, $\Sp_4$, $\Sp_6$, $\Sp_{14}$, 
$\SO_{14}$, or $G_2$. Since $[H:G] < \infty$, $G^\circ=H^\circ$.

Suppose first that $H^\circ = \SL_{14}(\C)$. In this case, $V$ is just the natural module for $H^\circ$, hence $M_{2,2}(V)$ takes the 
smallest possible value $2$ for all $H^\circ$, $H$, $G^\circ$, and $G$. But this contradicts Corollary \ref{m4}.

Next suppose that $G^\circ=H^\circ$ is of type $\SL_2$, $\Sp_4$, $\Sp_6$, $\Sp_{14}$, $\SO_{14}$, or $G_2$. In all these cases, 
the $G^\circ$-module $V$ is self-dual, and this contradicts Corollary \ref{cor:sd}.
%

\smallskip
(c) We have shown that $H$ is finite, and we are in the almost quasisimple case. In this case, $H_{\ari}$ is also finite, and
$H_{\ari}^{(\infty)} = H^{(\infty)} = G_{\ari}^{(\infty)} = G^{(\infty)} = L$. Let $\varphi$ denote the character of $H_{\ari}$ in the 
underlying representation $V$. Then the formula for the trace function of 
$\sG_q$ and the existence of an element with trace
$2\zeta_4$ show that 
\begin{equation}\label{c10}
  \Q(\varphi)=\Q(\zeta_4). 
\end{equation}
Since any element $z$ of $\CB_{H_{\ari}}(L)=\ZB(H_{\ari})$ acts as a root of unity 
$\gamma$ on $V$, this implies that $\gamma^4=1$. Since $H_{\ari}$ has trivial determinant, $\gamma^{14}=1$
and thus $\gamma = \pm 1$. It follows that 
\begin{equation}\label{c11}
  |\CB_{H_{\ari}}(L)| = |\ZB(H_{\ari})| \leq 2.
\end{equation}  
Now, using \eqref{c10} and the fact that $L$ acts irreducibly on $V = \C^{14}$, the classification 
results of \cite{HM} show that $L = \PSL_2(13)$, $\SL_2(13)$, $\ABS_7$, $\ABS_8$, $\SU_3(3)$, $G_2(3)$, $2\cdot J_2$, 
$\ABS_{15}$, or $\tw2 B_2(8)$.

In all but the last case, $\varphi|_L$ is real-valued; furthermore, $|\Out(L)| \leq 2$. As shown in the proof of Theorem \ref{prim}, 
$G$ has no subgroups of
index $2$, so $G$ can induce only inner automorphisms of $L$. But $G_{\ari}=\langle G,f\rangle$, implying
$G_{\ari}/G \inj C_{15}$; so the same conclusion holds for
$G_{\ari}$, and hence $G_{\ari} = \CB_{G_{\ari}}(L)L = \ZB(G_{\ari})L$, with $|\ZB(G_{\ari}| \leq 2$ by \eqref{c11}. 
On the other hand, as 
we saw in the computation for the proof of Corollary \ref{m4}, some $g \in G_{\ari}$ has trace $2\zeta_4$ on $V$.
Write $g=zh$ with $z \in \ZB(G_{\ari})$ and $h \in L$. It follows
that $2\zeta_4 = \varphi(g) = \pm \varphi(h)$, a contradiction since $\varphi(h) \in \R$.

We have therefore shown that $L = \tw2 B_2(8)$. As $|\Out(L)| = 3$, \eqref{c11} implies that $|H/L|$ divides $6$. On the other hand,
$H$ is the normal closure of the image of order $5$ of $I(0)$, so we conclude that $H=L$. Since $L \lhd G \lhd H$, we also 
get $G=L$.

Recalling again that $G_{\ari}/G \inj C_{15}$, we now see from \eqref{c11} that $\ZB(G_{\ari}) = 1$ and thus 
$G_{\ari} \inj \Aut(L)  \cong L \rtimes C_3$. A {\sc Magma} computation shows the existence of a Frobenius element with trace $\zeta_4$, and this implies that $G_{\ari} > L$,
whence $G_{\ari} = \Aut(L)$. Since $H_{\ari}$ is generated over $H$ by $f=\Frob_{1,\F_2}$, an element of order
$15$, $\ZB(H_{\ari})=1$ by \eqref{c11}. Thus $H_{\ari} \inj \Aut(L)$, and so $H_{\ari} = \Aut(L)$ as 
$H_{\ari} \geq G_{\ari}$.
\end{proof}

\section{Low-dimensional representations of classical groups}\label{low-dim}
In this section, we will extend the classification results obtained in \cite[Proposition 5.4.11]{KlL} and \cite[Theorem 5.1]{Lu}.
Even though the intended applications in the paper only need the complex case of these results, we establish them in the modular
case, which is interesting in its own right.
  
Let $\F$ be an algebraically closed field of characteristic $p \geq 0$ and 
let $G$ be a simple, simply connected, classical algebraic group of rank $r$ over 
$\F$. Fixing a maximal torus in $G$, we consider 
the set of simple roots $\{ \al_1, \ldots, \al_r\}$ and the corresponding
set of fundamental weights $\{ \om_1, \ldots ,\om_r\}$ (in the
ordering of \cite{OV}). Then the set 
$$\LP = \left\{ \sum^{r}_{i=1}a_i \om_i \mid a_i \in \Z,a_i \geq 0
  \right\}$$
of dominant weights admits the 
partial ordering $\succ$ where $\lam \succ \mu$ precisely when 
$\lam - \mu = \sum^{r}_{i=1}k_i\al_i$ for some non-negative integers $k_i$.
As usual, $W$ denotes the Weyl group. If 
$\lam \in \LP$, let $L(\lam)$ denote the irreducible $\F G$-module
with highest weight $\lam$. 

We will rely on the following two results.

\begin{thm}\label{premet} {\rm \cite{Pr}}
Let $G$ be a simple, simply connected algebraic group in 
characteristic $p > 0$. If the root system of $G$ has different root 
lengths, then we assume that $p \neq 2$, and if $G$ is of 
type $G_2$, then we 
also assume that $p \neq 3$. Let $\lam$ be a restricted dominant
weight. Then the set of weights $\Pi(\lam)$ of the irreducible $G$-module 
$L(\lam)$ is the union of the $W$-orbits of 
dominant weights $\mu$ with $\lam \succ \mu$.
\hfill $\Box$
\end{thm}

\begin{lem}\label{orb} {\rm \cite[Lemma 10.3B]{H}}
Let $\lam = \sum^r_{i=1}a_i\om_i$ be a dominant weight. Then the 
stabilizer of $\lam$ in the Weyl group is the Young subgroup 
generated by the reflections $\rho_i$ along the simple roots $\al_i$ for which 
$a_i = 0$.
\hfill $\Box$
\end{lem}

Our first result treats groups of type $A$ and includes a strengthening for $\SL_{22}$:

\begin{thm}\label{class-A}
Let $G = \SL_n(\F)$ with $n = r+1 \geq 8$, and let $L(\lam)$ be an  irreducible $\F G$-representation, which is 
restricted if $p = \Char(\F) > 0$. 
Suppose that $\dim L(\lam) \leq M$, where $M:=\binom{n}{4}$ if $n \neq 22$ and $M := 8176$ if $n = 22$.
Then $\lam = 0$, $a\om_1$ or $a\om_r$ with $1 \leq a \leq 3$, $\om_1+\om_r$, $\om_2$ or $\om_{r-1}$, 
$\om_3$ or $\om_{r-2}$, $\om_4$ or $\om_{r-3}$, $\om_1+\om_2$ or
$\om_{r-1}+\om_r$, $2\om_1+\om_r$ or $\om_1+2\om_r$, and $\om_2+\om_r$ or $\om_1+\om_{r-1}$.
\end{thm}

\begin{proof}
Write the highest weight $\lam$ as $\sum^r_{i=1}a_i\om_i$ with $a_i \in \Z_{\geq 0}$.

\smallskip
(a) First suppose that there is some weight $\mu=\sum^r_{i=1}b_i\om_i \in \LP$ with $\lam \succ \mu$ and $b_j \neq 0$ for 
some $s+1 \leq j \leq r-s$, where 
$$s := 3 \mbox{ if }n \neq 22 \mbox{ and }s:=4 \mbox{ if }n=22.$$ 
By Theorem \ref{premet}, $\Pi(\lam)$ contains the $W$-orbit $\sO(\mu)$ of $\mu$. By Lemma \ref{orb},
$$\Stab_W(\mu) \leq \langle \rho_1, \ldots,\rho_{j-1}, \rho_{j+1}, \ldots,\rho_r \rangle = \SBS_{j} \times \SBS_{n-j}$$
(where $\rho_i = (i,i+1) \in W = \SBS_n$). Hence 
$$\dim L(\lam) \geq |\sO(\mu)| \geq [\SBS_n:(\SBS_j \times \SBS_{n-j})] = \binom{n}{j} \geq \binom{n}{s+1} \geq M,$$
a contradiction if $n = 22$, or if $5 \leq j \leq r-4$. Suppose $j \in \{4,r-3\}$ and $n \neq 22$. Then 
$\dim L(\lam) = |\sO(\mu)| = \binom{n}{4}$, showing $\mu$ is the unique dominant weight of $L(\lam)$, whence $\mu=\lam$.
This also forces $\Stab_W(\lam) =  \SBS_{j} \times \SBS_{n-j}$, and so $\lam = a_j\om_j$. Now if $a \geq 2$, then 
$\lam \succ \lam-\al_j =\om_{j-1}+(a_j-2)\om_j+\om_{j+1}$, and the latter is another dominant weight of $L(\lam)$, a contradiction.
So $\lam=\om_4$ or $\om_{r-3}$ in such a case.

We may therefore assume that
\begin{equation}\label{a10}
  a_j = 0 \mbox{ for every }s+1 \leq j \leq r-s.
\end{equation}  
Next suppose that $\sum^s_{i=1}ia_i \geq s+1$. By \cite[Lemma 2.6]{GLT} (applied with $m=s$),  there is some weight 
$\mu=\sum^r_{i=1}b_i\om_i \in \LP$ with $\lam \succ \mu$ and $b_{s+1} > a_{s+1}$. This situation 
is already considered by the preceding analysis. 
Using the symmetry under the graph automorphism $\tau$ of $G$, 
\begin{equation}\label{a11}
  \sum^s_{i=1}ia_i \leq s,~\sum^s_{i=1}ia_{r+1-i} \leq s.  
\end{equation}

\smallskip
(b) Suppose $a_4 > 0$ or $a_{r-3} > 0$. By symmetry, we may assume $a_4> 0$. By \eqref{a11} $n=22$,
$a_4=1$, and $a_1=a_2=a_3=0$. Now if $\lam \neq \om_4$, then $a_j \geq 1$ for some $r-3 \leq j \leq r$ by \eqref{a10}. 
By Lemma \ref{orb}, 
$$\Stab_W(\lam) \leq \langle \rho_1,\rho_2,\rho_3, \rho_5, \rho_6, \ldots,\rho_{j-1},\rho_{j+1}, \ldots,\rho_r \rangle 
   \cong \SBS_4 \times Y_1,$$
with $Y_1$ a proper Young subgroup of $\SBS_{n-4}$ and hence $|Y_1| \leq (n-5)!$. It follows that 
$$\dim L(\lam) \geq [W:\Stab_W(\lam)] \geq (n-4)\binom{n}{4} > M,$$ 
a contradiction. Hence from now on we may assume that $a_4=a_{r-3}=0$.

\smallskip
(c) Suppose that $a_3 > 0$ or $a_{r-2} > 0$. By symmetry, we may assume $a_3 > 0$, and so $a_3=1$ by \eqref{a11}. 
Suppose $a_j \geq 1$ for some $r-3 \leq j \leq r$. By Lemma \ref{orb}, 
\begin{equation}\label{a11b}
  \Stab_W(\lam) \leq \langle \rho_1,\rho_2,\rho_4, \rho_5, \ldots,\rho_{j-1},\rho_{j+1}, \ldots,\rho_r \rangle 
   \cong \SBS_3 \times Y_2,
\end{equation}   
with $Y_2$ a proper Young subgroup of $\SBS_{n-3}$ and hence $|Y_2| \leq (n-4)!$. It follows that 
$$\dim L(\lam) \geq [W:\Stab_W(\lam)] \geq (n-3)\binom{n}{3} > M$$ 
(as $n \geq 8$), a contradiction. Now, if $n \neq 22$, then using \eqref{a10} and \eqref{a11} we see that $\lam=\om_3$.
If $n = 22$ but $\lam \neq \om_3$, then $\lam = \om_1+\om_3$. In this case, $\lam-(\al_1+\al_2+\al_3) = \om_4$ is 
also a weight of $L(\lam)$ by Theorem \ref{premet}, and since 
$$\Stab_W(\om_1+\om_3) = \langle \rho_2,\rho_4, \rho_5, \ldots,\rho_r \rangle 
   \cong \SBS_2 \times \SBS_{n-3},~\Stab_W(\om_4) = \SBS_4 \times \SBS_{n-4}$$
by Lemma \ref{orb}, 
\begin{equation}\label{a12}
   |\sO(\om_1+\om_3)|+|\sO(\om_4)| = \binom{22}{3}/2+\binom{22}{4} = 11935 > M,
\end{equation}
which contradicts $\dim L(\lam) < M$.

\smallskip
(d) We may now assume that $a_i=0$ for $3 \leq i \leq r-2$. Suppose that $a_2 > 0$ or $a_{r-1} > 0$. By symmetry, we may assume 
$a_2 > 0$. If $a_2 \geq 2$, then by \eqref{a11} $n=22$, $(a_2,a_1)= (2,0)$. Note that 
$\lam \succ \lam-\al_2 = \om_1+\om_3 +a_{r-1}\om_{r-1}+a_r\om_r$. It follows from Theorem \ref{premet} 
that $\Pi(\lam)$ contains $\om_1+\om_3 +a_{r-1}\om_{r-1}+a_r\om_r$ and $\om_4+a_{r-1}\om_{r-1}+a_r\om_r$. The lengths of
$W$-orbits of these two weights are at least $11395$ by \eqref{a12}, and so $\dim L(\lam) > M$.

If $a_2=1$ but $a_1 \geq 2$, then by \eqref{a11} $n=22$, $(a_2,a_1)= (1,2)$. Note that 
$$\lam \succ \lam-\al_1 = 2\om_2+a_{r-1}\om_{r-1}+a_r\om_r.$$ 
The preceding arguments show that  $\Pi(\lam)$ contains the weights
$\om_1+\om_3 +a_{r-1}\om_{r-1}+a_r\om_r$ and $\om_4+a_{r-1}\om_{r-1}+a_r\om_r$, and again $\dim L(\lam) > M$.

Hence $a_2=1$ and $a_1 \leq 1$. If $a_r=a_{r-1}=0$, then $\lam = \om_2$ or $\om_1+\om_2$.
So assume that $a_r+a_{r-1} > 0$.

Suppose that $a_{r-1} > 0$. Then by Lemma \ref{orb} 
\begin{equation}\label{a13}
\Stab_W(\lam) \leq \langle \rho_1,\rho_3,\rho_4, \ldots,\rho_{r-2},\rho_r \rangle 
   = \SBS_2 \times \SBS_{n-4} \times \SBS_2.
\end{equation}   
It follows that 
$$\dim L(\lam) \geq [W:\Stab_W(\lam)] \geq \frac{n!}{(n-4)! \cdot (2!)^2} > M,$$ 
a contradiction. Hence $a_{r-1}=0$. 

Suppose next that $a_r \geq 2$. Then 
$$\lam \succ \nu_1:= \lam-\al_r = \lam-(2\om_r-\om_{r-1}) = a_1\om_1+\om_2 +(a_{r-1}+1)\om_{r-1}+(a_r-2)\om_r.$$
Thus $\nu_1 \in \Pi(\lam) \cap \LP$, and \eqref{a13} applied to $\nu_1$ shows that 
$$\dim L(\lam) \geq |\sO(\nu_1)| \geq [\SBS_n:(\SBS_2 \times \SBS_{n-4} \times \SBS_2)] > M.$$
We may therefore assume that $a_r = 1$. If $a_1=0$ then $\lam = \om_2+\om_r$. Otherwise $\lam = \om_1+\om_2+\om_r$, in which
case $\lam \succ \nu_2:= \lam-(\al_1+\al_2) = \lam-(\om_1+\om_2-\om_3) = \om_3+\om_r$. But in such a case,
\eqref{a11b} applied to $\nu_2$ shows that
$$\dim L(\lam) \geq |\sO(\nu_2)| \geq (n-3)\binom{n}{3} > M.$$ 

\smallskip
(e) We may now assume that $a_i = 0$ for $2 \leq i \leq r-1$, i.e.\ $\lam = a\om_1+b\om_r$ with $s \geq a \geq b \geq 0$ (by symmetry).
Suppose $a \geq 4$, whence $n=22$ and $a=4$ by \eqref{a11}. Then $\lam \succ \lam-2\al_1 = 2\om_2+b\om_r$, and 
so $\dim L(\lam) > M$ by the first paragraph of (d).

Suppose $a=3$ but $b > 0$. Then $\lam \succ \nu_3:= \lam-(2\al_1+\al_2) = \lam-(3\om_1-\om_3) = \om_3+b\om_r$, and \eqref{a11b}
applied to $\nu_3$ shows that $\dim L(\lam) \geq |\sO(\nu_3)| > M$. Hence $\lam = \om_3$.

Suppose $a=b=2$, i.e.\ $\lam = 2\om_1+2\om_r$. Then $\lam \succ \nu_4:= \lam-(\al_1+\al_r) = \om_2+\om_{r-1}$. In such a case,
$$\dim L(\lam) \geq |\sO(\nu_4)| \geq \dfrac{n!}{(n-4)! \cdot (2!)^2} > M,$$ 
by using \eqref{a13} for $\nu_4$. So $a+b \leq 3$, and thus
$\lam = 2\om_1+\om_r$, $2\om_1$, $\om_1+\om_r$, $\om_1$, or $0$.
\end{proof}

To handle the other classical groups, we first consider a special case. Again, we use the weight labeling as in \cite{OV}.

\begin{prop}\label{bcd}
Let $G$ be a simply connected simple algebraic group over $\F$ of type $B_r$, $C_r$, or $D_r$, with $r \geq 7$. 
Then 
$$\dim L(\om_1+\om_2) \geq \left\{ \begin{array}{ll}4r(r^2-1)/3, & \mbox{if }p=3,\\ 
  4r(r-1)(2r-1)/3,& \mbox{if }p \neq 3.\end{array}\right.$$
\end{prop} 

\begin{proof}
Note that $|\sO(\om_3)| = 8\binom{r}{3}$ and $|\sO(\om_1+\om_2)| = 4r(r-1)$. Since $\lam:=\om_1+\om_2$ is the highest weight of
$L(\lam)$, it suffices to show that $\mu:=\om_3$ is a weight of $L(\lam)$, with multiplicity $m_{L(\lam)}(\mu)\geq 1$ if $p=3$ and 
$m_{L(\lam)}(\mu) \geq 2$ if 
$p \neq 3$. If $p=3$, then $\lam \succ \lam-(\al_1+\al_2) = \mu$, whence $\mu \in \Pi(\lam)$ by Theorem \ref{premet}, and we are done.

In what follows we may assume $p \neq 3$. We realize the roots and the weights of $G$ using an orthonormal basis 
$(e_i \mid 1\leq i \leq r)$ of $\R^r$ (with scalar product $(\cdot,\cdot)$); 
in particular, $\om_1=e_1$, $\om_2=e_1+e_2$, $\om_3=e_1+e_2+e_3$. 
Consider the simple (Weyl) module $V(\lam)$ of the corresponding algebraic group
over $\C$. Then $\nu \in \LP$ is a weight of $V(\lam)$ precisely when $\lam \succ \nu$. Writing $\om_i$ in terms of simple roots
(see \cite[Table 2]{OV}), it is straightforward to check that this is equivalent to
\begin{equation}\label{w10}
   \nu \in \left\{ \begin{array}{ll} \{\lam,\mu, \om_1\}, & \mbox{if }G = C_r \mbox{ or }D_r,\\
   \{\lam,\mu,\om_2,2\om_1,0\}, & \mbox{if }G = B_r.\end{array} \right.
\end{equation}    
Next we use Freudenthal's formula \cite[p.\ 122]{H} to find the multiplicity $m_{V(\lam)}(\mu)$ of $\mu$ as a weight of $V(\lam)$. 
Then we must 
find all multiples $l\al$ of positive roots $\al$, with $l \in \Z_{\geq 1}$, such that $\mu+l\al$ is a weight of $V(\lam)$, i.e.\ 
$W$-conjugate to one of the weights listed in \eqref{w10}. Again using \cite[Table 2]{OV}, we can check that this happens 
precisely when $l=1$ and $\al = e_1-e_2$, $e_1-e_3$, $e_2-e_3$. In all these cases, $\mu+l\al$ is $W$-equivalent to $\lam$,
and we readily obtain $m_{V(\lam)}(\mu)=2$.

Suppose now that 
\begin{equation}\label{w11}
  m_{L(\lam)}(\mu) \leq 1. 
\end{equation}
This can happen only when $\mu$ is a weight of some composition factor 
$L(\nu_0)$ of a reduction modulo $p$ of $V(\lam)$, where $\nu_0$ is listed in \eqref{w10}. Note that if $\nu$ is such a weight
and $\nu \neq \mu$, then $\nu \not{\succ} \mu$. It follows that $\nu_0=\mu$, 
i.e.\ $L(\mu)$ is a composition factor of
a reduction modulo $p$ of $V(\lam)$. By the linkage principle, see \cite{Jan}, this implies that $w \circ \lam = \mu$ for 
some $w \in W_p$.
Here, the affine Weyl group $W_{p}$ is generated by the map $\rho_{\alpha,l} \circ \lambda = \rho_{\alpha} \circ \lambda + lp \alpha$,
where $\alpha$ is a simple root, $l \in \Z$ and $\rho_{\alpha}$ denotes the 
reflection corresponding to $\alpha$; furthermore, $\rho_\alpha \circ \lambda = \rho_\alpha(\lambda+\delta) - \delta$,
where $\delta:=\sum^r_{i=1}\om_i$.

Write $|v|^2$ instead of $(v,v)$ for $v \in \R^r$. Then, for any root $\alpha$ and any $l \in \Z$ 
$$|\rho_{\alpha,l} \circ \lambda + \delta|^{2}
  = |\lambda + \delta|^{2} + l^{2}p^{2}|\alpha|^{2} +
   2lp(\rho_{\alpha} (\lambda + \delta),\alpha).$$
Observe that $|\alpha|^{2} \in \Z$ and 
$$\bigl(\rho_{\alpha} (\lambda + \delta),\alpha\bigr) = \left((\lambda + \delta) -
  \frac{2(\lambda + \delta,\alpha)}{(\alpha,\alpha)},\alpha\right) =
  - (\lambda + \rho,\alpha) \in \Z.$$
It follows that $|\rho_{\alpha,l} \circ \lambda+\delta|^2 \equiv 
|\lambda+\delta|^2 \pmod*{\gcd(2,p)}$. Thus we have shown
\begin{equation}\label{link}
 \mbox{If two weights }\lambda,~\mu \mbox{ are linked, then }|\lambda+\delta|^2 \equiv |\mu+\delta|^2 \pmod*{\gcd(2,p)p}.
\end{equation}   
(Note that a slightly weaker result than \eqref{link}, namely only modulo $p$, was obtained in \cite[Lemma 2.1]{T}.) In our case,
$|\lambda+\delta|^2 - |\mu+\delta|^2=6$. Applying \eqref{link}, we conclude that \eqref{w11} can happen only when $p=3$. 
\end{proof}

\begin{thm}\label{class-BCD}
Let $G$ be a simply connected simple algebraic group over $\F$ of type $B_r$, $C_r$, or $D_r$, with $r \geq 12$. 
Let $L(\lam)$ be an  irreducible $\F G$-representation, which is restricted if $p := \Char(\F) > 0$.
Suppose that $\dim L(\lam) \leq M$, where $M := 2(r+1)^3$ if $r \geq 14$, and $M:=r^3$ if $r=12,13$.
Then $\lam = a\om_1$ with $0 \leq a \leq 3$, $\om_2$, $\om_3$, or $\om_1+\om_2$. Moreover, 
if $\lam=\om_1+\om_2$ and $r \geq 16$, then $p=3$.
\end{thm}

\begin{proof}
Write the highest weight $\lam$ as $\sum^r_{i=1}a_i\om_i$ with $a_i \in \Z_{\geq 0}$ (note that $L(\om_1)$ is the natural module
for $G$).  

\smallskip
(a) First suppose that $a_i > 0$ for some $r-2 \leq i \leq r$. In this case, by Lemma \ref{orb}, the length of the $W$-orbit 
$\sO(\lam)$ of $\lam$ is 
at least $[W:\Stab_W(\lam)] \geq 2^{r-1} > M$.

Next suppose that $a_i > 0$ for some $4 \leq i \leq r-3$. In this case, if $G$ is of type $X_r$ so that $W=W(X_r)$, then
$\Stab_W(\lam)$ is contained in $W(A_{i-1}) \times W(X_{r-i})$ by Lemma \ref{orb}, hence
$$\dim L(\lam) \geq |\sO(\lam)| \geq 2^i \binom{r}{i} \geq \min\bigl( 2^4 \binom{r}{4},2^{r-3}\binom{r}{3} \bigr) > 2(r+1)^3,$$
since $r \geq 9$. 

Applying this argument to $\mu = \sum^r_{i=1}b_i\om_i \in \Pi(\lam) \cap \LP$, 
we deduce that
\begin{equation}\label{bcd1}
  b_i =0 \mbox{ for } 4 \leq i \leq r.
\end{equation} 

\smallskip
(b) Suppose $a_3 > 0$. If $a_1 > 0$ or $a_2 > 0$, then $\Stab_W(\lam) \leq W(A_1) \times W(X_{r-3})$, whence 
\begin{equation}\label{bcd2}
  |\sO(\lam)| \geq 4r(r-1)(r-2) > 2(r+1)^3
\end{equation}  
(since $r \geq 9$), contradicting the bound on $\dim L(\lam)$. So $\lam=a_3\om_3$. Now, if $a_3 \geq 2$ (and so $p \neq 2$ 
as $\lam$ is restricted), then by Theorem \ref{premet} 
$$\lam \succ \lam-\al_3 = \om_2+\om_4 \in \Pi(\lam) \cap \LP,$$ 
violating \eqref{bcd1}. Hence $\lam=\om_3$ in this case.

We have shown that $a_3=0$. Suppose $a_2 > 0$. Now, if $a_2 \geq 2$ (so again $p \neq 2$), then 
$$\lam \succ \lam-\al_2 = (a_1+1)\om_1 + (a_2-2)\om_2+\om_3 \in \Pi(\lam),$$
leading to a contradiction by applying \eqref{bcd2} to $\lam-\al_2$. If $a_2=1$ but $a_1 \geq 2$, then $p \neq 2$ 
and 
$$\lam \succ \lam-(\al_1+\al_2) = (a_1-1)\om_1 + \om_3 \in \Pi(\lam),$$
again yielding a contradiction by applying \eqref{bcd2} to $\lam-\al_1-\al_2$.  Hence $\lam \in \{\om_2,\om_1+\om_2\}$ in 
this case.

We are left with the case $\lam = a\om_1$. If $a \geq 4$, then again $p \neq 2$ and 
$$\lam \succ \lam-(2\al_1+\al_2) = (a_1-3)\om_1 +\om_3 \in \Pi(\lam),$$
leading to a contradiction by applying \eqref{bcd2} to $\lam-2\al_1-\al_2$. So $0 \leq a \leq 3$ as stated.

\smallskip
(c) We make some more comments about the cases $\lam \in \{3\om_1,\om_3,\om_1+\om_2\}$. Note 
that 
$$3\om_1-\al_1 = \om_1+\om_2,~(\om_1+\om_2)-(\al_1+\al_2)=\om_3.$$ 
So, assuming $p \neq 2$ when $\lam=\om_1+\om_2$, we may assume $\om_3 \in \Pi(\lam)$ in all these three cases. It now follows from Lemma \ref{orb} that in these cases
$$\dim L(\lam) \geq |\sO(\om_3)| = 4r(r-1)(r-2)/3 > r^3.$$ 
Proposition \ref{bcd} shows that if $p=2$, then $\dim L(\om_1+\om_2) > 8\binom{r}{3} > r^3$ (as $r \geq 12$), and 
if $r \geq 16$ and $p \neq 3$, then $\dim L(\om_1+\om_2) > 4r(r-1)(2r-1)/3 > 2(r+1)^3$.
Thus, 
\begin{equation}\label{bcd3}
  \mbox{If }r \geq 12\mbox{ and }\dim L(\lam) \leq r^3, \mbox{ then }\lam \in \{0,\om_1,2\om_1,\om_2\}.
\end{equation}
This statement \eqref{bcd3} was 
recorded in \cite[Theorem 5.1]{Lu}, but we note that the treatment of the weights $a_1\om_1+a_2\om_2$ therein is 
incorrect. Also note that $\dim L(\om_1+\om_2)$ may be smaller than $2(r+1)^3$ when $p=3$, see \cite{Lu} for examples.
\end{proof}

\section{A dichotomy for monodromy groups}

We will need some preliminary facts:

\begin{lem}\label{eqns}
Let $n \in \Z_{\geq 1}$ and let $D := 2^n(2^{2n+1}-1)$. Then none of the following equations
\begin{enumerate}[\rm(i)]
\item $D=x^2-1$,
\item $D= x(x-1)/2$,
\item $D= x(x-1)/2-1$ and $n \geq 2$,
\end{enumerate}
has a solution in the positive integers.
\end{lem}

\begin{proof}
(i) Suppose $x^2-1=D$ for some $x \in \Z_{\geq 1}$. Checking the cases $1 \leq n \leq 7$ directly, we may assume $n \geq 8$. 
Now $x > 1$ is odd, and $\gcd(x-1,x+1)=2$, but $2^n|(x^2-1)$. It follows that there is some $\eps= \pm 1$ such that
$2^{n-1}|(x-\eps)$. Write $x-\eps = 2^{n-1}y$ for some $y \in \Z_{\geq 1}$. Then 
$$2^{3n+1}-2^n = D = x^2-1=(2^{n-1}y+\eps)^2-1 = 2^{2n-2}y^2+2^n\eps y,$$
and so $\eps y+1 = 2^{n-2}(2^{n+3}-y^2)$. This implies $y > 1$, and $y+\eps$ is divisible by $2^{n-2}$. Hence 
$y+\eps = 2^{n-2}z$ for some $z \in \Z_{\geq 1}$.  In this case, $y \geq 2^{n-2}-1$, $x \geq 2^{2n-3}-2^{n-1}-1$,
and so $x^2-1 > 2^{4n-7} \geq 2^{3n+1} > D$ (as $n \geq 8$), a contradiction.

\smallskip
(ii) Suppose $x(x-1)/2=D$ for some $x \in \Z_{\geq 1}$. 
Then $\gcd(x-1,x)=1$, but $2^{n+1}|x(x-1)$. It follows that there is some $\eps \in \{0,1\}$ such that
$2^{n+1}|(x-\eps)$. Write $x-\eps = 2^{n+1}y$ for some $y \in \Z_{\geq 1}$. If $\eps=0$, then
$$2^{3n+2}-2^{n+1} = 2D = x(x-1)=2^{2n+2}y^2-2^{n+1}y,$$
and so $y-1 = 2^{n+1}(y^2-2^{n})$. This implies $y > 1$, and $y-1$ is divisible by $2^{n+1}$. Hence 
$y-1 \geq 2^{n+1}$, $x > 2^{2n+2}$,
and so $x(x-1) > 2^{4n+4} > 2^{3n+1}> D$, a contradiction.
If $\eps=1$, then
$$2^{3n+2}-2^{n+1} = 2D = x(x-1)=2^{2n+2}y^2+2^{n+1}y,$$
and so $y+1 = 2^{n+1}(2^{n}-y^2)$. This implies $1 \leq y < 2^{n/2}$, and $y+1$ is divisible by $2^{n+1}$, which is impossible. 

\smallskip
(iii) Suppose $n \geq 2$ and $D=x(x-1)/2-1 = (x+1)(x-2)/2$ for some $x \in \Z_{\geq 1}$. 
Then $\gcd(x+1,x-2)|3$, but $2^{n+1}|(x+1)(x-2)$. It follows that there is some $\eps \in \{-1,2\}$ such that
$2^{n+1}|(x-\eps)$. Write $x-\eps = 2^{n+1}y$ for some $y \in \Z_{\geq 1}$. If $\eps=-1$, then
$$2^{3n+2}-2^{n+1} = 2D = (x+1)(x-2)=2^{2n+2}y^2-3\cdot2^{n+1}y,$$
and so $3y-1 = 2^{n+1}(y^2-2^{n})$. This implies that $3y-1 \geq 2$ is divisible by $2^{n+1}$. Hence 
$y \geq (2^{n+1}+1)/3 > 2^{n-1}$, $x+1 > 2^{2n}$,
and so $(x+1)(x-2) > 2^{2n}(2^{2n}-3) >2^{3n+1}> D$ (as $n \geq 2$), a contradiction.
If $\eps=2$, then
$$2^{3n+2}-2^{n+1} = 2D = (x+1)(x-2)=2^{2n+2}y^2+ 3 \cdot 2^{n+1}y,$$
and so $3y+1 = 2^{n+1}(2^{n}-y^2)$. This implies $1 \leq y < 2^{n/2} \leq 2^{n-1}$, and $3y+1$ is divisible by $2^{n+1}$, which is impossible when $n \geq 2$.
\end{proof}

Note that if $(n,D) = (1,14)$, then $D = \binom{6}{2}-1 = \binom{6}{3}-6=(3^3+1)/2$.

\begin{thm}\label{inf-m}
Suppose the sheaf $\sF(q,f)$ in \eqref{airy-g}, of rank $D=2^n(2^{2n+1}-1)$, has infinite geometric 
monodromy group $G=G_\geo$. Then $G^\circ=\SL_D$. Under the more restrictive condition \eqref{airy-s} 
$G=\SL_D$; moreover, the sheaf $\sG(q,f,t(q))$ has geometric monodromy group equal to $G$.
\end{thm}

\begin{proof}
By Theorem \ref{cond-S}, $G$ satisfies $\CSP$. 
Thus $G^\circ \leq \SL_D$, and 
$\ZB(G)$ is finite, but $G^\circ$ is infinite. It follows from \cite[Lemma 1.4]{KT5} that $G^\circ$ is irreducible on 
the underlying representation $V$, i.e.\ $\sF(q,f)$ is 
Lie-irreducible. By Theorem \ref{comptsize}, $G^\circ$ is a simple algebraic group of dimension $\geq D$.

\smallskip
(a) Suppose $n \geq 3$, so that $D \geq 1016$. Then $G^\circ$ must be a classical group of rank say $r$, where 
$r(2r+1) \geq \dim G^\circ \geq D \geq 1016$, whence $r \geq 23$. Also, if $G^\circ$ is of type $A_r$, then
\begin{equation}\label{inf-a}
  D \leq \dim G^\circ = r(r+2) < \binom{r+1}{4},
\end{equation}  
and if $G^\circ$ is of type $B_r$, $C_r$, or $D_r$, then 
\begin{equation}\label{inf-b}
  D \leq r(2r+1) < r^3.
\end{equation}  
In the case of $A_r$, we can apply Theorem \ref{class-A} and, using the dimension formula for $L(\lam)$ given in 
\cite[Table 5]{OV} and the bound \eqref{inf-a}, we see that the highest weight $\lam$ of the $G^\circ$-module $V$ is, up to duality, 
$a\om_1$ with $a = 1,2$, $\om_2$, or $\om_1+\om_r$. If $\lam=\om_1$, then $D=r+1$, $G^\circ = \SL_D$,
and hence $G = \SL_D$ as stated. In the other cases, $D = \binom{r+1}{2}$, 
$\binom{r}{2}$, or $(r+1)^2-1$; all are impossible by Lemma \ref{eqns}.

In the case of types $B_r$, $C_r$, and $D_r$, we can apply Theorem \ref{class-BCD} (more precisely, \eqref{bcd3}), and, using the dimension formula for $L(\lam)$ given in \cite[Table 5]{OV} and the bound \eqref{inf-b}, we see that the highest weight $\lam$ of the $G^\circ$-module $V$ is $a\om_1$ with $a = 1,2$, or $\om_2$. If $\lam=\om_1$, then $G^\circ = \Sp(V)$ or $\SO(V)$, whence
$\sF(q,f)$ is Lie-self-dual, contrary to Corollary \ref{cor:sd}.
In the other cases, $D = \binom{m}{2}$ or $\binom{m}{2}-1$ for some integer $m \geq 2$,  
and this is impossible by Lemma \ref{eqns}.

\smallskip
(b) Suppose $n=2$, so that $D=124$. Then $G^\circ$ is of type $A_r$ with $r \geq 11$, $B_r$, $C_r$, or $D_r$ with $r \geq 8$, 
$E_7$, or $E_8$. Neither $E_7$ nor $E_8$ has irreducible representations of degree $124$, see \cite{Lu}, so $G^\circ$ is classical of 
rank $r$. If $G^\circ$ is of type $B_r$, $C_r$, or $D_r$ with $8 \leq r \leq 11$, then using \cite{Lu} we check that $G^\circ$ has no 
irreducible representation of degree $124$. So $r \geq 12$, and we apply Theorem \ref{class-A}, respectively \eqref{bcd3}, as 
above to conclude that $G = \SL_D$.

Finally, assume that $n=1$, so that $D=14$. The arguments in part (b) of the proof of Theorem \ref{main1} repeated verbatim show
that either $G^\circ = \SL_D$, or $V|_{G^\circ}$ is self-dual. In the former 
case $G = \SL_D$ as in (a), and the latter case
is ruled out by Corollary \ref{cor:sd}.

\smallskip
(c) Now assume \eqref{airy-s}. Then we can consider the descent $\sG(q,f,t(q))$, for which the trace function takes values in $\Q(\zeta_4)$
and all slopes are less than $1$. Since $p=2$, Theorem \ref{dets}(iii)
implies that $\sG(q,f,t(q))$ has trivial determinant, and thus $\sG(q,f,t(q))$ has geometric monodromy group $H \leq \SL_D$.
But $H \geq G$, so we conclude $H=G=\SL_D$.  
\end{proof}

\begin{thm}\label{finite-m}
Suppose the sheaf $\sF(q,f)$ in \eqref{airy-s}, of rank $D=2^n(2^{2n+1}-1)$, has finite geometric monodromy group 
$G=G_\geo$. Then $G=\tw2 B_2(q)$ with $q=2^{2n+1}$.
\end{thm}

\begin{proof}
(a) Let $V$ denote the underlying representation. By Theorem \ref{cond-S}, $G$ satisfies $\CSP$ on $V$. But the dimension
$D=\dim(V)=q_0(q-1)$, with $q_0:=2^n$, is not a prime power. Hence $G$ is almost quasisimple by \cite[Lemma 1.1]{KT5}:
$S \lhd G/\ZB(G) \leq \Aut(S)$ for some finite, non-abelian simple group $S$. 
Then the quasisimple subgroup $L:=G^{(\infty)}$ acts irreducibly on $V$ by \cite[Lemma 1.4]{KT5}. 

The condition \eqref{airy-s} allows us to consider the descent $\sG(q,f,t(q))$ on $\G_m$, with geometric monodromy group
$H$. Then $G \lhd H$ of finite index, whence $H$ is finite and satisfies $\CSP$, and $L=H^{(\infty)}$, 
as $H/G \inj C_{t(q)}$. The representation of $H$ on $V$ has $\infty$-slopes $\sigma:= (q_0+1)/D <1$, and
is not tame at $\infty$. Hence Theorem \ref{bound} applies to $G$ and $H$.
We collect some further facts about $(G,V)$ and the character $\varphi$ of $H$ on $V$ that we will use in the proof:

\begin{enumerate}[\rm(i)]
\item $\Q(\varphi|_G) = \Q(\varphi)=\Q(i)$. 
Indeed, by Theorem \ref{dets}(i), the arithmetic monodromy group $H_{\ari,k_0}$ of $\sG(q,f,t(q))$ 
over $k_0$ has finite determinant. But it normalizes the finite irreducible subgroup $H$, so finite determinant implies that $H_{\ari,k_0}$ is 
finite. Now, by Chebotarev density, the finiteness of $H_{\ari,k_0}$ implies that all elements of it are Frobenii, and   
all Frobenii have traces in $\Q(i)$. But $G \leq H \leq H_{\ari,k_0}$, so
$\Q(\varphi|_G) \subseteq \Q(\varphi) \subseteq \Q(i)$. 
But $V|_G$ is not self-dual by Theorem \ref{notselfdual}, hence $\Q(\varphi|_G)=\Q(i) = \Q(\varphi)$. Since 
each element of $\ZB(H)$ acts as a root of unity on $V$, and the only roots of unity in $\Q(i)$ are in $\mu_4$, 
both $|\ZB(G)|$ and $|\ZB(H)|$ divide $4$. 
\item If $n \neq 2$, then $G$ is perfect and hence $G=L$. 
Indeed, by \cite[Proposition 6]{Abh}, $\pi_1(\A^1/\overline{\F_p})$ has no nontrivial finite $p'$-quotient. Since $\sF(q,f)$ lives on $\A^1/\overline{\F_2}$, $G$ has no nontrivial quotient of odd order. 
On the other hand, the proof
of Theorem \ref{prim} shows that $G$ has no quotient of order $2$ when $n \neq 2$.
\item The image $J = QC$ of $I(\infty)$ has $Q = \OB_2(J)$ and $C = \langle g_\infty \rangle$, where the central order 
$\obar(g_\infty)$ is divisible by $q-1$. Indeed, since the $\infty$-slope of $\sF(q,f)$ is $1+1/D$, this implies, by
\cite[Proposition 1.14]{Ka-GKM}, that 
$I(\infty)$ acts irreducibly on $V$, of dimension $D=q_0(q-1)$. Since the image $J$ of $I(\infty)$ is cyclic of $p'$-order modulo 
the image $Q$ of $P(\infty)$, it follows that
$g_\infty$ permutes the pairwise non-isomorphic $q-1$ simple $Q$-summands on $V$, each of dimension $q_0$, transitively. 
\end{enumerate}

\smallskip
(b) Consider the case $n \geq 3$, so that $D \geq 1016$, and $G=L$ by (ii). 

\smallskip
(b1) First suppose that $S = \ABS_m$ for some $m \geq 3$. As $G/\ZB(G) = \ABS_m \inj \GL_{m-1}(\C)$, by Theorem \ref{bound}, 
$m-1 \geq 1/\sigma$. Note that $\lceil 1/\sigma \rceil = q-2q_0+1$, so $m \geq q-2q_0+2 \geq 114$. It follows
that $D = q_0(q-1) < (m^2-5m+2)/2$.  In this case, by \cite[Lemma 6.1]{GT} 
$m=D+1$, $G=\ABS_{D+1}$, and $V=\C^D$ is
the heart of the natural permutation module. But then $V$ is self-dual, contrary to (i).

\smallskip
(b2) Next suppose that $S$ is a sporadic simple group. By (iii), the maximum order $\meo(S)$ of elements in $S$ is at least 
$q-1 \geq 127$. On the other hand, $\meo(S) \leq 119$, as one can check using 
the \cite{ATLAS} (see also 
\cite[Table 2]{KT5}), a contradiction.

\smallskip
(b3) Consider the case $S$ is a simple group of Lie type in characteristic $r \neq 2$. By Theorem \ref{bound}, the degree $e$ of every nontrivial projective representation over $\overline{\F_r}$ of $S$ satisfies $e^2-1 \geq 1/\sigma$, in particular, $e \geq 11$. Similarly, the degree $d$ of every faithful linear representation over $\overline{\F_r}$ of $S$ satisfies $d \geq 1/\sigma$, in particular, $d \geq 114$.
This rules out all classical groups of types $A_m$ or $\tw2 A_m$ with $m \leq 9$, $B_m$ with $m \leq 56$, and
$C_m$ or $D_m$ with $m \leq 7$ (as $\PSp_{2m}(r^a)$ and $\mathrm{P}\Omega^\pm_{2m}(r^a)$ have faithful representations of degree $m(2m-1)$ over $\overline{\F_r}$). This also rules out exceptional groups of types $G_2$, $\tw2 G_2$, and $\tw3 D_4$. For
the remaining exceptional groups of type $F_4(s)$, $E_6(s)$, $\tw2 E_6(s)$, $E_7(s)$, and $E_8(s)$, with $s=r^a$, by 
\cite[Table A7]{KS} we have the following upper bounds 
$s(s+1)(s^2+1)$, $s(s^6-1)/((s-1)\gcd(3,s-1))$, $(s+1)(s^2+1)(s^3-1)/\gcd(3,s+1)$, $(s+1)(s^2+1)(s^4+1)/2$, and 
$(s+1)(s^2+s+1)(s^5-1)$ for $\meo(S)$, respectively. On the other hand, $q_0(q-1) = D$ is at least the smallest 
degree $\dl(S)$ of nontrivial projective complex representations of $S$, which in turn is at least 
$s^6(s^2-1)$, $s^9(s^2-1)$, $s^9(s^2-1)$, $s^{15}(s^2-1)$, $s^{27}(s^2-1)$, respectively, see e.g.\ \cite[Table I]{TZ}, and 
we arrive at a contradiction in all five cases, as $q-1 \leq \meo(S)$ by (iii) and $D < \sqrt{q^3/2}$.

If $S = \PSL_m(s)$ or $\PSU_m(s)$ with $m \geq 11$, then 
$m^2-1 \geq 1/\sigma$, and so $m^2 \geq q-2q_0+2$, by Theorem \ref{bound}. Since $m \geq 11$, by \cite[Theorem 1.1]{TZ} 
$$D \geq \frac{s^m-s}{s+1} \geq \frac{3^m-3}{4} > m^3 \geq (q-2q_0+2)^{3/2} > q_0(q-1),$$
a contradiction. For $S$ of type $BCD_m$ with $m \geq 7$, we have 
$\min(m(2m-1),2m+1) \geq 1/\sigma$, and so $2m^2 \geq q-2q_0+1$, by Theorem \ref{bound}. Since $m \geq 7$, by 
\cite[Theorem 1.1]{TZ}  
$$D \geq \frac{s^m-1}{2} \geq \frac{3^m-1}{2} > 3m^3 \geq (3/2^{3/2})(q-2q_0+1)^{3/2} > q_0(q-1),$$
again a contradiction.

\smallskip
(b4) We may now assume that $S$ is a simple group of Lie type in characteristic $2$, defined over a field $\F_s$ with $s=2^a$. 
Since $n \geq 3$, by Theorem \ref{prime-suzuki1}, $t(q)=q-2q_0+1$ admits a divisor $\ell = \ppd(2,4(2n+1))$. 
Next, we use the fact that the image of $I(0)$ in $H$ has order $q-2q_0+1$, which implies that $H$ has an element 
of prime order $\ell$ that normalizes $G=L$. But $\CB_H(G) = \ZB(H)$ has order dividing $4$ by (i), so 
\begin{equation}\label{div10}
  \ell \mbox{ divides } |\Aut(L)|.
\end{equation}   
 
First suppose that $S = \Sp_{2m}(s)$ with $m \geq 2$ or $\mathrm{P}\Omega^\pm_{2m}(s)$ with $m \geq 3$. Then, 
\eqref{div10} implies that $\ell \geq 4(2n+1)+1 \geq 29$ divides $a|S|$. If moreover $\ell \nmid a$, then $\ell$ divides
$\prod^m_{i=1}(s^{2i}-1)$, which implies $2ma \geq 4(2n+1)$ by primitivity of $\ell$. In either case, 
\begin{equation}\label{div11}
  s^m=2^{ma} \geq 2^{2(2n+1)} = q^2 \geq 2^{14}
\end{equation}  
Now, applying \cite[Theorem 1.1]{TZ} for $m \geq 3$ we obtain
$$q^{3/2} > D \geq \dl(S) \geq \frac{(s^m-1)(s^{m-1}-1)}{s^2-1} > s^{2m-3}/2 \geq s^m/2 \geq q^2/2,$$
a contradiction. If $S=\PSp_4(s)$, then $s \geq q \geq 2^7$ by \eqref{div11}, and so
$$q^{3/2} > D \geq \dl(S) = s(s-1)^2/2 \geq q(q-1)^2/2 > q^2,$$
again a contradiction.

Next suppose that $S = \PSL_m(s)$ or $\PSU_m(s)$ with $m \geq 2$. The same arguments as above show that 
\eqref{div11} still holds; in fact, $s^m \geq q^4$ for $S=\PSL_m(s)$. Assume that $S = \PSL_m(s)$ with $m \geq 2$; in particular, $s^{m-1} \geq s^{m/2} \geq q^2$. Then using \cite[Theorem 1.1]{TZ} we obtain that 
$$q^{3/2} > D \geq \dl(S) \geq \frac{s^m-s}{s-1} > s^{m-1} \geq q^{2},$$
a contradiction. Next, assume that $S = \PSU_m(s)$ with $m \geq 5$; in particular, $q^{3/2} \leq s^{3m/4} \leq s^{m-5/4}$. Then using \cite[Theorem 1.1]{TZ} we obtain that 
$$\frac{q^{3/2}}{\sqrt{2}} > D \geq \dl(S) \geq \frac{s^m-s}{s+1} \geq \frac{2}{3}(s^{m-1}-1) > \frac{2^{5/4}}{3}
   (s^{m-5/4}-1) > \frac{2^{5/4}}{3}(q^{3/2}-1) > \frac{q^{3/2}}{\sqrt{2}},$$
again a contradiction. When $S = \PSU_4(s)$, $\ell|a$ or $\ell$ divides 
$(s+1)(s^2-1)(s^3+1)(s^4-1)$, so instead of \eqref{div11} now $s^3 \geq q^2$. Again using \cite[Theorem 1.1]{TZ} we obtain
$$q^{3/2} > D \geq \dl(S) \geq \frac{s^4-1}{s+1} > s^3/2 \geq q^2/2,$$
a contradiction. Finally, if $S = \PSU_3(s)$, then $s^3 \geq q^2 \geq 2^{14}$ from \eqref{div11}.  But every irreducible 
character of $\SU_3(s)$ of even degree 
has degree divisible by $s \geq q^{2/3}$, and hence cannot be equal to
$D=q_0(q-1)$. 

Let $S$ be one of the exceptional groups of type $G_2(s)$ with $s>2$, $\tw3 D_4(s)$, $\tw2 F_4(s)'$ with $s>2$, $F_4(s)$ with 
$s>2$, $E_6(s)$, $\tw2 E_6(s)$, $E_7(s)$, and $E_8(s)$. That $\ell$ 
divides $|\Aut(S)|$ implies $q^4 \leq s^c$, where
$$c = 6,~12,~12,~12,~12,~18,~18,~30,$$
respectively. It follows that $D < q^{3/2} \leq s^{3c/8}$, with 
$$3c/8 = 9/4,~9/2,~9/2,~9/2,~9/2,~27/4,~27/4,~45/4,$$
respectively. But this contradicts the lower bounds 
$D \geq \dl(S) \geq s(s^2-1)$, $s^3(s^2-1)$, $s^4(s-1)\sqrt{s/2}$, $s^6(s^2-1)$, $s^9(s^2-1)$, $s^9(s^2-1)$, $s^{15}(s^2-1)$, $s^{27}(s^2-1)$, respectively, 
see e.g.\ \cite[Table I]{TZ}, and 
we arrive at a contradiction in all eight cases. The cases $\tw2 F_4(2)'$ and $F_4(2)$ are ruled out because
$|\Aut(S)|$ is not divisible by the prime $\ell \geq 29$. 

The only remaining case is that $S = \tw2 B_2(s)$. That $\ell$ divides
$|\Aut(S)|$ implies $q^4 \leq s^4$, i.e.\ $q\leq s$.  But 
$D=q_0(q-1)$ is the degree of some irreducible character of $G$, so $q=s$, i.e.\ $S = \tw2 B_2(q)$. In this
case also $G=L=\tw2 B_2(q)$, as stated.

\smallskip
(c) Now we consider the case $n=2$, i.e.\ $D = 124$. 

\smallskip
(c1) As the quasisimple subgroup $L=G^{(\infty)}$ acts irreducibly on
$V=\C^D$, by \cite{HM} we have the following possibilities for $S$: 
$$\PSL_2(125),~\SL_3(5),~\SL_5(2),~G_2(5),~\ABS_{125}, \mbox{ or }\tw2 B_2(32).$$ 
Here, a generator $g_0$ of the image of $I(0)$ in $H$ has order $t(q)=25$
and normalizes $G$ and $L$. Since $\CB_H(L) = \ZB(H) \leq C_4$ by (i), $25$ divides $|\Aut(S)|$. This rules out
the first three cases. 

In the two cases $S = G_2(5)$ and $\ABS_{125}$, $L=S$ and $V|_L$ is self-dual. Since $\CB_H(L) \leq C_4$ and 
$|\Out(S)| \leq 2$, it follows that $H/L$ is a $2$-group. On the other hand, $H$ has no quotient of order $2$, as shown 
in part (i) of the proof of Theorem \ref{prim} (which also works for $n=2$). 
Hence
$H=L=G$, and so $\sF(q,f)$ is self-dual. But this contradicts Theorem \ref{notselfdual}.

The only remaining case is $S=\tw2 B_2(32)$, in which case $L=S$. 

\smallskip
(c2) The rest of this paragraph applies to all $n \geq 1$, for which we know $L=S=\tw2 B_2(q)$. Since 
$\Out(S) \cong C_{2n+1}$, but $G = \OB^{2'}(G)$, we see that $G=\ZB(G)S$, with $\ZB(G) \leq C_4$. Recall that 
$H/G \inj C_{t(q)}$ and $\ZB(H) \leq C_4$. In particular, $\ZB(H)G/G \cong \ZB(H)/\ZB(G)$ has order dividing both $4$ and $t(q)$, whence
$\ZB(H)=\ZB(G)$. Next, $H/G$ acts trivially on $G/S \cong \ZB(G)$, i.e.\ $G/S \leq \ZB(H/S)$, and the quotient
$(H/S)/(G/S) \cong H/G$ is cyclic. It follows that $H/S$ is an abelian group. As mentioned above, $H$ has no quotient
of order $2$, so $2 \nmid |H/S|$, whence $\ZB(G)=1$; in particular $G=S$ when $n=2$. We have also shown
that $\CB_H(S)=\ZB(H)=1$, so 
\begin{equation}\label{div12}
  S \lhd H \leq \Aut(S) = S \rtimes C_{2n+1}
\end{equation}

\smallskip
(d) Finally, we consider the case $n=1$, i.e.\ $D = 14$. As the quasisimple group $G=L$ acts irreducibly on
$V=\C^D$ and $\Q(\varphi) = \Q(i)$, by \cite{HM} 
the only possibility is that
$G=\tw2 B_2(8)$.
\end{proof}

Now we are ready to prove the second main result of the paper:

\begin{thm}\label{main2}
For the geometric monodromy group $G=G_\geo$ of  the sheaf $\sF(q,f)$ in \eqref{airy-s}, of rank $D=2^n(2^{2n+1}-1)$, 
either $G =\SL_D$ or $G=\tw2 B_2(q)$ with $q=2^{2n+1}$. Furthermore, for any $r|t(q)$, the geometric monodromy 
group of the descent $\sG(q,f,r)$ is also equal to $G$.
\end{thm}

\begin{proof}
Suppose $G$ is infinite. Then the statements follow from Theorem \ref{inf-m}.

From now on, assume that $G$ is finite. Then $G = S= \tw2 B_2(q)$ by Theorem \ref{finite-m}. Since for any $r|t(q)$, the geometric monodromy group
of $\sG(q,f,r)$ contains $G$ and is contained in the geometric monodromy
group $H$ of $\sF(q,f,t(q))$, it suffices to show that $H=S$. 

First, note that the $\infty$-slope of $\sG(q,f,t(q))$ is $\sigma=(q_0+1)/D$ and 
$D+1=(q_0+1)t(q)$, so $\gcd(D,q_0+1)=1$. It follows from \cite[Proposition 1.14]{Ka-GKM} that $I(\infty)$ acts irreducibly on $V$, of dimension $D=q_0(q-1)$. Since the image $J$ of $I(\infty)$ is cyclic of $p'$-order modulo 
the image $Q$ of $P(\infty)$, $Q = \OB_2(J)$ and $J = \langle Q,g_\infty \rangle$, where 
the $p'$-element $g_\infty$ transitively permutes the pairwise non-isomorphic $q-1$ simple $Q$-summands on $V$, 
each of dimension $q_0$. 

Since $q=2^{2n+1}$, using \cite{Zs} we can find a primitive prime divisor $\ell = \ppd(2,2n+1)$, and fix a power
$h$ of $g_{\infty}$ that has order $\ell$. Clearly, the prime $\ell$ is at least $2n+3$, so it is coprime to 
$2n+1 = |\Out(S)|$. On the other hand, $S \leq H \leq \Aut(S)$ by \eqref{div12}. Hence $h \in S$ and $Q < S$.

We can write $h=x^{(q-1)/\ell}$, where on the natural module $U=\F_q^4$ for $S=\tw2 B_2(q) < \Sp(U)$, the spectrum of $x \in S$ 
consists of $4$ eigenvalues $\xi^{2^n}$, $\xi^{-2^n}$, $\xi^{1-2^n},~\xi^{2^n-1}$
where $\xi \in \overline{\F_2}^\times$ has order $q-1$ see \cite{Bur}, \cite{Suz}. We may write $\Aut(S) = \langle S,\theta \rangle$, where 
$\theta$ acts as the Galois automorphism $\lambda \mapsto \lambda^2$ of $\overline{\F_2}$. Suppose that for some $1 \leq a \leq 2n$ and for some $y \in S$, the element $y\theta^a$ centralizes $h$. Note that $\theta^a$ sends $x$ to $x^{2^a}$ and 
$h=x^{(q-1)/\ell}$ to $x^{2^a(q-1)/\ell}$. It follows that 
$$x^{2^a(q-1)/\ell} = \theta^a h \theta^{-a} = y^{-1}(y\theta^a)h(y\theta^a)^{-1}y = y^{-1}hy$$ 
is $S$-conjugate to $h=x^{(q-1)/\ell}$. On the other hand, it is known \cite{Bur} that if $b,c \in \Z$ then $x^b$ and $x^c$ are 
$S$-conjugate if and only if $c \equiv \pm b \pmod*{(q-1)}$. It follows that $\ell$ divides $2^a \pm 1$; in particular,
$\ell$ divides $2^{2a}-1$. The primitivity of $\ell$ then implies that $2n+1$ divides $2a$, a contradiction. 

We have shown that $\CB_{\Aut(S)}(h) \leq S$. As $g_\infty$ centralizes $h$, it follows that $g_\infty \in S$. Hence
$J = Q\langle g_\infty\rangle < S$. Now $S \lhd H$ and $H$ is finite, so $H=S$ by Theorem \ref{Iinftygen}.
\end{proof}

Theorems \ref{main1} and \ref{main2} imply the following.

\begin{cor}\label{main2c}
The sheaves $\sF_q$ and $\sG_q$ in \S\ref{suz-desc}, of rank $D=2^n(2^{2n+1}-1)$, have the same geometric monodromy group
$G=G_\geo$. Furthermore, either $G=\tw2 B_2(q)$ with $q=2^{2n+1}$, or $n \geq 2$ and $G = \SL_D$.
\end{cor}

\begin{rmk}\label{inf-ex}
It is plausible that for each $q = 2^{2n+1} \geq 32$, we can find a polynomial $f_1 \in \F_2[x]$ of degree $2^n+1$ such that the sheaf
$\sF(q,f)$ with $f(x)=f_1(x^{t(q)})$ as in \eqref{airy-s} has infinite geometric monodromy group, which then is $\SL_D$ by
Theorem \ref{main2}. Indeed, a {\sc Magma} calculation shows for each $2 \leq n \leq 25$ that 
$\Frob_{1,\F_{2^{k(n)}}}$ has non-integral trace on $\sF(q,x^{(1+2^n)t(q)})$ for some choice of $k(n) \in \Z_{\geq 1}$. 
For instance, $\Frob_{1,\F_{2^{k(n)}}}$ has 
trace $(2i-7)/2$ for $(n,(k(n))=(2,7)$, $(5i+3)/2$ for $(n,k(n))=(3,5)$, $(7-3i)/4$ for $(n,k(n))=(4,7)$, and 
$5/2$ for $(n,k(n))=(5,7)$. In fact, for infinitely many integers $n \geq 2$, we offer in Theorem \ref{inf-g} a construction of 
a sheaf $\sF(q,f)$ with $G_\geo=\SL_D$.  
\end{rmk}

\section{Arithmetic vs. geometric monodromy groups}
\subsection*{9A. Glauberman and Dade correspondences}
Our next results depend on the Glauberman correspondence. Recall that if $A$ is a solvable finite group acting
by automorphisms on another finite group $S$ with $(|A|,|S|)=1$, then there exists a canonical
bijection $^*:{\rm Irr}_A(S) \rightarrow \irr C$, where $C=\cent SA$ is the  fixed-points subgroup and ${\rm Irr}_A(S)$
is the set of  $A$-invariant irreducible characters of $S$.  (See \cite[Chapter 13]{Is}.)
Since the map is canonical, it is not difficult to see that $\Q(\theta)=\Q(\theta^*)$, where $\Q(\theta)$ is the field of values of $\theta$.
(See \cite[Problem 3.1]{Is}.)

\medskip
Suppose now that $S=\tw2 B_2(2^n)$, where $n$ is odd. Then it is well-known that
$S$ admits a field automorphism $a$ of order $n$. Assume further that $(m, |S|)=1$ for some divisor $m > 1$ of $n$.
Consider $A=\langle a^{n/m} \rangle \cong C_m$. If $m = n$, 
 then $C=\cent SA=\tw2 B_2(2) \cong C_5 \rtimes C_4$ has 5 irreducible
characters; two of its four linear characters are rational, the others have 
field of values $\Q(i)$. 
It also has a rational irreducible character of degree 4. 
In particular, it follows that $\irr S$ has exactly 5 irreducible $A$-invariant characters, and exactly two 
of them have field of values $\Q(i)$. On the other hand, if $m < n$, 
 then $C=\cent SA=\tw2 B_2(2^{n/m})$ has exactly 2 irreducible $A$-invariant characters with field of values $\Q(i)$, namely
 the ones of degree $(r-1)\sqrt{r/2}$ with $r:=2^{n/m}$;
 see \cite{Bur}.
This proves the following.

\begin{lem}\label{coprime1}
Suppose that $S=\tw2 B_2(2^n)$, where $n$ is odd.
Assume further that $(m, |S|)=1$ for a divisor $m>1$ of $n$, and let $A$ be the subgroup of field automorphisms of $S$ of order $m$.
Let $C=\cent SA$. If $\theta \in {\rm Irr}_A(S)$ has field of values $\Q(i)$, then
$\theta^*$ has degree $(r-1)\sqrt{r/2}$ with $r:=2^{n/m}$.
\end{lem}

\begin{lem}\label{coprime2}
Suppose that $A$ is a cyclic group of order $m$  acting faithfully and coprimely on $S$.
Let $\theta \in \irr S$ be $A$-invariant, $C=\CB_S(A)$, and let $\eta \in \Irr(C)$ be the $A$-Glauberman correspondent of $\theta$. 
Let $G=S \rtimes A$ be the semidirect product,  and let $\psi \in \irr G$ be an extension of $\theta$. 
Consider $x \in G \smallsetminus S$ such that $\langle x,S \rangle = G$.
\begin{enumerate}[\rm(i)]
\item If $\eta$ is linear, then $|\psi(x)|=1$.
\item In all cases, there exist a root of unity $\gamma \in \C^\times$ of order dividing $2m$ and 
$c \in C$ such that $\psi(x)=\gamma \eta(c)$.
\end{enumerate}
\end{lem}

\begin{proof}
Let $\psi \in \irr G$ be the canonical extension of $\theta$ to $G$.
(This is the unique extension $\psi$ such that the determinantal order is coprime with $|A|$, see \cite[Corollary 6.28]{Is}).
Since every extension of $\theta$ to $G$ is a multiple of $\psi$ by a linear character $\lambda$ of $G/S$, 
we may assume that $\psi = \lambda\psi_0$, where $\psi_0$ is the canonical extension.

Observe that $m$ divides the order of $x$. Since $m$ is coprime to $|S|$, we can write $x = cb=bc$, with
$b$ being a $\pi$-element and $c$ being a $\pi'$-element, where $\pi$ is the set of prime divisors of $m$.  
Also $|b|=m$, and $S \ni x^m = b^mc^m = c^m$, implying $c \in S$ since $\gcd(m,|S|)=1$. 
Moreover, $G = S \rtimes \langle b \rangle$, so without any loss we may replace $A$ by $\langle b \rangle$, and now 
$c \in C=\cent Sb=\cent SA$. 
By \cite[Theorem 13.6]{Is}
$$\psi_0(x)=\psi_0(cb)=\epsilon \eta(c) \, ,$$
where $\epsilon = \pm 1$.  
Hence, $\psi(x) = \gamma \eta(c)$, where $\gamma=\epsilon\lambda(x)$ and $\lambda(x)^m=1$.
\end{proof}

Since each of the two irreducible characters of $S:=\tw2 B_2(q)$ of degree $(q-1)\sqrt{q/2}$ has field of values $\Q(i)$ and is 
$\Aut(S)$-invariant, Lemmas \ref{coprime1} and \ref{coprime2} imply the following.

\begin{cor}\label{ext-suz}
Let $q=2^n$ with $2 \nmid n$ and let $\theta$ be either of the two irreducible characters of $S:=\tw2 B_2(q)$ of degree $(q-1)\sqrt{q/2}$. 
Then $\theta$ extends to $G = \Aut(S) \cong S \rtimes C_n$. If furthermore $n$ is coprime to $|S|$, then 
$|\psi(x)|=1$ for any extension $\psi$ of $\theta$ to $G$ and for any $x \in G \smallsetminus S$ with $G=\langle x,S \rangle$. 
\end{cor}

We remark that the character $\theta$ in Corollary \ref{ext-suz} has a canonical extension
to $G$. Indeed, by \cite[Theorem A]{N1}, there exists a unique $\psi \in {\rm Irr}(G)$ such that the field
of values of $\psi$ is $\Q(i)$. In particular, notice that the restriction $\psi$ to $C_n$ is rational-valued.

In some non-coprime situations, we can use the following statement. This also follows from results in \cite[\S9]{Da}, but in the situation under consideration, our approach is more straightforward.  

\begin{lem}\label{coprime3}
Suppose that $G=SA$, where $A=\langle a\rangle$,  $S\nor G$ and $A\cap S=1$.
Let $C=\cent SA$, and assume that $\chi \in {\rm Irr}_A(S)$ has an extension 
$\psi \in \irr G$ such that $\psi_A$ is rational-valued. 
Suppose that for every $x \in G \smallsetminus S$, there exists some  $g \in G$ such that $x^g=cb$ for some $c \in C$ and $b \in A$. Then there exist a character $\theta \in \irr C$ and a sign $\epsilon = \pm$
such that $\psi(cb)=\epsilon \theta(c)$ for every $c \in C$ and every generator 
$b$ of $A$. 
\end{lem}

\begin{proof}
By the proof of \cite[Lemma 13.5]{Is}, we can write
$$\psi_{CA}=\sum_{\beta \in \irr C, [\chi_C, \beta] \ne 0} \beta \boxtimes \psi_\beta \, ,$$
where $\psi_\beta$ is a rational-valued character of $S$.
Now, define $\theta(c)=\psi(ca)$ for $c \in C$.
Then
$$\theta=\sum_\beta \psi_\beta(s) \beta \, .$$
Thus $\theta$ is a virtual character of $C$. We claim that $[\theta,\theta]=1$. 
Let $T$ be a set of representatives for the right cosets of $C$ in $S$.
In order to use the proof of \cite[Theorem 13.6]{Is}, we claim that
$$Sa=\bigcup_{t \in T} (Ca)^t$$
is a disjoint union. By hypothesis, 
$$Sa=\bigcup_{s \in S} (Ca)^s=\bigcup_{t \in T} (Ca)^t \, .$$
Now,
$$|S|=|Sa|=|\bigcup_{t \in T} (Ca)^t|\le \sum_{t \in T} |(Ca)^t|=
|T||C|=|S| \, ,$$ and the claim follows. The rest of the proof follows as 
in \cite[Theorem 13.6]{Is}.
\end{proof}

 
Finally, to address the general situation in the Suzuki case, we must go much deeper into \cite[\S9]{Da}.

\begin{cor}\label{ext-suz2}
Let $q=2^n$ with $2 \nmid n$ and let $\theta$ be either of the two irreducible characters of $S:=\tw2 B_2(q)$ of degree $(q-1)\sqrt{q/2}$. 
Then $\theta$ extends to $G = \Aut(S) \cong S \rtimes C_n$ and 
$|\psi(x)|=1$ for any extension $\psi$ of $\theta$ to $G$ and for any $x \in G \smallsetminus S$ with $G=\langle x,S \rangle$. 
\end{cor}

\begin{proof}
We adapt the notation of \cite[\S9]{Da} to ours.
  We know that $\cent SA=B  \rtimes R$, where $B=\langle \beta \rangle$ has order 5,
 and $R$ is a cyclic group of order 4.  Now, let   $G'=\norm GB=S'\rtimes A$,
 where $S'=\norm SB$. Let $G_0$ be the set of $x \in G$ such that $\langle xS\rangle=G/S$,
 and let $G'_0$ be the set of $x \in G'$ such that $\langle xS'\rangle=G'/S'$. 
By \cite[Lemma 9.3]{Da},
 $S'=C \rtimes R$, where $C$ is a cyclic group described there.
 By \cite[Proposition 9.7]{Da}, $G'_0$ is a trivial intersection subset of $G$
 with normalizer $G'$. In Dade's language, it satisfies (6.4) of \cite{Da}. In particular, by 
 \cite[Lemma 6.5]{Da}, 
 $$G_0=\bigcup_{\tau \in T} (G'_0)^\tau$$ is a disjoint union,
 where $G=\bigcup_{\tau \in T} G' \tau$ is a disjoint union. (In 
that lemma
right cosets are used.)
 In particular, if $x \in G'_0$ and $\eta$ is a class function of $G'$, 
then 
$$\eta^G(x)=\eta(x) \, .$$
 
   By \cite[Theorem 9.8]{Da}, $\theta$ naturally
 corresponds to some irreducible character $\eta \in {\rm Irr}_A(S')$, which therefore has field of values
 $\Q(i)$. Now, $S$ only has two $A$-invariant irreducible characters with field of values $\Q(i)$, so using the inverse of Dade's natural correspondence 
from that theorem, 
necessarily $\eta$ is one of the two linear characters of $S'/C=R$.
  
By hypothesis $\langle x, S\rangle=G$, and so $x \in G_0$.
 Since $$G_0=\bigcup_{\tau \in T} (G'_0)^\tau \, ,$$
we may assume that $x \in G'_0$. By \cite[Lemma G]{Is2}, 
 $\psi(x)=\epsilon\gamma(x)$, where $\gamma \in \irr{S'A}$ is an extension of $\eta$, and $\epsilon$ is a root of unity. In particular,
 $\gamma$ is linear. We conclude that $|\psi(x)|=1$, as desired. 
  \end{proof}
  
\begin{lem}\label{action2}
Let $S \lhd G < \GL(V)\cong \GL_d(\C)$, where $G$ is finite, $S = \tw2 B_2(q)$, $D = q_0(q-1)$ with $q_0=2^n,\ q=2^{2n+1}$, 
and $V = \C^D$ is irreducible over $S$.
Suppose that $G$ contains two elements $g_0, g_1$ such that $\Trace(g_0) = \pm i$, $\Trace(g_1) = \pm 1$, and
$g_1 \in g_0S$. 
Suppose in addition that the trace of every element in $\ZB(G)$ belongs to $\Q(i)$. Then $g_0$ and $g_1$ both induce
non-inner automorphisms of $S$.
\end{lem}
 
\begin{proof}
Assume the contrary. As $g_1 \in g_0S$, $g_0$ induces an inner automorphism of $S$. Hence we can write 
$g_0 = zs$ for some $s \in S$ and some $z \in \CB_{G}(S) = \ZB(G)$. Since $G$ is finite, $z$ has finite order,
whence $z = \zeta \cdot \id_V$ for some root of unity $\zeta$. By assumption, $\zeta \in \Q(i)$, whence $\zeta^4=1$.

First suppose that $\zeta = \pm 1$. Then $\Trace(s) = \zeta^{-1} \Trace(g_0) = \pm i$, but this is impossible for any element in 
$S$, see \cite{Suz}. Hence, $\zeta = \pm i$. Recalling that $g_1 \in g_0S$, we can write $g_1=g_0t = zst$ for some
$t \in S$. In such a case, $\Trace(st) = \zeta^{-1}\Trace(g_1) = \pm i$, and we again arrive at a contradiction since 
$st \in S$.
\end{proof} 

\subsection*{9B. Traces of Frobenii and arithmetic monodromy groups}
 We now give a lemma on traces for those local systems $\sF(q,f)$ in which the polynomial $f(x)$ of degree $(q_0+1)t(q)$ lies in $\F_2[x]$ and has $f(0)=0$.
Recall that $q_0=2^n,\ n \ge 1, q=2q_0^2$, $t(q) =q+1-2q_0$; for $k/\F_2$ a finite extension and $i=\zeta_4$, the trace function of  $\sF(q,f)$ is
$$ t \in k \mapsto  \frac{-1}{\bigl(1-(-1)^n\,i\bigr)^{\deg(k/\F_2)}}\sum_{x \in k}\psi_2(\Trace_{W_2(k)/W_2(\F_2)}([x^{t(q)},f(x)+tx])).$$

\begin{lem}\label{01traces}For $\sF(q,f)$ as above, i.e.\ with $f(x)\in \F_2[x]$ of degree $(q_0+1)t(q)$ and $f(0)=0$,  define
$$A:= {\rm \ the \ number \ of\  nonzero\ monomials\ in\ }f(x),$$
so that $A \equiv f(1) \pmod*{2}$.
Then the traces at points of $\F_2$ are given as follows, with $i = \zeta_4$:
$$\begin{array}{c}
    \Trace\bigl(\Frob_{0,\F_2}|\sF(q,f)\bigr) = \left\{ \begin{array}{rl}-1,& {\rm if \ }2 \nmid (A-n),\\
    -(-1)^ni,& {\rm if \ }2|(A-n), \end{array} \right.\\
  \Trace\bigl(\Frob_{1,\F_2}|\sF(q,f)\bigr) = \left\{ \begin{array}{rl}-(-1)^ni,& {\rm \ if \ }2\nmid (A-n),\\
   -1,& {\rm  \ if\ }2|(A-n). \end{array}\right. \end{array}$$
Furthermore, for any finite extension $k/\F_2$ and any $t \in k$, $|\Trace(\Frob_{t,k}|\sF(q,f))| \leq \sqrt{\#k}$.   
\end{lem}
\begin{proof}The trace at time $0 \in \F_2$ is  $\frac{-1}{(1-(-1)^n\,i)}$ times the sum
$$\psi_2[0,0]+\psi_2[1,A] = 1+(-1)^Ai,$$
while the trace at time $1 \in \F_2$ is $\frac{-1}{(1-(-1)^n\,i)}$ times the sum
$$\psi_2[0,0]+\psi_2[1,A+1] = 1+(-1)^{A+1}i.$$
If $(-1)^A = -(-1)^n$, these traces are respectively $-1$ and $-(-1)^n\,i$. If $(-1)^A = (-1)^n$, these traces are respectively $-(-1)^ni$ and $-1$.
For the second statement, note that $\Trace(\Frob_{t,k}|\sF(q,f))$ is a sum of $\#k$ terms, each of absolute value $1$, divided by a 
clearing factor of absolute value $\sqrt{\#k}$, hence it has absolute value $\leq \sqrt{\#k}$.
\end{proof}

\begin{thm}\label{main2d}
Let $n \in \Z_{\geq 1}$ and $q=2^{2n+1}$. Suppose the Airy sheaf $\sF(q,f)$ defined in \eqref{airy-s}
with $q=2^{2n+1}$, $f(x) \in \F_2[x]$, and $f(0)=0$, has finite geometric monodromy group $S$. Then the following statements hold for 
$S$ and the arithmetic monodromy group $G:=G_{\ari,\F_2}$ of $\sF(q,f)$ over $\F_2$.
\begin{enumerate}[\rm(i)]
\item $S=\tw2 B_2(q)$. Furthermore, $G$  induces non-inner automorphisms of $S$,  
$G = \ZB(G) \times R$ for some $S < R \leq \Aut(S) = S \rtimes C_{2n+1}$, and $\ZB(G) \leq C_4$. 
\item If $2n+1$ is coprime to $|S|$, then $G \cong \ZB(G) \times \Aut(S)$. 
\item Suppose that  $|\Trace(\Frob_{0,\F_{q^2}})|=q_0:= 2^n$.
Then $G \cong \ZB(G) \times \Aut(S)$. Moreover, $|\ZB(G)| \leq 2$ if $\Trace(\Frob_{0,\F_{q^2}}) = -q_0$ or 
$\Trace(\Frob_{0,\F_{q}}) =  \pm q_0 i$, and 
$\ZB(G) = 1$ if $\Frob_{1,\F_2}$ has odd order.
\end{enumerate}
\end{thm} 
 
\begin{proof}
(i) By Theorem \ref{main2}, the geometric monodromy group of $\sF(q,f)$ is $S = \tw2 B_2(q)$ in its irreducible representation 
of degree $D=q_0(q-1)$. By Lemma \ref{01traces}, for the images $g_0$ of $\Frob_{0,\F_2}$ and $g_1$ of $\Frob_{1,\F_2}$ in
$G$, one has trace $\pm i$ and the other has trace $\pm 1$. Note that $g_1 \in g_0S$, and since $G$ is finite,
the traces of all elements in $G$ belong to $\Q(i)$ by Chebotarev density. 
By Lemma \ref{action2}, $g_0$ induces a non-inner automorphism of $S$. 

Recall that $\Out(S) \cong C_{2n+1}$. Then $G$ projects onto a subgroup $S \rtimes C_m$ of $\Aut(S)$ with kernel
$\CB_G(S) = \ZB(G)$, for some divisor $m > 1$ of $2n+1$. 
Again, each element in $\ZB(G)$ acts as a scalar $\al$ on $\C^D$ and has trace a root of unity belonging to 
$\Q(i)$, whence $\al^4=1$, and thus $\ZB(G) \leq C_4$. Now $(G/S)/(\ZB(G)S/S) \cong G/\ZB(G)S \cong C_m$ is cyclic, and 
$\ZB(G)/S \leq \ZB(G/S)$. It follows that $G/S$ is abelian of order $m|\ZB(G)|$, with cyclic quotient of odd order $m$ and 
a cyclic $2$-subgroup $\ZB(G)S/S$ of order $|\ZB(G)|$. Hence $G/S = R/S \times \ZB(G)/S$, with $R \cong S \rtimes C_m$.
The composition factors of $R$ are $S$ and cyclic groups of odd order (dividing $m$), so $R \cap \ZB(G)=1$ and 
$G = \ZB(G) \times R$.

\smallskip
(ii) Now assume that $2n+1$ is coprime to $|S|$, but $m < 2n+1$. Let $\psi$ be the character of $R = S \rtimes C_m$ 
acting on the sheaf, which extends the character $\theta$ of $S$. Let $\eta$ denote the Glauberman correspondent of $\theta$
as in Lemma \ref{coprime2}, in particular, it has degree $(r-1)\sqrt{r/2}$ for $r:=2^{(2n+1)/m} \geq 8$. 
We can write $g_0 = zh_0$ and $g_1 = zh_1$ for some $h_0,h_1 \in R$ and $z \in \ZB(G)$
(recall $g_0S=g_1S$). Now $z$ acts as a root of unity $\beta$ with $\beta^4=1$, and $\psi(h_j) = \gamma_j\eta(c_j)$ with
$c_j \in C = \tw2 B_2(r)$ and $\gamma_j^{2m}=1$ for $j= 0,1$ by Lemma \ref{coprime2}. Since $\eta(c_j),\beta,\psi(g_j) \in \Q(i)$, 
we see that $\gamma_j \in \Q(i)$. But $\gamma_j^{2m}=1$ and $2 \nmid m$, so $\gamma_j = \pm 1$.

Note that, since $r \geq 8$, the character $\eta$ does not take value $\pm i$, see \cite{Suz}.
Without loss we may assume $\psi(g_0) = \pm 1$ and $\psi(g_1) = \pm i$. Now, if $\beta = \pm i$, then 
$\eta(g_0) = \psi(g_0)/(\beta\gamma_0) = \pm i$, which is impossible.  
On the other hand, if $\beta = \pm 1$, then 
$\eta(g_1) = \psi(g_1)/(\beta\gamma_1) = \pm i$, which is again impossible. Hence $m=2n+1$, as stated. 

\smallskip
(iii) By assumption, $|\varphi(g_0^{4n+2})| = q_0$, where $\varphi$ is the 
character of $G$ acting on $\sF(q,f)$. Using $G = \ZB(G) \times R$, we again write $g_0 = zh_0$
with $z \in \ZB(G)$ acting on $\sF(q,f)$ as a root of unity $\beta$, and $h_0 \in R$. Then $\ZB(G) = \langle z \rangle$ since  
$G = \langle g_0,S \rangle$, and $\beta^4=1$. Now $s:=h_0^{2n+1} \in S$ as $R/S \hookrightarrow C_{2n+1}$,  
and $|\varphi(s^2)| = q_0$. Checking the character table of $S$ \cite{Bur}, we see that
$|s^2|=2$ or $4$ and thus $s$ is a $2$-element. But every $2$-element of $S$ has order dividing $4$, so $|s|=4$. 
Hence we can write $|h_0| = 4e$ for some $e|(2n+1)$. 

Suppose that $e < 2n+1$, whence $e \leq (2n+1)/3$. Then $g_0^e$ is the image of $\Frob_{0,2^e}$, and so
$$q_0 =2^n > 2^{(2n+1)/6} \geq 2^{e/2} \geq |\varphi(g_0^e)| = |\varphi(h_0^e)|$$ 
by Lemma \ref{01traces}. 
On the other hand, since $R/S \hookrightarrow C_{2n+1}$ and $h_0^e$ has order $4$,
$h_0^e \in S$ and hence $|\varphi(h_0^e)| = q_0$, a contradiction. Thus $|h_0| = 4(2n+1)$. Note that $|\CB_S(s^2)| = q^2$, hence 
$|\CB_G(s^2)|$ has order dividing $q^2|\ZB(G)| \cdot |R/S|$. But $h_0$ belongs to $\CB_G(s^2)$ and has order $4(2n+1)$; thus
$(2n+1)$ divides $|R/S|$ and so $R=\Aut(S)$. 

Since $G=\langle g_1,S \rangle$, the assumption that the image $g_1$ of $\Frob_{1,\F_2}$ has odd order implies 
that $2 \nmid |G/S|$, and so $\ZB(G)=1$.  Next suppose that $|\ZB(G)| > 2$ but $\varphi(g_0^{4n+2}) = -q_0$ or 
$\varphi(g_0^{2n+1}) = \pm q_0i$. Then $\beta = \pm i$. In the former case,   
$$\varphi(g_0^{4n+2}) = \varphi(\beta^{4n+2}h_0^{4n+2}) = -\varphi(h_0^{4n+2}) = -\varphi(s^2),$$ 
and thus the involution $s^2 \in S$ has 
trace $q_0$, which is impossible, cf. \cite{Bur}. In the latter case,   
$$\varphi(g_0^{2n+1}) = \varphi(\beta^{2n+1}h_0^{4n+2}) =  \pm i \varphi(h_0^{2n+1}) = \pm i \varphi(s),$$
and thus $s \in S$ has 
trace $\pm q_0$, which is again impossible, cf. \cite{Bur}. 
\end{proof} 

\begin{prop}\label{traces-f0a}
We consider the sheaf $\sF_q$, $q=2^{2n+1}$, as defined in the Introduction. Thus $\sF_q:=\sF(q,f)$ with $f(x):=f_1(x^{t(q)}), \ f_1(x):=\sum_{i=1}^n x^{1+2^i}$ as in 
\eqref{airy-s}. For a finite extension $k/\F_2$, define
$$\Ker(k):=\bigl{\{}x \in k\  \mid\  \sum_{i=0}^{2n}x^{2^i} =0\bigr{\}}.$$
Then we have the following results.
\begin{itemize}
\item[(i)]For any subfield $k$ of $\F_{q^2}$, $|\Trace(\Frob_{0,k}|\sF_q)|^2$ is either $0$ or $\#\Ker(k)$. 
\item[(ii)]$|\Trace(\Frob_{0,\F_{q^2}}|\sF_q)|^2=\#\Ker(\F_{q^2})=\#\Ker(\F_q)=q/2$.
\end{itemize}
\end{prop}

\begin{proof}
We first observe that $\gcd(t(q),q^2-1)=1$. To see this, note that $t(q)=q+1-2q_0$ divides $q^2+1$ (indeed $(q+1-2q_0)(q+1+2q_0)=q^2+1$),
while $\gcd(q^2-1,q^2+1) =\gcd(q^2-1,2)=1$. Thus for any subfield $k$ of of $\F_{q^2}$, the map $x \mapsto x^{t(q)}$ is bijective on $k$. 

The sheaf $\sF_q$ was built out of the Witt vector 
$$\bigl[x^{t(q)},\sum_{i=1}^n x^{t(q)(1+2^i)}\bigr].$$
Let us denote by $\sH_q$ the sheaf built by the same recipe, with same clearing factor, out of the Witt vector
$$\bigl[x,\sum_{i=1}^n x^{1+2^i}\bigr].$$
Then for any subfield $k$ of of $\F_{q^2}$
\begin{equation}\label{fq-hq1}
  \Trace(\Frob_{0,k} |\sF_q) =\Trace(\Frob_{0,k}|\sH_q),
\end{equation}  
precisely because the map $x \mapsto x^{t(q)}$ is bijective on $k$. 

Let us rewrite the input Witt vector for $\sH_q$ as
\begin{equation}\label{fq-hq2}
  \bigl[x,\sum_{i=1}^n x^{1+2^i}\bigr] =[x,xR(x)]\ {\rm\ with\ }R(x):=\sum_{i=1}^n x^{2^i}.
\end{equation}  
With this rewriting, we apply the idea of van der Geer-van der Vlugt, cf. \cite[\S5]{vdG-vdV}, as follows.
Let us define
\begin{equation}\label{fq-hq3}
  V(x):=[x,xR(x)].
\end{equation}  
In Witt vector addition in $\F_2$-algebras, using the fact that $R(x)$ is an additive polynomial, we get
$$\begin{aligned}V(x+y)-V(x)-V(y) & = [x+y,(x+y)(R(x)+R(y))] +[x,xR(x)+x^2] +[y,yR(y)+y^2]\\
 & =[y,(x+y)x+(x+y)(R(x)+R(y)) + xR(x)+x^2] +[y,yR(y)+y^2]\\
 & =[0, y^2 + (x+y)x +(x+y)(R(x)+R(y)) + xR(x)+x^2+yR(y)+y^2]\\
& = [0, xy+xR(y)+yR(x)]\\
 &  =[0,\langle x,y \rangle]\end{aligned}$$
for 
\begin{equation}\label{fq-hq4}
  \langle x,y \rangle:=xy+xR(y)+yR(x).
\end{equation}  
The key point is that $\langle x,y \rangle$ on $k\times k$ is a symmetric $\F_2$-bilinear map to $k$, and $\Trace_{k/\F_2}(\langle x,y \rangle)$ is a symmetric $\F_2$-bilinear form on $k \times k$ as $\F_2$ vector space.

Then 
$$|\Trace(\Frob_{0,k}|\sH_q)|^2 = (1/\#k)\sum_{x,y \in k}\psi_2\bigl(\Trace_{k/\F_2}(V(x)-V(y))\bigr)$$
(by the shearing transformation $(x,y)\mapsto (x+y,y)$)
$$= (1/\#k)\sum_{x,y \in k}\psi_2\bigl(\Trace_{k/\F_2}(V(x+y)-V(y))\bigr)$$
$$=(1/\#k)\sum_{x,y \in k}\psi_2\bigl(\Trace_{k/\F_2}(V(x) +[0,\langle x,y \rangle])\bigr)$$
$$ =\sum_{x \in k}\psi_2\bigl(\Trace_{k/\F_2}(V(x))\bigr)\biggl((1/\#k)\sum_{y \in k}\psi\bigl(\Trace_{k/\F_2}(\langle x,y \rangle)\bigr)\biggr).$$
The second summand vanishes unless the given $x \in k$ has $\Trace_{k/\F_2}(\langle x,y \rangle)=0$ for all $y \in k$, in which case it is $1$.
But $x \in k$ has $\Trace_{k/\F_2}(\langle x,y \rangle)=0$ for all $y \in k$ if and only if $x \in \Ker(k)$. To see this, note that for $x,y
\in k$, 
$$\langle x,y \rangle= xy +xR(y)+yR(x)=xy + \sum_{i=1}^n xy^{2^i} +\sum_{i=1}^n yx^{2^i} $$
has the same $\Trace_{k/\F_2}$ as
$\bigl(x +\sum_{i=1}^n x^{1/2^i} + \sum_{i=1}^n x^{2^i}\bigr)y$.
So by nondegeneracy of the trace, $x \in k$ has $\Trace_{k/\F_2}(\langle x,y \rangle)=0$ for all $y \in k$ if and only if
$$ x +\sum_{i=1}^n x^{1/2^i} + \sum_{i=1}^n x^{2^i} =0,\ {\rm i.e.\ if\ and\ only\ if\ }\sum_{i=0}^{2n}x^{2^i}=0.$$
Thus for $k$ a subfield of $\F_{q^2}$, 
$$\bigl{|}\Trace(\Frob_{0,k}\sH_q)\bigr{|}^2=\sum_{x \in \Ker(k)}\psi_2\bigl(\Trace_{k/\F_2}(V(x))\bigr).$$

To show (i), notice that on $\Ker(k)$, $x \mapsto \Trace_{k/\F_2}(V(x))$ is additive, i.e.\ it is a linear form on $\Ker(k)$. If it is nontrivial,
the sum giving $|\Trace(\Frob_{0,k}\sH_q)|^2$ vanishes. If it is trivial, this sum is $\#\Ker(k)$.

To show (ii), notice that for any $k \subseteq \overline{\F_2}$, $\Ker(k)$ is precisely the set of elements $x \in k \cap \F_q$ with $\Trace_{\F_q/\F_2}(x)=0$. (Indeed, 
if $x \in \Ker(k)$, then 
$$0=\sum^{2n}_{i=0}F^i(x) = F\bigl(\sum^{2n}_{i=0}F^i(x)\bigr) = \sum^{2n+1}_{i=1}F^i(x),$$
and so $x=F^{2n+1}(x)$, i.e. $x \in \F_q$.)
In particular, $\Ker(\F_{q^2}) = \Ker(\F_q)$ and $\#\Ker(\F_q)=q/2$. Now for $x \in \Ker(\F_{q^2} )$ we have
$$\begin{aligned}\Trace_{\F_{q^2}/\F_2}(V(x))& =\Trace_{\F_{q}/\F_2}\bigl(\Trace_{\F_{q^2}/\F_q}(V(x)))\\
  & = \Trace_{\F_{q}/\F_2}(V(x)+V(x))\\ & = \Trace_{\F_{q}/\F_2}([0,x^2])\\
  & = [0,\Trace_{\F_{q}/\F_2}(x^2)]\\ & =[0,0],\end{aligned}$$
precisely because every $x \in \Ker(\F_{q^2})$ is an element of $\F_q$ of trace zero. So we see directly that each of the summands in the sum giving $|\Trace(\Frob_{0,\F_{q^2}}|\sH_q)|^2$ is simply $1$.
%
\end{proof}

\begin{prop}\label{traces-f0b}
For $q=2^{2n+1}$ and the Airy sheaf $\sF_q$, 
$$\Trace(\Frob_{0,\F_q}|\sF_q)=-\eps2^ni,$$
where $\epsilon := (-1)^{n(n+1)/2}$ is the Jacobi symbol 
$$\eps_{2n+1}=\biggl(\frac{2}{2n+1}\biggr)= 
\left\{ \begin{array}{rl} 1, & {\rm if\ }2n+1 \equiv \pm 1 \pmod* 8,\\ -1, &{\rm if\ }2n+1 \equiv \pm 3 \pmod*{8}.\end{array}\right.$$
\end{prop}
\begin{proof}
(i) Proposition \ref{traces-f0a} implies that $\Trace(\Frob_{0,\F_q}|\sF_q)$ lies in $\Z[i]$ (because the only possible non-integrality is at the unique place of $\Q(i)$ over $2$, where this trace and its complex conjugate have the same $2$-adic ord). Furthermore, we can work
with traces over $\sH_q$ instead of $\sF_q$.

Let us denote by $F$ the absolute Frobenius $x \mapsto x^2$, and define
$$R_n(x):=\sum_{i=1}^n F^i(x),\ \ V_n(x):=[x,xR_n(x)]=[x,x\bigl(\sum_{i=1}^n F^i(x)\bigr)].$$
Consider the non-normalized sum
$$\RawTrace(\Frob_{0,\F_q}|\sF_q):=-\sum_{x \in \F_q}\psi_2\bigl(\Trace_{\F_q/\F_2}(V_n(x))\bigr),$$
so that
$$\RawTrace(\Frob_{0,\F_q}|\sF_q)=(1-(-1)^ni)^{2n+1}\Trace(\Frob_{0,\F_q}|\sF_q).$$

\smallskip
(ii) To explain the idea of the proof, consider first the case when $2n+1$ is an odd prime $p$. Then $\F_{2^p}/\F_2$ has degree $p$, and  $\F_{2^p} \smallsetminus \F_2$ is the disjoint union of $F$-orbits of length $p$. On each such $F$-orbit, the value of $\Trace_{\F_q/\F_2}(V_n(x))$ is constant, and this constant value is then repeated $p$ times as we sum over this orbit. So we have a congruence modulo $p\Z[i]$:
$$\begin{aligned}\RawTrace(\Frob_{0,\F_q}|\sF_q) & \equiv -\sum_{x\in \F_2}\psi_2\bigl(\Trace_{\F_q/\F_2}(V_n(x))\bigr)\\
   & \equiv-\psi_2\bigl(\Trace_{\F_q/\F_2}(V_n(0))\bigr)-\psi_2\bigl(\Trace_{\F_q/\F_2}(V_n(1))\bigr)\\
   & \equiv-\psi_2\bigl(\Trace_{\F_q/\F_2}[0,0]\bigr)-\psi_2\bigl(\Trace_{\F_q/\F_2}[1,n]\bigr)\end{aligned}$$
(remembering that $\F_q/\F_2$ has odd degree $p$, and both $V_n(0), V_n(1)$ are already $\F_2$-rational)
$$\equiv-\psi_2(p[0,0])-\psi_2(p[1,(p-1)/2])=-1 -\psi_2([1,(p-1)/2])^p= -1 -(i^{1+(p-1)})^p=-1-i^{p^2}=-1-i.$$
So when $2n+1=p$, we have a congruence modulo $p\Z[i]$:
$$(1-(-1)^{(p-1)/2}i)^p\Trace(\Frob_{0,\F_q}|\sF_q) \equiv-1-i.$$
Multiplying both sides by $(1+(-1)^{(p-1)/2}i)^p$, we get a congruence modulo $p\Z[i]$:
$$2^p\Trace(\Frob_{0,\F_q}|\sF_q) \equiv -(1+i)(1+(-1)^{(p-1)/2}i)^p.$$
If $p \equiv 1 \pmod* 4$, the right side  is 
$$-(1+i)^{p+1}=-(2i)(2i)^{(p-1)/2}=-(2i)2^{(p-1)/2}i^{(p-1)/2}.$$
If  $p \equiv 3 \pmod* 4$, the right side  is 
$$-(1+i)(1-i)^p=-2(1-i)^{p-1}=-2(-2i)^{(p-1)/2}=2^{(p+1)/2}i^{(p-1)/2}.$$
Recalling that $2^p \equiv 2 \pmod* p$, and that $(p-1)/2=n$, we get a congruence modulo $p\Z[i]$:
$$
\Trace(\Frob_{0,\F_q}|\sF_q) \equiv 
  \left\{ \begin{array}{ll} -i2^{n}i^{(p-1)/2},& {\rm if\ }p \equiv 1 \pmod* 4,\\ 2^{n}i^{(p-1)/2}, & {\rm if\ }p \equiv 3 \pmod* 4.\end{array}\right.
$$ 
Thus
\begin{equation}\label{eq-124}
  \Trace(\Frob_{0,\F_q}|\sF_q) \equiv -\epsilon_p2^ni,
\end{equation}
where $\epsilon_p$ is given by $i^{(p-1)/2}=(-1)^{(p-1)/4} = (-1)^{(p^2-1)/8}$ when $p \equiv 1  \pmod* 4$, and by 
$i^{(p+1)/2} = (-1)^{(p+1)/4} = (-1)^{(p^2-1)/8}$ when $p \equiv 3 \pmod* 4$. Thus in both cases 
$\epsilon_p$ is the Legendre symbol $\bigl(\frac{2}{p}\bigr)=(-1)^{(p^2-1)/8}$.

In view of Proposition \ref{traces-f0a}, $\Trace(\Frob_{0,\F_q}|\sF_q)$ is either $0$ or an element of $\Z[i]$ of absolute value $2^n$.
It cannot be $0$ because of \eqref{eq-124}. So it must be one of $\pm 2^n$ or $\pm 2^ni$. Of
these four possibilities, only $-\epsilon_p2^ni$ is congruent modulo  $p\Z[i]$ to $-\epsilon_p2^ni$ (this is just the statement that for an odd prime $p$,
the four powers of $i$ are distinct modulo  $p\Z[i]$).

\smallskip
(iii) We now turn to the general case, where we proceed by induction on the total number (counting multiplicity) of primes dividing $2n+1$. Thus we write
$$2n+1 =ps,\ s=2a+1,\ p=2b+1, \ {\rm with\ } a,b \ge 1 \ {\rm and\ }p \ {\rm prime}.$$
We will need to deal with both $\sF_{2^{ps}}$ and $\sF_{2^s}$. To simplify notation, let us write
$$Q:=2^{ps},~q':=2^s.$$
Here 
$$ps =s(2b+1)=2bs+s = 2(sb+a)+1.$$

The idea is that $\F_Q \smallsetminus \F_{q'}$ is the disjoint union of $F^s$-orbits, each of length $p$, and on each of these orbits, the
value of  $\Trace_{\F_Q/\F_2}(V_{sb+a}(x))$ is constant, and this constant value is then repeated $p$ times as we sum over this orbit. 
Thus we get a congruence modulo $p\Z[i]$:
$$\RawTrace(\Frob_{0,\F_Q}|\sF_Q) \equiv -\sum_{x \in \F_{q'}}\psi_2\bigl(\Trace_{\F_Q/\F_2}(V_{sb+a}(x))\bigr).$$
But for $x \in \F_{q'}$, the trace from $\F_Q$ down to $\F_{q'}$ is just multiplication by $p$, so 
$$\Trace_{\F_Q/\F_2}(V_{sb+a}(x))=\Trace_{\F_{q'}/\F_2}(pV_{sb+a}(x)).$$

Suppose first that $b$ is even, i.e.\ that $p \equiv 1 \pmod* 4$. Recall $4[x,y]=0$. So for $x \in \F_{q'}$,
$$\begin{aligned}pV_{sb+a}(x)& =V_{sb+a}(x)=[x,xR_a(x)+x\bigl(\sum_{i=a+1}^{sb+a}F^i(x)\bigr)]\\
    & =[x,xR_a(x) +x\bigl(b\Trace_{\F_{q'}/\F_2}(x)\bigr)]=[x,xR_a(x)]=V_a(x),\end{aligned}$$
where the last equality holds because $b$ is even. So in this even $b$ case, 
$$\RawTrace(\Frob_{0,\F_Q}|\sF_Q) \equiv \RawTrace(\Frob_{0,\F_a}|\sF_{q'}) \pmod*{p\Z[i]}.$$

Suppose now that $b$ is odd, i.e.\ that $p \equiv -1 \pmod* 4$. Then for $x \in \F_{q'}$, 
$$\begin{aligned}pV_{sb+a}(x)& =-V_{sb+a}(x)=-[x,xR_a(x) +x\bigl(\sum_{i=a+1}^{sb+a}F^i(x)\bigr)]\\
   & =-[x,xR_a(x) +x\Trace_{\F_{q'}/\F_2}(x)] = [x,xR_a(x)+x^2 +x\Trace_{\F_{q'}/\F_2}(x)]\\
   & = [x,xR_a(x)]+[0,x^2 +x\Trace_{\F_{q'}/\F_2}(x)] = V_a(x) +[0,x^2 +x\Trace_{\F_{q'}/\F_2}(x)].\end{aligned}$$
But the term $[0,x^2 +x\Trace_{\F_{q'}/\F_2}(x)]$ has  
$$\Trace_{\F_{q'}/\F_2}\bigl([0,x^2 +x\Trace_{\F_{q'}/\F_2}(x)]\bigr)=[0,\Trace_{\F_{q'}/\F_2}\bigl(x^2 +x\Trace_{\F_{q'}/\F_2}(x)\bigr)]=[0,0],$$
where the last equality holds because $\Trace_{\F_{q'}/\F_2}(x^2) =\Trace_{\F_{q'}/\F_2}(x)= (\Trace_{\F_{q'}/\F_2}(x))^2.$
So in this odd $b$ case as well, also 
$$\RawTrace(\Frob_{0,\F_Q}|\sF_Q) \equiv \RawTrace(\Frob_{0,\F_a}|\sF_{q'}) \pmod*{p\Z[i]}.$$

We have shown that
\begin{equation}\label{eq-tr20}
  (1-(-1)^{sb+a}i)^{ps}\Trace(\Frob_{0,\F_Q}|\sF_Q) \equiv (1-(-1)^{a}i)^{s}\Trace(\Frob_{0,\F_{q'}}|\sF_q) \pmod*{p\Z[i]}.
\end{equation}  
By the induction hypothesis,
\begin{equation}\label{eq-tr21}
  \Trace(\Frob_{0,\F_{q'}}|\sF_{q'})= -\eps_s2^ai.
\end{equation}  

The clearing factors are invertible modulo $p\Z[i]$. We next show that the clearing factors are equal modulo $p\Z[i]$.  Their ratio is
$$\frac{(1-(-1)^{sb+a}i)^{ps}}{(1-(-1)^{a}i)^{s}}=\frac{(1-(-1)^{sb+a}i)^{ps}(1+(-1)^{a}i)^{s}}{2^s}.$$
If $b$ is odd, then $(p+1)/2$ is even, and  this ratio is
$$\frac{(1+(-1)^{a}i)^{ps+s}}{2^s}=\frac{(2(-1)^ai)^{s(p+1)/2}}{2^s}= 2^{s(p-1)/2}((-1)^ai)^{s(p+1)/2}=                                                                 2^{s(p-1)/2}i^{s(p+1)/2}.$$
If $b =(p-1)/2$ is even, this ratio is 
$$(1-(-1)^{a}i)^{(p-1)s} =(-2(-1)^{a}i)^{s(p-1)/2} =2^{s(p-1)/2}i^{s(p-1)/2}.$$
Let $\chi_{p,{\rm quad}}$ be the quadratic character of $\F_p^\times$, so 
$2^{(p-1)/2} \equiv \chi_{p,{\rm quad}}(2) \pmod*{p\Z[i]}$.
As $s$ is odd, the ratio of clearing factors modulo $p\Z[i]$ is  $\chi_{p,{\rm quad}}(2)i^{s(p-1)/2}$ when $p \equiv 1 \pmod* 4$, and 
$\chi_{p,{\rm quad}}(2)i^{s(p+1)/2}$ when $p \equiv 3 \pmod* 4$. Next, when $p \equiv 1 \pmod* 4$, 
$$i^{(p-1)/2} = \chi_{p,{\rm quad}}(2).$$
When  $p \equiv 3 \pmod* 4$, 
$$i^{(p+1)/2}= \chi_{p,{\rm quad}}(2).$$
As $s$ is odd, we find that in all cases the ratio of clearing factors is $1 \pmod*{p\Z[i]}$, as stated.                                  
Hence \eqref{eq-tr20} and \eqref{eq-tr21} imply the congruence
$$\Trace(\Frob_{0,\F_Q}|\sF_Q) \equiv \Trace(\Frob_{0,\F_{q'}}|\sF_{q'})  = -\epsilon_s2^ai \pmod*{p\Z[i]}.$$
Note that
$$2^{sb}=(2^{(p-1)/2})^s \equiv \chi_{p,{\rm quad}}(2)^s = \chi_{p,{\rm quad}}(2)  = \eps_p \pmod*{p}.$$ 
Hence
$$\eps_{ps}2^{sb+a} \equiv \eps_{ps}\eps_p2^a = \eps_s2^a \pmod*{p}$$
by multiplicativity of the Jacobi symbol.  It follows that 
$$\Trace(\Frob_{0,\F_Q}|\sF_Q) \equiv -\epsilon_{ps}2^{sb+a}i \pmod*{p\Z[i]}.$$
This congruence shows that $\Trace(\Frob_{0,\F_Q}|\sF_Q)$ is nonzero, so by Proposition \ref{traces-f0a} it is one of $\pm 2^{sb+a}$ or $\pm 2^{sb+a}i$. Again, of these four possibilities, only $-\epsilon_{ps}2^ni$ is congruent modulo  $p\Z[i]$ to $-\epsilon_{ps}2^ni$, 
and the induction step is complete.
\end{proof}

\begin{prop}\label{absvaluetraces-f1}
For $q=2^{2n+1}$, consider the Airy sheaf $\sF_q$. 
Then for any subfield $k$ of $\F_q$ 
$$|\Trace(\Frob_{1,k}|\sF_q)|^2 =1.$$
\end{prop}

\begin{proof}
Let $k$ be a subfield of $\F_{q^2}$. Note that $N:=(q+2q_0+1)q^2/2$ divides $(q^2+1)q^2$ and so is coprime to $q^2-1$. Hence the map $x \to x^N$ is a 
bijection on $k$, and it sends $x^{t(q)}$ to $x^{(q^2+1)q^2/2} = x$ for any $x \in k$.
As in the proof of  Proposition \ref{traces-f0b}, let us denote by $F$ the absolute Frobenius $x \mapsto x^2$.
For each integer $j \ge 0$, we define
$$R_j(x):=\sum^j_{i=1}F^i(x).$$
Consider the non-normalized sum
$$\RawTrace(\Frob_{1,k}|\sF_q):=-\sum_{x \in k}\psi_2\bigl(\Trace_{k/\F_2}(V(x))\bigr),$$
where $V(x)=[x,xR_n(x)+x^N]$. Then we have
$$\RawTrace(\Frob_{1,k}|\sF_q)=(1-(-1)^ni)^{\deg(k/\F_2)}\Trace(\Frob_{1,k}|\sF_q).$$

We now take $k$ to be a subfield of $\F_q$. We examine the function $x \mapsto x^N$ on $k$. This function depends only on $
N \pmod*{(q-1)}$. Then
$$N \equiv (1+1+2q_0)q/2=(2+2q_0)q/2 =q+qq_0 \equiv 1+q_0=1+2^n \pmod*{(q-1)}.$$
Thus if $x \in k$, then $x^N=xF^n(x)$,  and hence
$$V(x) = [x,xR_n(x)+x^N]=[x,xR_n(x)+xF^n(x)]=[x,xR_{n-1}(x)].$$

At this point, we repeat the van der Geer-van der Vlugt argument of Theorem \ref{traces-f0a}.  We find that if $k$ is a subfield of $\F_q$, then
$$|\Trace(\Frob_{1,k}|\sF_q)|^2=\sum_{x \in \Ker'(k)}\psi_2\bigl(\Trace_{k/\F_2}(V(x))\bigr),$$
with 
$$\Ker'(k):=\bigl{\{}x \in k| \sum_{i=0}^{2n-2}F^i(x)=0 \bigr{\}}.$$

The key observation is that for $k$ a subfield of $\F_q$, $\Ker'(k)=\{0\}.$ Indeed, since 
$$0 \in \Ker'(k) \subseteq \Ker'(\F_q),$$ 
it suffices to show that $\Ker'(\F_q)=\{0\}.$ But for $x\in \Ker'(\F_q)$, 
$$\Trace_{\F_q/\F_2}(x) = \bigr(\sum_{i=0}^{2n-2}F^i(x)\bigl) +F^{2n-1}(x) +F^{2n}(x)=F^{2n-1}(x) +F^{2n}(x).$$
Thus for $x\in \Ker'(\F_q)$, $\Trace_{\F_q/\F_2}(x) =F^{2n-1}(x+F(x))$. As $\Trace_{\F_q/\F_2}(x) \in \F_2$, this gives
$$\Trace_{\F_q/\F_2}(x) =x+x^2.$$
If $\Trace_{\F_q/\F_2}(x) =0$, then $x+x^2=0$, i.e.\ $x \in \F_2$, so $x$ is $0$ or $1$. Of these, only $x=0$ has $\Trace_{\F_q/\F_2}(x) =0$.
If $\Trace_{\F_q/\F_2}(x) =1$, then $x+x^2=1$ and so $\F_2(x)=\F_4$. But $\F_4$ is not a subfield of $\F_q$, which has odd degree $2n+1$ over $\F_2$. So in this second case, there are no possible $x$.

Thus
$$\bigl{|}\Trace(\Frob_{1,k}|\sF_q)\bigr{|}^2=\sum_{x \in \Ker'(k)}\psi_2\bigl(\Trace_{k/\F_2}(V(x))\bigr)=\psi_2(\Trace_{k/\F_2}(V(0)))=\psi_2([0,0])=1.$$
\end{proof}

\begin{prop}\label{finaltraces-f1}
For $q=2^{2n+1}$ and the Airy sheaf $\sF_q$, 
$$\Trace(\Frob_{1,\F_q}|\sF_q)=-1.$$
\end{prop}

\begin{proof}
(i) The quantity $\Trace(\Frob_{1,k}|\sF_q)$ lies in $\Q(i)$ and is integral outside the unique place over $2$, so by Proposition \ref{absvaluetraces-f1}, $\Trace(\Frob_{1,\F_q}|\sF_q)$ is one of $\{1,-1,i,-i\}$. For any odd prime $p$, these four elements are distinct 
modulo $p\Z[i]$.
We proceed by induction on the total number (counting multiplicity) of primes in the factorization of $2n+1$. 

For the induction base, suppose that $k=\F_q=\F_{2^p}$. Then $\F_q \smallsetminus \F_2$ is the disjoint union of $F$-orbits of length 
divisible by $p$. On each such $F$-orbit, the value of $\Trace_{\F_q/\F_2}(V(x))$ is constant, and this constant value is then repeated a multiple of $p$ times as we sum over this orbit. So we have a congruence modulo $p\Z[i]$:
$$\begin{aligned}\RawTrace(\Frob_{1,\F_q}|\sF_q) & \equiv -\sum_{x\in \F_2}\psi_2\bigl(\Trace_{\F_q/\F_2}(V(x))\bigr)\\
   & \equiv-\psi_2\bigl(\Trace_{\F_q/\F_2}(V(0))\bigr)-\psi_2\bigl(\Trace_{\F_q/\F_2}(V(1))\bigr)\\
   & \equiv-\psi_2\bigl(\Trace_{\F_q/\F_2}[0,0]\bigr)-\psi_2\bigl(\Trace_{\F_q/\F_2}[1,n+1]\bigr)\\
  & \equiv-\psi_2(s[0,0])-\psi_2(p[1,(p+1)/2])=-1 -\psi_2([1,(p+1)/2])^p\\
  & = -1 -(i^{1+(p+1)})^p=-1-i^{p^2+2p}=-1+i.\end{aligned}$$
So for $\Frob_{1,\F_q}$, we have a congruence modulo $p\Z[i]$:
$$(1-(-1)^{(p-1)/2}i)^p\Trace(\Frob_{1,\F_q}|\sF_q) = \RawTrace(\Frob_{1,\F_q}|\sF_q) \equiv-1+i.$$
Multiplying both sides by $(1+(-1)^{(p-1)/2}i)^p$, we get a congruence modulo $p\Z[i]$:
$$2^p\Trace(\Frob_{0,\F_q}|\sF_q) \equiv -(1-i)(1+(-1)^{(p-1)/2}i)^p.$$
If $p \equiv 3 \pmod* 4$, the right side  is 
$$-(1-i)^{p+1}=-(2i)^{(p+1)/2}=-2^{(p+1)/2}(-1)^{(p+1)/4} \equiv -2(-1)^{(p^2-1)/8+(p+1)/4} = -2 \pmod*{p}.$$
If  $p \equiv 1 \pmod* 4$, the right side  is $-(1-i)(1+i)^p$, which is 
$$-2(1+i)^{s-1}=-2(2i)^{(p-1)/2}=-2^{(p+1)/2}(-1)^{(p-1)/4}  \equiv -2(-1)^{(p^2-1)/8+(p-1)/4} = -2 \pmod*{p}.$$
Thus in both cases $q\Trace(\Frob_{1,\F_q}|\sF_q) \equiv -2 \pmod*{p\Z[i]}$. Using Proposition \ref{absvaluetraces-f1}, we conclude
that $\Trace(\Frob_{1,\F_q}|\sF_q)=-1$ in this case.

\smallskip
(ii) Suppose now that $2n+1=ps$ with an odd prime $p$ and an odd integer $s$. Then just as in (iii) of the proof of Proposition \ref{traces-f0b} we write
$$s=2a+1,~p=2b+1,~ps = 2(sb+a)+1.$$
We will need to deal with both $\sF_{2^{ps}}$ and $\sF_{2^s}$. To simplify notation, let us write
$$Q:=2^{ps},~q':=2^s.$$
Then, just as in the proof of Proposition \ref{absvaluetraces-f1}, 
$$\RawTrace(\Frob_{1,\F_Q}|\sF_Q)= -\sum_{x \in \F_Q}\psi_2(\Trace_{\F_Q/\F_2}([x,R_{sb+a-1}(x)])).$$
The elements of $\F_Q \smallsetminus \F_{q'}$ fall into $F$-orbits of length $p$, and the elements from each of these orbits have the same 
$\Trace_{\F_Q/\F_2}([x,R_{sb+a-1}(x)])$.
So we get a congruence modulo $p\Z[i]$
$$\begin{aligned}\RawTrace(\Frob_{1,\F_Q}|\sF_Q) & \equiv -\sum_{x \in \F_{q'}}\psi_2(\Trace_{\F_Q/\F_2}([x,R_{sb+a-1}(x)]))\\
   & =-\sum_{x \in \F_{q'}}\psi_2(\Trace_{\F_{q'}/\F_2}(p[x,R_{sb+a-1}(x)])).\end{aligned}$$
Notice that for $x\in \F_{q'}$, 
$$R_{sb+a-1}(x)=R_{a-1}(x)+b\Trace_{\F_{q'}/F_2}(x).$$
Suppose first that $p \equiv 1 \pmod* 4$. Then $b$ is even, 
$$p[x,R_{sb+a-1}(x)]=[x,R_{sb+a-1}(x)]=[x,xR_{a-1}(x)],$$ and hence 
$$\RawTrace(\Frob_{1,\F_Q}|\sF_Q) \equiv \RawTrace(\Frob_{1,\F_{q'}}|\sF_{q'}) \pmod*{p\Z[i]}$$
when $p \equiv 1 \pmod* 4$.

Suppose next that $p \equiv -1 \pmod* 4$. Since $b$ is odd,
$$\begin{aligned}p[x,R_{sb+a-1}(x)] & = -[x,R_{sb+a-1}(x)]=[x,x^2 + R_{sb+a-1}(x)]\\
    & =[x,x^2+x\Trace_{\F_{q'}/\F_2}(x) +xR_{a-1}(x)]\\
   & =[x,xR_{a-1}(x)] + [0,x^2+x\Trace_{\F_{q'}/\F_2}(x)].\end{aligned}$$
But for $x \in \F_{q'}$, $\Trace_{\F_{q'}/\F_2}$ annihilates $x^2+x\Trace_{\F_{q'}/\F_2}(x)$, so again
$$\RawTrace(\Frob_{1,\F_Q}|\sF_Q) \equiv \RawTrace(\Frob_{1,\F_{q'}}|\sF_{q'}) \pmod*{p\Z[i]}$$
when $p \equiv -1 \pmod* 4$.

It remains only to show that the ratio of clearing factors is $1 \pmod*{p\Z[i]}$ in both cases.
When $p \equiv 1 \pmod* 4$, $(p-1)/2$ is even, and this ratio is
$$\begin{aligned}(1-(-1)^ai)^{ps}/(1-(-1)^ai)^s & =(1-(-1)^ai)^{(p-1)s} =(-2(-1)^ai)^{s((p-1)/2)}\\
   & = \bigr(2^{(p-1)/2}(-1)^{(p-1)/4}\bigr)^s \equiv 1 \pmod*{p},\end{aligned}$$
since $2^{(p-1)/2} \equiv (-1)^{(p^2-1)/8} \pmod*{p}$. When $p \equiv -1 \pmod* 4$, $(p+1)/2$ is even, and 
the ratio is
$$\begin{aligned}(1+(-1)^ai)^{ps}/(1-(-1)^ai)^s & =(1+(-1)^ai)^{ps}(1+(-1)^ai)^{s}/2^s = (1+(-1)^ai)^{(p+1)s}/2^s\\
 &  =(2(-1)^ai)^{s(p+1)/2}/2^s
=  \bigr(2^{(p-1)/2}(-1)^{(p+1)/4}\bigr)^s  \equiv 1 \pmod* p,\end{aligned}$$
again because $2^{(p-1)/2} \equiv (-1)^{(p^2-1)/8} \pmod*{p}$.
\end{proof}

\begin{lem}\label{order}
For $q=2^{2n+1}$ and $S := \tw2 B_2(q)$, suppose that $s \in \Aut(S)$ has odd order and that 
$\langle s,S \rangle = \Aut(S)$. Then $|s|$ divides $5(2n+1)$. 
\end{lem}

\begin{proof}
Let $\sG := \Sp_4(\overline{\F_2})$. It is well known that there is a Steinberg endomorphism $\sigma: \sG \to \sG$ such that 
$\sigma^2$ is the standard Frobenius map $(a_{ij}) \mapsto (a_{ij}^2)$ on $\sG$, and then we can identify $S$ with 
$\sG^{\sigma^{2n+1}} := \{ X \in \sG \mid \sigma^{2n+1}(X) = X \}$, which is then $\sigma$-invariant. Letting $\sigma$ also denote its
action on $S$, we can write $\Aut(S) = \langle \sigma,S \rangle$. Without loss, we may assume $Ss = S\sigma^{-1}$, and write
$s=x\sigma^{-1}$. By the Lang-Steinberg theorem, there is some $a \in \sG$ such that 
\begin{equation}\label{ord10}
  x=a\sigma(a)^{-1}.
\end{equation}
We also note that in $\Aut(S)$ 
\begin{equation}\label{ord11}
  \bigl(x\sigma^{-1}\bigr)^{2n+1} = x\cdot\sigma(x)\cdot \sigma^2(x) \cdot \ldots \cdot \sigma^{2n}(x).
\end{equation}  
As in the proof of \cite[Theorem 2.16]{GMPS}, define $t:=a^{-1}s^{(2n+1)}a=\bigl(x\sigma^{-1}\bigr)^{2n+1}a$. 
Then using \eqref{ord10} and 
\eqref{ord11} we obtain
$$\begin{aligned}\sigma(t) & =\sigma(a)^{-1} \cdot \bigl( \sigma(x)\cdot\sigma^2(x)\cdot \ldots \cdot \sigma^{2n}(x)\cdot\sigma^{2n+1}(x)
     \bigr) \sigma(a)\\ 
     & = \bigl(\sigma(a)^{-1}x^{-1}\bigr) \cdot \bigl( x\cdot\sigma(x) \cdot \ldots \cdot \sigma^{2n}(x) \bigr) \cdot \bigl(x\sigma(a))\\
     & = a^{-1}\bigl(x\sigma^{-1}\bigr)^{2n+1}a = t.\end{aligned}$$ 
Thus $t \in \sG^\sigma \cong \tw2 B_2(2)$, and note that $|\tw2 B_2(2)| = 20$. Since $s$ has odd order and $t,s^{2n+1}$ are conjugate
in $\sG$, $|t|$ is odd, whence $|s^{2n+1}| = |t|$ divides $5$. Thus $|s|$ divides $5(2n+1)$, as stated.
\end{proof}

\begin{thm}\label{main3}
Suppose the Airy sheaf $\sF_q$ for $q=2^{2n+1}$ has finite geometric monodromy group $G_\geo$. Then the following statements
hold.
\begin{enumerate}[\rm(i)]
\item The arithmetic monodromy group over $\F_2$ of $\sF_q$ is $G_{\ari,\F_2} = C \times \Aut(\tw2 B_2(q))$, with $|C| \leq 2$. 
\item Moreover, $C=1$ if $2n+1 \equiv \pm 3 \pmod*{8}$, and $C \cong C_2$ if $2n+1 \equiv \pm 1 \pmod*{8}$.
\item Suppose $2n+1 \equiv \pm 3 \pmod*{8}$. For the arithmetic monodromy group $G_{\ari,k}$ of $\sF_q$ over a finite extension 
$k/\F_2$, we have $G_{\ari,\F_2} = \Aut(\tw2 B_2(q))$, $G_{\ari,k} = G_\geo=\tw2 B_2(q)$ when $k \supseteq \F_q$, and 
$[G_{\ari,k}:G_\geo] = \deg(\F_q/k)$ when $k \subseteq \F_q$.
\item Suppose $2n+1 \equiv \pm 1 \pmod*{8}$. 
\begin{enumerate}
\item[\rm{($\alpha$)}] For the arithmetic monodromy group $G_{\ari,k}$ of $\sF_q$ over a finite extension 
$k/\F_2$, we have $G_{\ari,k} = G_\geo=\tw2 B_2(q)$ when $k \supseteq \F_{q^2}$, and 
$[G_{\ari,k}:G_{\geo}] = \deg(\F_{q^2}/k)$ when $k \subseteq \F_{q^2}$. 
\item[\rm{($\beta$)}] For the arithmetic monodromy group $G_{\ari,\tilde\sF_q,k}$ of the sheaf 
$\tilde\sF_q := \sF_q \otimes (-1)^{\deg/\F_2}$ over $k/\F_2$, we have $G_{\ari,\tilde\sF_q,\F_2} = \Aut(\tw2 B_2(q))$,  
$G_{\ari,\tilde\sF_q,k}=G_{\geo,\tilde\sF_q} = \tw2 B_2(q)$ when $k \supseteq \F_q$, and 
$[G_{\ari,\tilde\sF_q,k}:G_{\geo,\tilde\sF_q}] = \deg(\F_q/k)$ when $k \subseteq \F_q$.
\end{enumerate}
\end{enumerate}    
\end{thm}

\begin{proof}
Part (i) follows from Propositions \ref{traces-f0a}, \ref{traces-f0b}, and Theorem \ref{main2d}(iii). 

\smallskip
(ii) Let $g$ denote the image of $\Frob_{1,\F_2}$ in $G_{\ari,\F_2}$. As $\langle g,G_\geo \rangle = C \times \Aut(S)$ 
for $S := \tw2 B_2(q)$, we can write $g =zs$ with $C = \langle z \rangle \leq C_2$ and $s \in \Aut(S) \cong S \rtimes C_{2n+1}$. 
The central element $z$ acts on $\sF_q$ as the scalar $\xi$, where $\xi=1$ if $|C|=1$ and $\xi=-1$ if $|C|=2$. 
The finiteness of $G_\geo$ implies 
that $G_{\ari,\F_2}$ is finite, and so $\Trace(g^m)$ is a Gaussian integer for any $m \in \Z$. Now Proposition \ref{finaltraces-f1} implies that
\begin{equation}\label{ord20}
  \Trace(s^{2n+1})=\xi\Trace(g^{2n+1}) = -\xi,
\end{equation}  
and note that $t:=s^{2n+1} \in S$. 

Also note that since $\Aut(S) = \langle s,S \rangle$, $s$ has order $2n+1$ modulo $S$, whence $|s|=e(2n+1)$ for some
$e \in \Z_{\geq 1}$, whence $|t|=e$. Inspecting the character table of $S$ \cite{Bur} and using \eqref{ord20}, we see that 
$e=|t|$ is odd and greater than $1$; in fact $e$ divides $q-2q_0+1$ if $\xi=1$ and $e$ divides $q+2q_0+1$ if $\xi=-1$, where $q_0:=2^n$. 
It follows that $|s|=e(2n+1)$ is odd. By Lemma \ref{order}, $e>1$ divides $5$, so $e=5$; in particular, $|s| = 5(2n+1)$. 
 An easy computation shows 
that $5|(q-2q_0+1)$ precisely when $n \equiv 1,2 \pmod*{4}$ (equivalently, $2n+1 \equiv \pm 3 \pmod* 8$), and $5|(q+2q_0+1)$ precisely when $n \equiv 0,3 \pmod*{4}$ (equivalently, $2n+1 \equiv \pm 1 \pmod* 8$).
Hence the statement follows, and we have also proved that the image $g$ of $\Frob_{1,\F_2}$ has order
\begin{equation}\label{ord21}
  \left\{\begin{array}{ll} 10(2n+1), & \mbox{if }2n+1 \equiv \pm 1 \pmod*{8},\\
  5(2n+1), & \mbox{if }2n+1 \equiv \pm 3 \pmod*{8}.\end{array}\right.
\end{equation}    

\smallskip
Parts (iii) and (iv)$(\alpha)$  follow from (i), (ii), and the facts that $G_\geo = \tw2 B_2(q)$ and $G_{\ari,\F_2}/G_\geo$ is cyclic of order $|C|(2n+1)$. 

To prove (iv)$(\beta)$, note that when $k \supseteq \F_4$, the images of $\pi_1(\A^1/k)$ on $\sF_q$ and 
$\tilde\sF_q$ are the same. Hence $G_{\ari,\tilde\sF_q,k} = G_{\ari,k} = \tw2 B_2(q)$ whenever $k \supseteq \F_{q^2}$, whence
$G_{\geo,\tilde\sF_q} = \tw2 B_2(q) = S$. Now $G_{\ari,\tilde\sF_q,\F_2} = \langle \tilde{g},S \rangle$, where $\tilde{g}$ is the image of 
$\Frob_{1,\F_2}$ on $\tilde\sF_q$. By its definition, $\tilde{g} = -g$, where $g$ is the image of 
$\Frob_{1,\F_2}$ on $\sF_q$, and the proof of (ii) shows (recalling $2n+1 \equiv \pm 1 \pmod*{8}$) that $g =-s$ with 
$\Aut(S) = \langle s,S \rangle$. Hence $G_{\ari,\tilde\sF_q,\F_2} = \langle s,S \rangle = \Aut(S) = S \rtimes C_{2n+1}$, and the assertion
follows. 
\end{proof}

\begin{rmk}
Computations in {\sc Magma} suggest that, for the Airy sheaf $\sF_q$ with $q=2^{2n+1}$, the Frobenii $\Frob_{a,\F_{2^j}}$ with
$a \in \F_{2^j}$ and $\gcd(j,2n+1) = 1$, all have traces of absolute value $1$. 
If one knew that $\sF_q$ has geometric monodromy group $G_{\geo,\sF_q}=\tw2 B_2(q)$, then this ``absolute value one" property agrees with Corollary \ref{ext-suz2}. Also, for $n = 1,2$, computations show that $\Frob_{1,\F_2}$ has order $5(2n+1)$, and this again 
agrees with \eqref{ord21}. 


On the other hand, the infinite case of the dichotomy, namely $G_\geo=\SL_D$ would imply by Deligne's equidistribution theorem  
\cite[Theorem 9.7.13]{Ka-Sar} that when $j$ is large enough (compared to $q$), some (in fact most) Frobenii $\Frob_{a,\F_{2^j}}$ would
have traces of absolute value $\neq 1$. This again gives some evidence in support of the geometric part of \cite[Conjecture 2.2]{Ka-ERS}
asserting that $G_{\geo,\sF_q}=\tw2 B_2(q)$. However, Theorem \ref{main3}(iii) shows that the arithmetic part of 
\cite[Conjecture 2.2]{Ka-ERS} stating that $G_{\ari,\sF_q,\F_2} = \Aut(\tw2 B_2(q))$ is false when $q=2^m$ with 
$m \equiv \pm 1 \pmod*{8}$; it should be corrected in that case by replacing $\sF_q$ by $\tilde\sF_q$ as in (iv)($\beta$) of Theorem \ref{main3}.
\end{rmk}

\subsection*{9C. Local systems with infinite monodromy groups} 
Theorem \ref{main2d} allows us to prove the following criterion for infinite monodromy.
 
\begin{prop}\label{inf-g1}
Let $n \in \Z_{\geq 1}$, $q=2^{2n+1}$, and consider the Airy sheaf $\sF(q,f)$ defined in \eqref{airy-s}
with $q=2^{2n+1}$, $f(x) \in \F_2[x]$, and $f(0)=0$. Suppose that 
\begin{enumerate}[\rm(i)]
\item $2n+1$ is coprime to $|\tw2 B_2(q)|$, and 
\item there is some odd integer $m$ coprime to $2n+1$ such that 
$|\Trace\bigl(\Frob_{0,\F_{2^m}}|\sF(q,f)\bigr)| \neq 1$. 
\end{enumerate}
Then $\sF(q,f)$ has infinite geometric monodromy group. 
\end{prop} 

\begin{proof}
Assume the contrary: $\sF(q,f)$ has finite geometric monodromy group $S$. By Theorem \ref{main2d}, 
$S=\tw2 B_2(q)$ and the arithmetic monodromy group of $\sF(q,f)$ over $\F_2$ is $G = \ZB(G) \times \Aut(S)$,
where $\ZB(G) \leq C_4$. In particular, $|G/S|$ divides $4(2n+1)$, and so $S$ equals the arithmetic monodromy group
of $\sF(q,f)$ over $\F_{q^4}$. Moreover, $G = \langle g,S \rangle$ for the image $g$ of $\Frob_{0,\F_2}$ in $G$, and 
we can write $g=zs$ with $z \in \ZB(G)$, $s \in \Aut(S)$, and $\Aut(S) = \langle s,S \rangle$. Applying 
Corollary \ref{ext-suz2}, we see that $|\Trace(g)| = |\Trace(s)|=1$, a contradiction.
\end{proof}
 
\begin{prop}\label{traces-f0c}
Let $n \in \Z_{\geq 2}$, $q:=2^{2n+1}$, 
$r:=[(n-1)/2]$, and define 
$$f_1(x):= \sum^r_{i=0}x^{1+2^{n-2i}}.$$
Consider the sheaf $\sF(q,f)$ with $f(x):=f_1(x^{t(q)})$ as in \eqref{airy-s}. 
Then for $m:= 2[n/2]+1$ we have
$|\Trace\bigl(\Frob_{0,\F_{2^m}}|\sF(q,f)\bigr)| \neq 1$.
\end{prop}

\begin{proof}
As in the proof of Proposition \ref{traces-f0a}, the starting point is that $\gcd(t(q),2^m-1)=1$. To see this, note that $t(q)=q+1-2q_0$ divides $q^4-1$, and $\gcd(q^4-1,2^m-1) =2^{\gcd(4(2n+1),m)}-1=1$. Thus for the field $k:=\F_{2^m}$, 
the map $x \mapsto x^{t(q)}$ is bijective on $k$. 

The sheaf $\sF(q,f)$ was built out of the Witt vector $\bigl[x^{t(q)},f_1(x^{t(q)}\bigr]$.
Let us denote by $\sH(q,f_1)$ the sheaf built by the same recipe, with same clearing factor, out of the Witt vector
$[x,f_1(x)]$. Then we have
$$\Trace(\Frob_{0,k} |\sF(q,f)) =\Trace(\Frob_{0,k}|\sH(q,f_1)),$$
precisely because the map $x \mapsto x^{t(q)}$ is bijective on $k$. 

Next, as in \eqref{fq-hq2}, we write the input vector $[x,f_1(x)]$ as $[x,xR(x)]$, with 
$$R(x):=\sum^r_{i=0}x^{2^{n-2i}} = F^n(x)+F^{n-2}(x)+ \ldots + F^{n-2r}(x),$$
and $F$ denotes the absolute Frobenius. Now we can repeat the arguments in the proof of Proposition \ref{traces-f0a}, and 
compute the form $\Trace_{k/\F_2}(\langle x,y \rangle)$ on $k \times k$ with $\langle x,y \rangle := xy+xR(y)+yR(x)$ as in \eqref{fq-hq4}.
For any $x,y \in k$ we have $F^m(x)=x$, $F^m(y)=y$. If $2|n$, then $n=2r+2$, $m=2r+3$, and  
$$\langle x,y \rangle =xy +y\bigl(\sum^{r+1}_{i=1}F^{2i}(x)\bigr) +\sum^{r+1}_{i=1}xF^{2i}(y)$$ 
has the same trace over $\F_2$ as 
$$\begin{aligned}xy +y\bigl(\sum^{r+1}_{i=1}F^{2i}(x)\bigr) + \sum^{r+1}_{i=1}F^{2r+3-2i}\bigl(xF^{2i}(y)\bigr) & = xy +y\bigl(\sum^{r+1}_{i=1}F^{2i}(x)\bigr) + \sum^{r+1}_{i=1}yF^{2r+3-2i}(x)\\
  & =y\bigl( \sum^{2r+2}_{j=0}F^j(x)\bigr)=y\Trace_{k/\F_2}(x).\end{aligned}$$
%
If $2 \nmid n$, then $n=m=2r+1$, then 
$$\langle x,y \rangle = xy +y\bigl(\sum^{r}_{i=0}F^{2i+1}(x)\bigr) +\sum^r_{i=0}xF^{2i+1}(y) = 
   3xy+ y\bigl(\sum^{r-1}_{i=0}F^{2i+1}(x)\bigr) +\sum^{r-1}_{i=0}xF^{2i+1}(y)$$ 
has the same trace over $\F_2$ as of
$$\begin{aligned}xy +y\bigl(\sum^{r-1}_{i=0}F^{2i+1}(x)\bigr) +\sum^{r-1}_{i=0}F^{2r-2i}\bigl(xF^{2i+1}(y) \bigr) & 
   =xy +y\bigl(\sum^{r-1}_{i=0}F^{2i+1}(x)\bigr) +\sum^{r-1}_{i=0}yF^{2r-2i}(x)\\
   & = y\bigl(\sum^{2r}_{j=0}F^j(x)\bigr) 
     =  y\Trace_{k/\F_2}(x).\end{aligned}$$
So in both cases, the symmetric bilinear form $\Trace_{k/\F_2}(\langle x,y \rangle)$ on $k \times k$ has kernel consisting 
of the elements $x \in k$ with $\Trace_{k/\F_2}(x)=0$, that is, of exactly $2^{m-1}$ elements.
Then the proof of Proposition 
\ref{traces-f0a} shows that $|\Trace\bigl(\Frob_{0,\F_{2^m}}|\sF(q,f)\bigr)|$ is either $0$ or $2^{(m-1)/2} >1$, and hence can never be 
equal to $1$ (since by hypothesis $n \geq 2$, and hence $m \ge 3$).
\end{proof} 

\begin{thm}\label{inf-g}
Let $n \in \Z_{\geq 2}$, $q:=2^{2n+1}$, and consider the sheaf $\sF(q,f)$ of rank $D:=2^n(2^{2n+1}-1)$, 
with $f(x)=f_1(x^{t(q)})$ and $f_1(x)$ as defined in
Proposition \ref{traces-f0c}. Assume in addition that $2n+1$ is coprime to $|\tw2 B_2(q)|$; for instance, take $2n+1=\ell^a$
for any odd prime $\ell \neq 5$ and any
$a \in \Z_{\geq 1}$.
Then the geometric monodromy group of $\sF(q,f)$ is $\SL_D$.
\end{thm}
 
\begin{proof}
We first apply Proposition \ref{traces-f0a}, and note that $m$ is coprime to $2(2n+1)$. It then follows from Proposition 
\ref{inf-g1} that $\sF(q,f)$ has infinite geometric monodromy group $G_\geo$. By Theorem \ref{main2}, $G_\geo=\SL_D$.

Suppose that $2n+1=\ell^a$ for a prime $\ell$, but $\ell$ divides $|\tw2 B_2(q)|$. Then $\ell$ divides $q^4-1 = 2^{4\ell^a}-1$. 
Since $\ell|(2^{\ell-1}-1)$, $\ell$ divides $\gcd(2^{4\ell^a}-1,2^{\ell-1}-1) = 2^{\gcd(4\ell^a,\ell-1)}-1$, and so $\ell$ divides $2^4-1=15$. 
But $3 \nmid |\tw2 B_2(q)|$, so $r=5$.
\end{proof} 

More generally, to ensure that $2n+1$ is coprime to $|\tw2 B_2(q)|$ for $q=2^{2n+1}$, we can take any $n$ such that
$2n+1 = p_1^{a_1}p_2^{a_2} \ldots p_t^{a_t}$, where $p_1 < p_2 < \ldots < p_t$ are primes, $p_i \neq 5$, $a_i \in \Z_{\geq 1}$, and 
$p_i \nmid (p_j-1)$ whenever $i < j$. 
Indeed, suppose $p_j$ divides $|\tw2 B_2(q)|$ for some $j$. Then $p_j$ divides both $2^{p_j-1}-1$ and $q^4-1 = 2^{4(2n+1)}-1$. 
Since $\gcd(p_j-1,4(2n+1)) = \gcd(p_j-1,4\prod^t_{i=1}p_i^{a_i})$ divides $4$, it follows that $p_j$ divides $2^4-1=15$.
But $p_j \neq 5$ by assumption, and $p_j \neq 3$ since $3 \nmid |\tw2 B_2(q)|$, a contradiction.

\end{document}